
\newcount\chapno
\newcount\parno
\newcount\secno
\newcount\subsecno
\newdimen\partisize


\magnification 1200

\hsize=156 true mm
\vsize=8.9 true in

\tolerance=300
\pretolerance=100

\parindent=15pt

\mathsurround=0pt

\normallineskiplimit=.5pt
\normalbaselineskip=15pt
\normallineskip=1pt plus .5 pt minus .5 pt
\normalbaselines

\abovedisplayskip = 8pt plus 3pt minus 3pt
\abovedisplayshortskip = 1pt plus 2pt
\belowdisplayskip = 8pt plus 3pt minus 3pt
\belowdisplayshortskip = 5pt plus 2pt

\partisize=14pt

\newdimen\firstident
\newdimen\secondident
\newdimen\thirdident
\newbox\identbox
\newcount\assertioncount


\font\fractur=cmfrak


\newskip\chapskipamount              
 \chapskipamount=30pt plus6pt minus4pt
 \def\chapskip{\vskip\chapskipamount}
\newskip\parskipamount               
 \parskipamount=21pt plus4pt minus4pt
 \def\parskip{\vskip\parskipamount}
\newskip\secskipamount               
 \secskipamount=18pt plus4pt minus2pt
 \def\secskip{\vskip\secskipamount}


\font\titlefont=cmr12 at 18pt
\font\chaptifont=cmbx10 at 15pt
\font\partifont=cmcsc12
\font\namefont=cmr12


\outer\def\chapter#1{
   \global\parno=0
   \global\advance\chapno by 1
   \chapskip
   \centerline{{\chaptifont \the\chapno~}
   {\chaptifont #1}}
   \par\nobreak}


\outer\def\endchapter{
  \newpage}


\outer\def\paragraph#1{\goodbreak
   \global\secno=0
   \global\advance\parno by 1
   \parskip
   \hangitem{\partifont \the\chapno.\the\parno~}
   {\partifont #1}
   \par\nobreak}

%

\outer\def\Notation#1: #2\par{\nobreak\secskip
  \itemitem{}{\bf Notation#1}:\enspace{\sl#2\par}\penalty50}

%

\outer\def\secstart#1{\penalty-100
   \global\subsecno=0
   \global\advance\secno by 1
   \secskip\noindent
   {\bf (\the\chapno.\the\parno.\the\secno)~#1}}

%

\outer\def\proof#1: {\smallbreak {\it Proof\ #1}\/:\enspace}

%

\outer\def\claim#1: #2\par{\smallbreak
   {\bf #1}:\enspace{\sl#2\par}\penalty50}


\outer\def\bibliography{
   \newpage
   \tolerance=1000
   \parindent=0pt
   {\partifont Bibliography}
   \parskip
   \frenchspacing
   }




\let\pz=\S

%

\let\eps=\varepsilon

%

\def\Abf{{\bf A}}

%

\def\Fbar{{\bar F}}

\def\hbar{{\bar h}}

\def\kbar{{\bar k}}
\def\Kbar{{\bar K}}

\def\Lbar{{\bar L}}

%

\def\fhat{{\hat f}}

\def\Hhat{{\hat H}}

\def\Khat{{\hat K}}

\def\That{{\hat T}}

%

%

\def\hline{\underline{h}}

\def\nline{\underline{n}}

\def\xline{\underline{x}}

%

\def\Wtilde{{\tilde W}}

\def\Xtilde{{\tilde X}}

%

%

%

%

%

\def\Qdbar{{\bar {\Q}}}

%

%

%

%

\def\hgbar{{\bar \eta}}

\def\kgbar{{\bar \kappa}}

\def\ygbar{{\bar \psi}}

%

%

\def\Sgtilde{{\tilde \Sigma}}

%

%

%

%

%


\def\star{{}^*}

\def\vdual{{}^{\mathord{\vee}}}               
\def\cross{{}^{\times}}
\def\dlbrack{\lbrack{\mskip-2mu}\lbrack}      
\def\drbrack{\rbrack{\mskip-2mu}\rbrack}      

\def\Qp{\Q_{p}}
\def\Fp{{\F_{p}}}

\def\Zp{\Z_{p}}


\def\overneq#1#2{\lower0.5pt\vbox{\lineskiplimit\maxdimen\lineskip-.5pt
                \ialign{$#1\hfil##\hfil$\crcr#2\crcr\not=\crcr}}}


\def\stimes{\mathbin{\raise1pt\hbox{$\scriptscriptstyle \bf
            \vert$}\mkern-5mu\times}}   


\def\ad{\mathop{\rm ad}\nolimits}

\def\Aut{\mathord{\rm Aut}}

\def\det{\mathop{\rm det}\nolimits}
\def\diag{\mathop{\rm diag}\nolimits}

\def\dim{\mathop{\rm dim}\nolimits}
\def\End{\mathop{\rm End}\nolimits}

\def\et{{\rm \acute et}}

\def\exp{\mathop{\rm exp}\nolimits}

\def\Gal{\mathop{\rm Gal}\nolimits}

\def\Hom{\mathop{\rm Hom}\nolimits}
\def\Homline{\mathop{\underline {\Hom}}\nolimits}

\def\id{\mathop{\rm id}\nolimits}

\def\inf{\mathop{\rm inf}}
\let\integral=\int
\def\int{\mathop{\rm int}\nolimits}

\def\inv{\mathop{\rm inv}\nolimits}

\def\Ker{\mathop{\rm Ker}\nolimits}
\def\length{\mathop{\rm length}\nolimits}
\def\Lie{\mathop{\rm Lie}\nolimits}
\def\log{\mathop{\rm log}\nolimits}

\def\Re{\mathop{\rm Re}\nolimits}

\def\Spec{\mathop{\rm Spec}\nolimits}
\def\Spf{\mathop{\rm Spf}\nolimits}

\def\vol{\mathop{\rm vol}\nolimits}

\def\setback(#1){\mathrel{\mkern-#1mu}}
\def\varfill#1{$\smash{#1} \mkern-8mu
   \cleaders\hbox{$\mkern-3mu \smash{#1} \mkern-3mu$}\hfill
   \mkern-8mu #1$}
\def\equalfill{\varfill{=}}


\let\ar=\rightarrow
\let\al=\leftarrow

\let\air=\hookrightarrow

\let\asr=\mapsto

\let\arr=\longrightarrow
\let\all=\longleftarrow
\def\aerr{\arr\setback(27)\arr}        

\def\arriso{\buildrel\sim\over\arr}    
\def\arrover#1{\buildrel#1\over\arr}   





\def\allover#1{\buildrel#1\over\all}   

\def\arvar(#1){\hbox to #1pt{\rightarrowfill}}
\def\alvar(#1){\hbox to #1pt{\leftarrowfill}}
\def\arvarover(#1)#2{\mathop{\arvar(#1)}\limits^{#2}}
\def\alvarover(#1)#2{\mathop{\alvar(#1)}\limits^{#2}}


\let\ad=\downarrow
\let\au=\uparrow
\def\add{\Big\ad}                      
\def\auu{\Big\au}                      
\def\addleft#1{\llap{$\vcenter{\hbox{$\scriptstyle #1$}}$}\add}
\def\auuleft#1{\llap{$\vcenter{\hbox{$\scriptstyle #1$}}$}\auu}
\def\addright#1{\add\rlap{$\vcenter{\hbox{$\scriptstyle #1$}}$}}







\def\aqrvar(#1){\hbox to #1pt{\equalfill}}


\let\implies=\Rightarrow

\let\afr=\implies

\let\iff=\Leftrightarrow






%
%

\def\set#1#2{\{\,#1\>\vert\>\hbox{#2}\,\}}

%
%

\def\bigset#1#2#3#4{\eqalign{#1\ \bigl\{\,#2\ \bigm| \>
    & \hbox{#3} \cr
    & \hbox{#4} \,\bigr\} \cr}}

%
%

\def\powerseries over #1 in #2{{#1 \dlbrack #2 \drbrack}}

%
%

\def\laurentseries over #1 in #2{{#1 (\!( #2 )\!)}}

%
%

\def\smallmatrix(#1,#2;#3,#4){\left({{#1\atop #3}\>{#2\atop #4}}\right)}

%
%

\def\twomatrix#1#2#3#4{\left({{#1\atop #3}\>{#2\atop #4}}\right)}

%
%

\def\tensor#1{\otimes_{#1}}

%
%

%
%

%

%

\def\restricted#1{\vert_{#1}}

%
%

\def\subscript#1\atop#2{_{\scriptstyle #1 \atop #2}}

%
%

\def\limind{\mathop{\mathop{\lim}\limits_{\arr}}}

%
%

\def\limproj{\mathop{\mathop{\lim}\limits_{\all}}}




\newfam\frac
\textfont\frac=\fractur
\def\afr{{\fam\frac a}}

\def\mfr{{\fam\frac m}}

\def\pfr{{\fam\frac p}}

\def\Pfr{{\fam\frac P}}

%
%
%
%
\font\tenss=cmss10
\newfam\ssfam %
\textfont\ssfam=\tenss
\catcode`\_=11
\def\suf_fix{}
\def\scaled_rm_box#1{%
 \relax
 \ifmmode
   \mathchoice
    {\hbox{\tenrm #1}}%
    {\hbox{\tenrm #1}}%
    {\hbox{\sevenrm #1}}%
    {\hbox{\fiverm #1}}%
 \else
  \hbox{\tenrm #1}%
 \fi}
\def\suf_fix_def#1#2{\expandafter\def\csname#1\suf_fix\endcsname{#2}}
\def\I_Buchstabe#1#2#3{%
 \suf_fix_def{#1}{\scaled_rm_box{I\hskip-0.#2#3em #1}}
}
\def\rule_Buchstabe#1#2#3#4{%
 \suf_fix_def{#1}{%
  \scaled_rm_box{%
   \hbox{%
    #1%
    \hskip-0.#2em%
    \lower-0.#3ex\hbox{\vrule height1.#4ex width0.07em }%
   }%
   \hskip0.50em%
  }%
 }%
}
\I_Buchstabe B22
\rule_Buchstabe C51{34}
\I_Buchstabe D22
\I_Buchstabe E22
\I_Buchstabe F22
\rule_Buchstabe G{525}{081}4
\I_Buchstabe H22
\I_Buchstabe I20
\I_Buchstabe K22
\I_Buchstabe L20
\I_Buchstabe M{20em }{I\hskip-0.35}
\I_Buchstabe N{20em }{I\hskip-0.35}
\rule_Buchstabe O{525}{095}{45}
\I_Buchstabe P20
\rule_Buchstabe Q{525}{097}{47}
\I_Buchstabe R21 
\rule_Buchstabe U{45}{02}{54}
\suf_fix_def{Z}{\scaled_rm_box{Z\hskip-0.38em Z}}
\catcode`\"=12
\newcount\math_char_code
\def\suf_fix_math_chars_def#1{%
 \ifcat#1A
  \expandafter\math_char_code\expandafter=\suf_fix_fam
  \multiply\math_char_code by 256
  \advance\math_char_code by `#1
  \expandafter\mathchardef\csname#1\suf_fix\endcsname=\math_char_code
  \let\next=\suf_fix_math_chars_def
 \else
  \let\next=\relax
 \fi
 \next}
%
%
%
%
\def\font_fam_suf_fix#1#2 #3 {%
 \def\suf_fix{#2}
 \def\suf_fix_fam{#1}
 \suf_fix_math_chars_def #3.
}
\font_fam_suf_fix
 0rm
 ABCDEFGHIJKLMNOPQRSTUVWXYZabcdefghijklmnopqrstuvwxyz
\font_fam_suf_fix
 2scr
 ABCDEFGHIJKLMNOPQRSTUVWXYZ
\font_fam_suf_fix
 \slfam sl
 ABCDEFGHIJKLMNOPQRSTUVWXYZabcdefghijklmnopqrstuvwxyz
\font_fam_suf_fix
 \bffam bf
 ABCDEFGHIJKLMNOPQRSTUVWXYZabcdefghijklmnopqrstuvwxyz
\font_fam_suf_fix
 \ttfam tt
 ABCDEFGHIJKLMNOPQRSTUVWXYZabcdefghijklmnopqrstuvwxyz
\font_fam_suf_fix
 \ssfam
 ss
 ABCDEFGHIJKLMNOPQRSTUVWXYZabcdefgijklmnopqrstuwxyz
\catcode`\_=8
%
%
%
%
%
%
%
%
%
%
%
%
%
%
%
%
%
%
%
%
%
%
%
%
%
%



\def\en_item#1#2{%
 \par
 \setbox\identbox=\hbox #1{#2}%
 \noindent
 \hangafter=1%
 \hangindent=\wd\identbox
 \box\identbox
 \ignorespaces
}

\def\ennopar_item#1#2{%
 \setbox\identbox=\hbox #1{#2}%
 \noindent
 \hangafter=1%
 \hangindent=\wd\identbox
 \box\identbox
 \ignorespaces
}

%
%

\def\hangitem#1{\en_item{}{#1\enspace}}

%
%
%

%

\def\hanghangitem#1{\en_item{}{\kern\firstident #1\enspace}}

%
%
\def\hanghanghangitem#1{\en_item{}{\kern\secondident #1\enspace}}


\def\newpage{\vfil\eject}

%
%
\def\indention#1{%
 \setbox\identbox =\hbox{\kern\parindent{#1}\enspace}%
 \firstident=\wd\identbox
}

%
%
\def\subindention#1{%
 \setbox\identbox=\hbox{{#1}\enspace}%
 \secondident=\wd\identbox
 \advance \secondident by \firstident
}

%
%
\def\subsubindention#1{%
 \setbox\identbox=\hbox{{#1\ }\enspace}%
 \thirdident=\wd\identbox
 \advance \thirdident by \secondident
}%

%
%
\def\litem#1{\en_item{to\firstident}{\kern\parindent#1\hfil \enspace }}
\def\llitem#1{\en_item{to\secondident}{\kern\firstident#1\hfil\enspace}}
\def\lllitem#1{\en_item{to\thirdident}{\kern\secondident#1\hfil\enspace}}
\def\ritem#1{\en_item{to\firstident}{\kern\parindent\hfil#1\enspace}}
\def\rritem#1{\en_item{to\secondident}{\kern\firstident\hfil#1\enspace}}
\def\rrritem#1{\en_item{to\thirdident}{\kern\secondident\hfil#1\enspace}}
\def\citem#1{\en_item{to\firstident}{\kern\parindent\hfil#1\hfil\enspace}}
\def\ccitem#1{\en_item{to\secondident}{\kern\firstident\hfil#1\hfil\enspace}}
\def\cccitem#1{\en_item{to\thirdident}{\kern\secondident\hfil#1\hfil\enspace}}
\def\rmlitem#1{\en_item{to\firstident}{\kern\parindent{\rm #1}\hfil \enspace }}
\def\rmllitem#1{\en_item{to\secondident}{\kern\firstident{\rm #1}\hfil\enspace}}
\def\rmlllitem#1{\en_item{to\thirdident}{\kern\secondident{\rm #1}\hfil\enspace}}

%
%

\def\lnoitem#1{\ennopar_item{to\firstident}{\kern\parindent#1\hfil \enspace }}
\def\llnoitem#1{\ennopar_item{to\secondident}{\kern\firstident#1\hfil\enspace}}
\def\lllnoitem#1{\ennopar_item{to\thirdident}{\kern\secondident#1\hfil\enspace}}

%
%

\def\assertionlist{
   \assertioncount=0
   \indention{(2)}}

\def\assertionitem{
   \advance\assertioncount by 1
   \litem{{\rm(\number\assertioncount)}}}

%
%

\def\bulletlist{
   \indention{$\bullet$}}

\def\bulletitem{
   \litem{$\bullet$}}

%
%

%


\global\chapno=0
\global\parno=0
\global\secno=0
\global\subsecno=0





%

\def\lformno{
   \global\advance\subsecno by 1
   \leqno(\the\chapno.\the\parno.\the\secno.\the\subsecno)
   }


\def\displayno{
   \global\advance\subsecno by 1
   (\the\chapno.\the\parno.\the\secno.\the\subsecno)
   }


\def\actualsecno{(\the\chapno.\the\parno.\the\secno)}


\def\actualsubsecno{(\the\chapno.\the\parno.\the\secno.\the\subsecno)}

%
%
%

\def\label#1{\xdef#1{\actualsecno}}

%
%
%

\def\sublabel#1{\xdef#1{\actualsubsecno}}


\def\rec{\mathord{\rm rec}}
\outer\def\reasoning#1: {\smallbreak {\it Reasoning\ #1}\/:\enspace}
\outer\def\hint#1: {\smallbreak {\it Hint\ #1}\/:\enspace}
\def\JL{\mathord{\rm JL}}
\def\Art{\mathord{\rm Art}}
\def\GL{\mathord{\rm GL}}
\def\SL{\mathord{\rm SL}}
\def\Rep{\mathord{\rm Rep}}


\ 
\vskip 3 true cm
\centerline{\titlefont The local Langlands correspondence}
\bigskip
\centerline{\titlefont for GL(n) over p-adic fields}
\parskip
\centerline{\namefont Torsten Wedhorn\footnote*{wedhorn@mi.uni-koeln.de}}
\vskip 2 true cm
\centerline{\it Mathematisches Institut der Universit\"at zu K\"oln}
\centerline{\it Weyertal 86-90, D-50931 K\"oln, Germany}
\vskip 9 true cm
\centerline{\it Lecture given at the}
\centerline{\it School on Automorphic forms on $\GL(n)$}
\centerline{\it Trieste, 31 July -- 11 August, 2000}

\newpage

\centerline{\bf Abstract}

\secskip

This work is intended as an introduction to the statement and the
construction of the local Langlands correspondence for $\GL(n)$ over
$p$-adic fields. The emphasis lies on the statement and the explanation of
the correspondence.

\secskip

{\it Keywords}: Langlands Program, Representations of $p$-adic Groups,
Weil-Deligne Representations, Formal Groups, Vanishing Cycles.

{\it AMS classifications}: 11S37, 11F70, 11F80, 11S40, 14L05, 20G25, 22E50

\newpage

{\partifont Contents}

\parskip

\indention{2.}
Introduction\par
Notations
\litem{1.} The local Langlands correspondence
\subindention{2.2}
\llitem{1.1} The local Langlands correspondence for GL(1)
\llitem{1.2} Formulation of the local Langlands correspondence
\litem{2.} Explanation of the GL(n)-side
\llitem{2.1} Generalities on admissible representations
\llitem{2.2} Induction and the Bernstein-Zelevinsky classification for
GL(n)
\llitem{2.3} Square integrable and tempered representations
\llitem{2.4} Generic representations
\llitem{2.5} Definition of $L$- and epsilon-factors
\litem{3.} Explanation of the Galois side
\llitem{3.1} Weil-Deligne representations
\llitem{3.2} Definition of $L$- and epsilon-factors
\litem{4.} Construction of the correspondence
\llitem{4.1} The correspondence in the unramified case
\llitem{4.2} Some reductions
\llitem{4.3} A rudimentary dictionary of the correspondence
\llitem{4.4} The construction of the correspondence after Harris and Taylor
\litem{5.} Explanation of the correspondence
\llitem{5.1} Jacquet-Langlands theory
\llitem{5.2} Special $p$-divisible $O$-modules
\llitem{5.3} Deformation of $p$-divisible $O$-modules
\llitem{5.4} Vanishing cycles
\llitem{5.5} Vanishing cycles on the universal deformation of special
$p$-divisible $O$-modules\par
Bibliography

\newpage

\centerline{\chaptifont Introduction}
\parskip

Let $K$ be a local field, i.e.\ $K$ is either the field of real or complex
numbers (in which case we call $K$ archimedean) or it is a finite extension
of $\Q_p$ (in which case we call $K$ $p$-adic) or it is isomorphic to
$\F_q(\!(t)\!)$ for a finite field $\F_q$ (in which case we call $K$ a local
function field). The local Langlands conjecture for $GL_n$ gives a
bijection of the set of equivalence classes of admissible representations
of $\GL_n(K)$ with the set of equivalence classes of $n$-dimensional
Frobenius semisimple representations of the
Weil-Deligne group of $K$. This bijection should be compatible with $L$- and
$\eps$-factors. For the precise definitions see chap.\ 2 and chap.\ 3.

If $K$ is archimedean, the local Langlands conjecture is known for a long
time and follows from the classification of (infinitesimal) equivalence
classes of admissible representations of $\GL_n(K)$ (for $K = \C$ this is
due to \v Zelobenko and Na\u\i mark and for $K = \R$ this was done by
Langlands). The archimedean case is particularly simple because all
representations of $\GL_n(K)$ can be built up from representations of
$\GL_1(\R)$, $\GL_2(\R)$ and $\GL_1(\C)$. See the survey article of Knapp
[Kn] for more details about the local Langlands conjecture in the
archimedean case.

If $K$ is non-archimedean and $n = 1$, the local Langlands conjecture is
equivalent to local abelian class field theory and hence is known for a
long time (due originally to Hasse [Has]). Of course, class field theory
predates the general Langlands conjecture. For
$n = 2$ the local (and even the global Langlands conjecture) are also known
for a couple of years (in the function field case this is due to Drinfeld
[Dr1][Dr2], and in the $p$-adic case due to Kutzko [Kut] and Tunnel
[Tu]). Later on Henniart [He1] gave also a proof for the $p$-adic case for
$n = 3$.

If $K$ is a local function field, the local Langlands conjecture for arbitrary
$n$ has been proved by Laumon, Rapoport and Stuhler [LRS] generalizing
Drinfeld`s methods. They use certain moduli spaces of ``$\Dscr$-elliptic
sheaves'' or ``shtukas'' associated to a global function field.

Finally, if $K$ is a p-adic field, the local Langlands conjecture for all
$n$ has been proved by Harris and Taylor [HT] using Shimura varieties,
i.e.\ certain moduli spaces of abelian varieties. A few months later
Henniart gave a much simpler and more elegant proof [He4]. On the other
hand, the advantage of the methods of Harris and Taylor is the geometric
construction of the local Langlands correspondence and that it establishes
many instances of compatibility between the global and the local
correspondence.

Hence in all cases the local Langlands conjecture for $\GL(n)$ is now a
theorem! We remark that in all cases the proof of the local Langlands
conjecture for $n > 1$ uses global methods although it is a purely local
statement. In fact, even for $n = 1$ (i.e.\ the case of local class field
theory) the first proof was global in nature.

\medskip

This work is meant as an introduction to the local Langlands
correspondence in the $p$-adic case. In fact, approximately half of it
explains the precise statement of the local Langlands conjecture as
formulated by Henniart. The other half gives the construction of the
correspondence by Harris and Taylor. I did not make any attempt to explain
the connections between the local theory and the global theory of
automorphic forms. In particular, nothing is said about the proof that the
constructed map satisfies all the conditions postulated by the local
Langlands correspondence, and this is surely a severe shortcoming. Hence
let me at least here briefly sketch the idea roughly:

Let $F$ be a number field which is a totally imaginary extension of a
totally real field such that there exists a place $w$ in $F$ with $F_w =
K$. The main idea is to look at the cohomology of a certain projective
system $X = (X_m)_m$ of projective $(n-1)$-dimensional
$F$-schemes (the $X_m$ are ``Shimura varieties of PEL-type'', i.e.\ certain
moduli spaces of abelian varieties with polarizations and a level structure
depending on $m$). This system is associated to a
reductive group $G$ over $\Q$ such that $G \tensor{\Q} \Qp$ is
equal to
$$\Q\cross_p \times GL_n(K) \times \hbox{anisotropic mod center factors}.$$
More precisely, these anisotropic factors are algebraic groups associated to
skew fields. They affect the local structure of the $X_m$ only in a minor
way, so let us ignore them for the rest of this overview. By the general
theory of Shimura varieties, to every absolutely irreducible representation
$\xi$ of $G$ over $\Q$ there is associated a smooth $\Qdbar_{\ell}$-sheaf
$\Lscr_{\xi}$ on $X$ where $\ell \not= p$ is some fixed prime. The cohomology
$H^i(X,\Lscr_{\xi})$ is an infinite-dimensional $\Qdbar_{\ell}$-vector
space with an action of $G({\bf A}_f) \times \Gal(\Fbar/F)$ where ${\bf
  A}_f$ denotes the ring of finite adeles of $\Q$. We can choose $\xi$ in
such a way that $H^i(X,\Lscr_{\xi}) = 0$ for $i \not= n-1$.

We have a map from the set of equivalence classes of irreducible admissible
representations $\Pi$ of $G({\bf A}_f)$ to the set of finite-dimensional
representations of $W_K \subset \Gal(\Fbar/F)$ where $W_K$ denotes the Weil
group of $K$ by sending $\Pi$ to
$$R_{\xi}(\Pi) = \Hom_{G(\Abf_f)}(\Pi,H^{n-1}(X,\Lscr_{\xi})).$$
For every such $\Pi$ the decomposition
$$G(\Abf_f) = \Q\cross_p \times GL_n(K) \times \hbox{remaining components}$$
gives a decomposition
$$\Pi = \Pi_0 \otimes \Pi_w \otimes \Pi^w.$$
If $\pi$ is a supercuspidal representation of $GL_n(K)$ then we can find a
$\Pi$ as above such that $\Pi_w \cong \pi\chi$ where $\chi$ is an
unramified character of $K\cross$, such that $\Pi_0$ is unramified and such
that $R_{\xi}(\Pi) \not= 0$.

Now we can choose a model $\Xtilde_m$ of $X_m \tensor{F} K$ over $O_K$ and
consider the completions $R_{n,m}$ of the local rings of a certain stratum
of the special fibre of $\Xtilde_m$. These completions carry canonical
sheaves $\psi^i_m$ (namely the sheaf of vanishing cycles) and their limits
$\psi^i$ are endowed with a canonical action $\GL_n(K) \times D\cross_{1/n}
\times W_K$ where $D\cross_{1/n}$ is the skew field with invariant $1/n$
and center $K$. If $\rho$ is any irreducible representation of $D\cross_{1/n}$,
$\psi^i(\rho) = \Hom(\rho,\psi^i)$ is a representation of $\GL_n(K) \times
W_K$. Via Jacquet-Langlands theory we can associate to every supercuspidal
representation $\pi$ of $\GL_n(K)$ an irreducible representation $\rho =
jl(\pi\vdual)$ of $D\cross_{1/n}$. Now there exists an $n$-dimensional
representation $r(\pi)$ of $W_K$ which satisfies
$$[\pi \otimes r(\pi)] =
\sum_{i=0}^{n-1}(-1)^{n-1-i}[\psi^i(jl(\pi\vdual))]$$
and
$$n\cdot[R_{\xi}(\Pi) \otimes \chi(\Pi_0 \circ {\rm Nm}_{K/\Qp})] \in
\Z[r(\pi)]$$
where $[\ ]$ denotes the associated class in the Grothendieck group. To
show this one gives a description of $H^{n-1}(X,\Lscr_{\xi})^{\rm
  \Z\cross_p}$ in which the $[\psi^i(\rho)]$ occur.
This way one gets sufficient information to see that the map
$$\pi \asr r(\pi\vdual \otimes \vert\ \vert^{1-n \over 2})$$
defines the local Langlands correspondence.

\medskip

I now briefly describe the contents of the various sections. The first
chapter starts with an introductory section on local abelian class field
theory which is reformulated to give the local Langlands correspondence for
$\GL_1$. The next section contains the formulation of the general
correspondence. The following two chapters intend to explain all terms and
notations used in the formulation of the local Langlands correspondence.
We start with some basic definitions in the theory of representations of
reductive $p$-adic groups and give the Langlands classification of
irreducible smooth representations of $\GL_n(K)$. In some cases I did not
find references for the statements (although everything is certainly well
known) and I included a short proof. I apologize if some of those proofs
are maybe somewhat laborious. After a short interlude
about generic and square-integrable representations we come to the
definition of $L-$ and $\eps$-factors of pairs of representations. In the
following chapter we explain the Galois theoretic side of the
correspondence.

The fourth chapter starts with the proof of the correspondence in the
unramified case. Although this is not needed in the sequel, it might be an
illustrating example. After that we return to the general case and give a
number of sketchy arguments to reduce the statement of the existence of a
unique bijection satisfying certain properties to the statement of the
existence of a map satisfying these properties. The third section contains
a small ``dictionary'' which translates certain properties of irreducible
admissible representations of $\GL_n(K)$ into properties of the associated
Weil-Deligne representation. In the fourth section the
construction of the correspondence is given. It uses Jacquet-Langlands
theory, and the cohomology of the sheaf of vanishing cycles on a certain
inductive system of formal schemes. These notions are explained in
the last chapter.

\medskip

Nothing of this treatise is new. For each of the topics there is a number
of excellent references and survey articles. In many instances I just
copied them (up to reordering). In addition to original articles my main
sources, which can (and should) be consulted for
more details, were [CF], [AT], [Neu], [Ta2] (for the number theoretic
background), [Ca], [BZ1], [Cas1], [Ro] (for the background on representation
theory of $p$-adic groups), and [Kud] (for a
survey on ``non-archimedean local Langlands''). Note that
this is of course a personal choice. I also benefited from the opportunity
to listen to the series of lectures of M.\ Harris and G.\ Henniart on the local
Langlands correspondence during the automorphic semester at the IHP in
Paris in spring 2000. I am grateful to the European network in Arithmetic
Geometry and to M. Harris for enabling me to participate in this semester.

This work is intended as a basis for five lectures at the summer school on
``Automorphic Forms on $\GL(n)$'' at the ICTP in Triest. I am grateful to
M.S.\ Raghunathan and G.\ Harder for inviting me to give these
lectures. Further thanks go to U.\ G\"ortz, R.\ Hill, N.\ Kr\"amer and C.\
M\"uller who made many helpful remarks on preliminary versions und to
C.\ Boyallian for fruitful discussions. Finally I
would like to thank the ICTP to provide a pleasant atmosphere during the
summer school.


\bigskip\bigskip

{\partifont Notations}

\parskip

Throughout we fix the following notations and conventions
\bulletlist
\bulletitem $p$ denotes a fixed prime number.
\bulletitem $K$ denotes a $p$-adic field, i.e.\ a finite extension of
$\Q_p$.
\bulletitem $v_K$ denotes the discrete valuation of $K$ normalized such
that it sends uniformizing elements to 1.
\bulletitem $O_K$ denotes the ring of integers of $K$. Further $\pfr_K$ is the
maximal ideal of $O_K$ and $\pi_K$ a chosen generator of $\pfr_K$.
\bulletitem $\kappa$ denotes the residue field of $O_K$, and $q$ the number
of elements in $\kappa$.
\bulletitem $\vert\ \vert_K$ denotes the absolute value of $K$ which takes
$\pi_K$ to $q^{-1}$.
\bulletitem $\psi$ denotes a fixed non-trivial additive character of $K$
(i.e.\ a continuous homomorphism $K \arr \set{z \in \C}{$\vert z \vert = 1$}$).
\bulletitem $n$ denotes a positive integer.
\bulletitem We fix an algebraic closure $\Kbar$ of $K$ and denote by
$\kgbar$ the residue field of the ring of integers of $\Kbar$. This is
an algebraic closure of $\kappa$.
\bulletitem $K^{\rm nr}$ denotes the maximal unramified extension of $K$ in
$\Kbar$. It is also equal to the union of all finite unramified extensions
of $K$ in $\Kbar$. Its residue field is equal to $\kgbar$ and the canonical
homomorphism $\Gal(K^{\rm nr}/K) \arr \Gal(\kgbar/\kappa)$ is an
isomorphism of topological groups.
\bulletitem $\Phi_K \in \Gal(\kgbar/\kappa)$ denotes the geometric
Frobenius $x \asr x^{1/q}$ and $\sigma_K$ its inverse, the arithmetic
Frobenius $x \asr x^q$. We also denote by $\Phi_K$ and $\sigma_K$ the
various maps induced by $\Phi_K$ resp.\ $\sigma_K$ (e.g.\ on $\Gal(K^{\rm
nr}/K)$).
\bulletitem If $G$ is any Hausdorff topological group we denote by
$G^{\rm ab}$ its maximal abelian Hausdorff quotient, i.e.\ $G^{\rm ab}$ is
the quotient of $G$ by the closure of its commutator subgroup.
\bulletitem If $\Ascr$ is an abelian category, we denote by ${\rm
  Groth}(\Ascr)$ its Grothendieck group. It is the quotient of the free
abelian group with basis the isomorphism classes of objects in $\Ascr$
modulo the relation $[V'] + [V''] = [V]$ for objects $V$, $V'$ and $V''$ in
$\Ascr$ which sit in an exact sequence $0 \arr V' \arr V \arr V'' \arr 0$.
For any abelian group $X$ and any
function $\lambda$ which associates to isomorphism classes of objects in
$\Ascr$ an element in $X$ and which is additive (i.e.\ $\lambda(V) =
\lambda(V') + \lambda(V'')$ if there exists an exact sequence $0 \arr V'
\arr V \arr V'' \arr 0$) we denote the induced homomorphism of abelian
groups ${\rm Groth}(\Ascr) \arr X$ again by $\lambda$.

\endchapter


\chapter{The local Langlands correspondence}

\paragraph{The local Langlands correspondence for GL(1)}

\secstart{} In this introductory section we state the local Langlands
correspondence for $\GL_1$ which amounts to one of the main theorems of
abelian local class field theory. For the sake of brevity we use Galois
cohomology without explanation. Galois cohomology will not be needed in the
sequel.

\secstart{} For any finite extension
$L$ of $K$ of degree $m$ and for $\alpha \in K\cross$ we denote by
$$(\alpha, L/K) \in \Gal(L/K)^{\rm ab}$$
the norm residue symbol of local class field
theory. Using Galois cohomology it can be defined as follows (see e.g.\
[Se1] 2, for an alternative more elementary description see [Neu] chap.\
IV, V):

The group $H^2(\Gal(L/K),L\cross)$ is cyclic of order $m$ and up to a sign
canonically isomorphic to ${1 \over m}\Z/\Z$. We use now the sign convention of
[Se1]. Let $v_{L/K}$ be the generator of $H^2(\Gal(L/K),L\cross)$
corresponding to $-{1 \over m}$. By a theorem of Tate (e.g. [AW] Theorem
12) we know that the map $\Hhat^q(\Gal(L/K),\Z) \arr \Hhat^{q+2}(\Gal(L/K),
L\cross)$ which is given by cup-product with $v_{L/K}$ is an
isomorphism. Now we have
$$\Hhat^{-2}(\Gal(L/K),\Z) = H_1(\Gal(L/K),\Z) = \Gal(L/K)^{\rm ab}$$
and
$$\Hhat^0(\Gal(L/K),L\cross) = K^*/N_{L/K}(L\cross)$$
where $N_{L/K}$ denotes the norm of the extension $L$ of $K$. Hence we get
an isomorphism
$$\varphi_{L/K}\colon \Gal(L/K)^{\rm ab} \arriso K^*/N_{L/K}(L\cross).$$
We set
$$(\alpha, L/K) = \varphi_{L/K}^{-1}([\alpha])$$
where $[\alpha] \in K\cross/N_{L/K}(L\cross)$ is the class of $\alpha \in
K\cross$.

\secstart{}\label{\normsymbolunram} If $L$ is a finite unramified extension
of $K$ of degree $m$ we also have the following description of the norm
residue symbol (cf.\ [Se1] 2.5): Let $\Phi_K \in \Gal(L/K)$ be the geometric
Frobenius (i.e.\ it induces on residue fields the map $\sigma_K^{-1}\colon
x \asr x^{-q}$). Then we have for $\alpha \in K\cross$
$$(\alpha,L/K) = \Phi_K^{v_K(\alpha)}.$$

\secstart{} In the sequel we will only need the isomorphisms
$\varphi_{L/K}$. Nevertheless let us give the main theorem of abelian local
class field theory:

\claim Theorem: The map
$$L \asr \phi(L) := N_{L/K}(L\cross) = \Ker(\ ,L/K)$$
defines a bijection between finite abelian extensions $L$ of $K$ and closed
subgroups of $K\cross$ of finite index. If $L$ and $L'$ are finite abelian
extensions of $K$, we have $L \subset L'$ if and only if $\phi(L) \supset
\phi(L')$. In this case $L$ is characterized as the fixed field of
$(\phi(L),L'/K)$.

\proof: See e.g.\ [Neu] chap. V how to deduce this theorem from the
isomorphism $\Gal(L/K) \cong K\cross/N_{L/K}(L\cross)$ using Lubin Tate
theory. We note that this is a purely local proof.

\secstart{} If we go to the limit over all finite extensions $L$ of
$K$, the norm residue symbol defines an isomorphism
$$\limproj_L \Gal(L/K)^{\rm ab} = \Gal(\Kbar/K)^{\rm ab} \arriso \limproj_L
K\cross/N_{L/K}(L\cross).$$
The canonical homomorphism $K\cross \arr  \limproj_L
K\cross/N_{L/K}(L\cross)$ is injective with dense image and hence we get an
injective continuous homomorphism with dense image, called the {\it
  Artin reciprocity homomorphism}
$$\Art_K\colon K\cross \arr \Gal(\Kbar/K)^{\rm ab}.$$

\secstart{}\label{\defineinertia} Let $O_{\Kbar}$ be the ring of
integers of the algebraic closure $\Kbar$ of $K$. Every element of
$\Gal(\Kbar/K)$ defines an automorphism of $O_{\Kbar}$ which reduces to an
automorphism of the residue field $\kgbar$ of $O_{\Kbar}$. We get a
surjective map $\pi\colon \Gal(\Kbar/K) \arr \Gal(\kgbar/\kappa)$ whose
kernel is by definition the inertia group $I_K$ of $K$. The group
$\Gal(\kgbar/\kappa)$ is topologically generated by the arithmetic Frobenius
automorphism $\sigma_K$ which sends $x \in \kgbar$ to $x^q$. It contains
the free abelian group $\langle \sigma_K \rangle$ generated by $\sigma_K$
as a subgroup.

The fixed field of $I_K$ in $\Kbar$ is $K^{\rm nr}$, the union of all
unramified extensions of $K$ in $\Kbar$. By definition we have an
isomorphism of topological groups
$$\Gal(K^{\rm nr}/K) \arriso \Gal(\kgbar/\kappa).$$

\secstart{} The reciprocity homomorphism is already characterized as follows
(cf.\ [Se1] 2.8): Let $f\colon K\cross \arr \Gal(\Kbar/K)^{\rm ab}$ be a
homomorphism such that:
\indention{(b)}
\litem{(a)} The composition
$$K\cross \arrover{f} \Gal(\Kbar/K) \arr \Gal(\kgbar/\kappa)$$
is the map $\alpha \asr \Phi_K^{v_K(\alpha)}$.
\litem{(b)} For $\alpha \in K\cross$ and for any finite abelian extension
$L$ of $K$ such that $\alpha \in N_{L/K}(L\cross)$, $f(\alpha)$ is trivial
on $L$.

Then $f$ is equal to the reciprocity homomorphism $\Art_K$.

\secstart{}\label{\defWeilgroup} We keep the notations of
\defineinertia. The {\it Weil group of $K$} is the inverse image of
$\langle \sigma_K \rangle$ under $\pi$. It is denoted by $W_K$ and
it sits in an exact sequence
$$0 \arr I_K \arr W_K \arr \langle \sigma_K \rangle \arr 0.$$
We endow it with the unique topology of a locally compact group such that
the projection $W_K \arr \langle \sigma_K \rangle \cong \Z$ is continuous
if $\Z$ is endowed with the discrete topology and such that the induced
topology on $I_K$ equals the the profinite topology induced by the
topology of $\Gal(\Kbar/K)$. Note that this topology is different from the
one which is induced by $\Gal(\Kbar/K)$ via the inclusion $W_K \subset
\Gal(\Kbar/K)$. But the inclusion is still continuous, and it has dense
image.

\secstart{}\label{\Weilaltdef} There is the following alternative
definition of the Weil group: As classes in $H^2$
correspond to extensions of groups, we get for every finite extension $L$
of $K$ an exact sequence
$$1 \arr L\cross \arr W(L/K) \arr \Gal(L/K) \arr 1$$
corresponding to the class $v_{L/K}$. For $L \subset L'$ we get a diagram
$$\matrix{1 \arr & L\cross & \arr & W(L/K) & \arr & \Gal(L/K) & \arr 1 \cr
& \auuleft{N_{L'/L}} && && \auu \cr
1 \arr & L'\cross & \arr & W(L'/K) & \arr & \Gal(L'/K) & \arr 1 \cr}$$
which can be commutatively completed by an arrow $W(L'/K) \arr W(L/K)$ such
that we get a projective system $(W(L/K))_L$ where $L$ runs through the set
of finite extensions of $K$ in $\Kbar$. Its projective limit is the Weil
group of $K$ and the projective limit of the homomorphisms $W(L/K) \arr
\Gal(L/K)$ is the canonical injective homomorphism $W_K \arr \Gal(\Kbar/K)$
with dense image.

\secstart{}\label{\reciprocity} Denote by $W_K^{\rm ab}$ the maximal
abelian Hausdorff quotient of $W_K$, i.e.\ the quotient of $W_K$ by the
closure of its commutator subgroup. As the map $W_K \arr \Gal(\Kbar/K)$ is
injective with dense image, we get an induced injective map
$$W_K^{\rm ab} \air \Gal(\Kbar/K)^{\rm ab}.$$
It follows from \Weilaltdef\ and from the definition of
$$\Art_K\colon K\cross \arr \Gal(\Kbar/K)^{\rm ab}$$
that the image of $\Art_K$ is $W_K^{\rm ab}$.

We get an isomorphism of topological groups
$$\Art_K\colon K\cross \arriso W_K^{\rm ab}.$$
This isomorphism maps $O\cross_K$ onto the abelianization $I_K^{\rm ab}$ of
the inertia group and a uniformizing element to a geometric Frobenius
elements, i.e.\ if $\pi_K$ is a uniformizer, the image of ${\rm
  Art}_K(\pi_K)$ in $\Gal(\kgbar/\kappa)$ is $\Phi_K$.

\secstart{}\label{\LocLang1} We can reformulate \reciprocity\ as follows:
Denote by $\Ascr_1(K)$ the set of isomorphism classes of irreducible complex
representations $(\pi,V)$ of $K\cross = \GL_1(K)$ such that the stabilizer
of every vector in $V$ is an open subgroup of $K$. It follows from the
general theory of admissible representations that every $(\pi,V)$ in
$\Ascr_1(K)$ is one-dimensional (see paragraph 2.1 below). Hence
$\Ascr_1(K)$ is equal to the set of continuous homomorphisms $K\cross \arr
\C\cross$ where we endow $\C$ with the discrete topology.

On the other hand denote by $\Gscr_1(K)$ the set of continuous
homomorphisms $W_K \arr \C\cross = \GL_1(\C)$ where we endow $\C\cross$
with its usual topology. Now a homomorphism $W_K \arr \C\cross$ is
continuous if and only if its restriction to the inertia group $I_K$ is
continuous. But $I_K$ is compact and totally disconnected hence its image
will be a compact and totally disconnected subgroup of $\C\cross$ hence it
will be finite. It follows that a homomorphism $W_K \arr \C\cross$ is
continuous for the usual topology of $\C\cross$ if and only if it is
continuous with respect to the discrete topology of $\C\cross$.

Therefore \reciprocity\ is equivalent to:

\claim Theorem {\rm (Local Langlands for $\GL_1$)}: There is a natural
bijection between the sets $\Ascr_1(K)$ and $\Gscr_1(K)$.

\medskip

The rest of these lectures will deal with a generalization of this theorem
to $\GL_n$.


\paragraph{Formulation of the local Langlands correspondence}

\secstart{} Denote by $\Ascr_n(K)$ the set of equivalence classes
of irreducible admissible representations of $\GL_n(K)$.
On the other hand denote by $\Gscr_n(K)$ the set of equivalence classes of
Frobenius semisimple $n$-dimensional complex Weil-Deligne representations
of the Weil group $W_K$ (see chap.\ 2 and chap.\ 3 for a definition of
these notions).

\secstart{THEOREM} (Local Langlands conjecture for $GL_n$ over $p$-adic
fields):\label{\localLanglands} {\sl There is a unique collection of bijections
$$\rec_{K,n} = \rec_n\colon \Ascr_n(K) \arr \Gscr_n(K)$$
satisfying the following properties:
\assertionlist
\assertionitem For $\pi \in \Ascr_1(K)$ we have
$$\rec_1(\pi) = \pi \circ {\rm Art}^{-1}_K.$$
\assertionitem For $\pi_1 \in \Ascr_{n_1}(K)$ and $\pi_2 \in
\Ascr_{n_2}(K)$ we have
$$\eqalign{L(\pi_1 \times \pi_2,s) &= L(\rec_{n_1}(\pi_1) \otimes
  \rec_{n_2}(\pi_2), s),\cr
\eps(\pi_1 \times \pi_2,s,\psi) &= \eps(\rec_{n_1}(\pi_1) \otimes
\rec_{n_2}(\pi_2), s,\psi).\cr}$$
\assertionitem For $\pi \in \Ascr_n(K)$ and $\chi \in \Ascr_1(K)$ we have
$$\rec_n(\pi\chi) = \rec_n(\pi) \otimes \rec_1(\chi).$$
\assertionitem For $\pi \in \Ascr_n(K)$ with central character
$\omega_{\pi}$ we have
$$\det \circ\,\rec_n(\pi) = \rec_1(\omega_{\pi}).$$
\assertionitem For $\pi \in \Ascr_n(K)$ we have $\rec_n(\pi\vdual) =
\rec_n(\pi)\vdual$ where $(\ )\vdual$ denotes the contragredient.

This collection does not depend on the choice of the additive character
$\psi$.}

\secstart{} As the Langlands correspondence gives a bijection between
representations of $\GL_n(K)$ and Weil-Deligne representations of $W_K$
certain properties of and constructions with representations on the one
side correspond to properties and constructions on the other side. Much of
this is still an open problem. A few ``entries in this dictionary'' are
given by the following theorem. We will prove it in chapter 4.\label{\Langdic}

\claim Theorem: {\sl Let $\pi$ be an irreducible admissible representation
of $\GL_n(K)$ and denote by $\rho = (r,N)$ the $n$-dimensional Weil-Deligne
representation associated to $\pi$ via the local Langlands correspondence.
\assertionlist
\assertionitem The representation $\pi$ is supercuspidal if and only if
$\rho$ is irreducible.
\assertionitem We have equivalent statements
\subindention{{\rm (iii)}}
\llitem{{\rm (i)}} $\pi$ is essentially square-integrable.
\llitem{{\rm (ii)}} $\rho$ is indecomposable.
\llitem{{\rm (iii)}} The image of the Weil-Deligne group $W'_F(\C)$ under
$\rho$ is not contained in any proper Levi subgroup of $\GL_n(\C)$.
\assertionitem The representation $\pi$ is generic if and only if $L(s,{\rm
Ad} \circ \rho)$ has no pole at $s = 1$ (here ${\rm Ad}\colon \GL_n(\C)
\arr \GL(M_n(\C))$ denotes the adjoint representation).}

\secstart{} We are going to explain all occuring notations in the following
two chapters.

\endchapter


\chapter{Explanation of the GL(n)-side}

\paragraph{Generalities on admissible representations}

\secstart{} Throughout this chapter let ${\bf G}$ be a connected reductive
group over $K$ and set $G = {\bf G}(K)$. Then $G$ is a locally compact
Hausdorff group such that the compact open subgroups form a basis for the
neighborhoods of the identity (this is equivalent to the fact that $G$ has
compact open subgroups and they are all profinite). In particular, $G$ is
totally disconnected.

To understand \localLanglands\ and \Langdic\ we will only need
the cases where $\Gbf$ is either a product of $\GL_n$'s or the reductive
group associated to some central skew field $D$ over $K$. Nevertheless, in
the first sections we will consider the general case of reductive groups to
avoid case by case considerations. In fact, almost everywhere we could even
work with an arbitrary locally compact totally disconnected group (see
e.g.\ [Vi] for an exposition).

\secstart{} In the case ${\bf G} = \GL_n$ and hence $G =
\GL_n(K)$, a fundamental system of open neighborhoods of the identity
is given by the open compact subgroups $C_m = 1 + \pi_K^mM_n(O_K)$ for $m
\geq 1$. They are all contained in $C_0 = \GL_n(O_K)$, and it is not
difficult to see that $C_0$ is a maximal open compact subgroup and that any
other maximal open compact subgroup is conjugated to $C_0$ (see e.g.\ [Moe]
2).

\secstart{Definition}:\label{\defadmissible} A representation $\pi\colon G
\arr GL(V)$ on a vector space $V$ over the complex numbers is called {\it
admissible} if it satisfies the following two conditions:
\indention{(b)}
\litem{(a)} $(V,\pi)$ is {\it smooth}, i.e.\ the stabilizer of each vector
$v \in V$ is open in $G$.
\litem{(b)} For every open subgroup $H \subset G$ the space $V^H$ of
$H$-invariants in $V$ is finite dimensional.

We denote the set of equivalence classes of irreducible admissible
representations of $G$ by $\Ascr(G)$. For $G = \GL_n(K)$ we define
$$\Ascr_n(K) = \Ascr(GL_n(K)).$$

Note that the notions of ``smoothness'' and ``admissibility'' are purely
algebraic and would make sense if we replace $\C$ by an arbitrary field. In
fact for the rest of this survey article we could replace $\C$ by an
arbitrary non-countable algebraically closed field of characteristic
zero. We could avoid the ``non-countability assumption'' if we worked
consequently only with admissible representations. Further, most elements
of the general theory even work over algebraically closed fields of
characteristic $\ell$ with $\ell \not= p$ ([Vi]).

\secstart{} As every open subgroup of $G$ contains a compact open subgroup,
a representation $(\pi,V)$ is smooth if and only if
$$V = \bigcup_CV^C$$
where $C$ runs through the set of open and compact subgroups of $G$, and it
is admissible if in addition all the $V^C$ are finite-dimensional.

\secstart{Example}: A smooth one-dimensional representation of $K\cross$ is a
{\it quasi-character of $K\cross$} or by abuse of language a {\it
  multiplicative quasi-character of $K$}, i.e.\ a homomorphism of abelian
groups $K\cross \arr \C\cross$ which is continuous for the discrete or
equivalently for the usual topology of $\C\cross$ (cf.\ \LocLang1).

\secstart{}\label{\defineHecke} Let $\Hscr(G)$ be the {\it Hecke
  algebra of $G$}. Its underlying vector space is the space of locally
constant, compactly supported measures $\phi$ on $G$ with complex
coefficients. It becomes an associative $\C$-algebra (in general without
unit) by the convolution
product of measures. If we choose a Haar measure $dg$ on $G$ we can identify
$\Hscr(G)$ with the algebra of all locally constant complex-valued
functions with compact support on $G$ where the product is given by
$$(f_1 * f_2)(h) = \integral_Gf_1(hg^{-1})f_2(g)\,dg.$$

\secstart{}\label{\Heckecompact} If $C$ is any compact open subgroup of
$G$, we denote by $\Hscr(G//C)$ the subalgebra of $\Hscr(G)$ consisting of
those $\phi \in \Hscr(G)$ which are left- and right-invariant under $C$. If
we choose a Haar measure of $G$, we can identify $\Hscr(G//C)$ with the set
of maps $C\backslash G/C \arr \C$ with finite support. The algebra
$\Hscr(G//C)$ has a unit, given by
$$e_C := \vol(C)^{-1}1_C$$
where $1_C$ denotes the characteristic function of $C$.

If $C' \subset C$ is an open compact subgroup, $\Hscr(G//C)$ is a
$\C$-subalgebra of $\Hscr(G//C')$ but with a different unit element if $C
\not= C'$. We have
$$\Hscr(G) = \bigcup_C \Hscr(G//C).$$

\secstart{}\label{\repasHecke} If $(\pi,V)$ is a smooth representation
of $G$, the space $V$ becomes an $\Hscr(G)$-module by the formula
$$\pi(\phi)v = \integral_G\pi(g)d\phi$$
for $\phi \in \Hscr(G)$. This makes sense as the integral is essentially a
finite sum by \Heckecompact.

As $V = \bigcup_C V^C$ where $C$ runs through the open compact subgroups of
$G$, every vector $v \in V$ satisfies $v = \pi(e_C)v$ for some $C$. In
particular, $V$ is a non-degenerate $\Hscr(G)$-module, i.e.\
$\Hscr(G)\cdot V = V$.

We get a functor from the category of smooth representations of $G$ to the
category of non-degenerate $\Hscr(G)$-modules. This functor is an
equivalence of categories [Ca] 1.4.

\secstart{}\label{\Heckemodule} Let $C$ be an open compact subgroup of $G$
and let $(\pi,V)$ be a smooth representation of $G$. Then the space of
$C$-invariants $V^C$ is stable under $\Hscr(G//C)$. If $V$ is an irreducible
$G$-module, $V^C$ is zero or an irreducible $\Hscr(G//C)$-module. More
precisely we have

\claim Proposition: The functor $V \asr V^C$ is an equivalence of the
category of admissible representations of finite length such that every
irreducible subquotient has a non-zero vector fixed by $C$ with the category
of finite-dimensional $\Hscr(G//C)$-modules.

\proof: This follows easily from [Cas1] 2.2.2, 2.2.3 and 2.2.4.

\secstart{Corollary}:\label{\repirred} {\sl Let $(\pi,V)$ be an admissible
representation of $G$ and let $C$ be an open compact subgroup of $G$ such
that for every irreducible subquotient $V'$ of $V$ we have $(V')^C \not=
0$. Then the following assertions are equivalent:
\assertionlist
\assertionitem The $G$-representation $\pi$ is irreducible.
\assertionitem The $\Hscr(G//C)$-module $V^C$ is irreducible.
\assertionitem The associated homomorphisms of $\C$-algebras $\Hscr(G//C)
\arr \End_{\C}(V^C)$ is surjective.}

\proof: The equivalence of (1) and (2) is immediate from \Heckemodule. The
equivalence of (2) and (3) is a standard fact of finite-dimensional modules
of an algebra (see e.g.\ [BouA] chap. VIII, \pz 13, 4, Prop. 5).

\secstart{Corollary}:\label{\unramone} {\sl Let $(\pi,V)$ be an irreducible
admissible representation of $G$ and let $C \subset G$ be an open compact
subgroup such that $\Hscr(G//C)$ is commutative. Then $\dim_{\C}(V^C) \leq 1$.}

\secstart{} For $G = GL_n(K)$ the hypothesis that $\Hscr(G//C)$ is
commutative is fulfilled for $C = \GL_n(O_K)$. In this case we have
$$\Hscr(G//C) = \C[T_1^{\pm 1},\ldots,T_n^{\pm 1}]^{S_n}$$
where the symmetric group $S_n$ acts by permuting the variables $T_i$.

More generally, let ${\bf G}$ be unramified, which means that there exists
a reductive model of ${\bf G}$ over $O_K$, i.e.\ a flat affine group
scheme over $O_K$ such that its special fibre is reductive and such that
its generic fibre is equal to ${\bf G}$. This is equivalent to the
condition that ${\bf G}$ is quasi-split and split over an unramified
extension [Ti] 1.10. If $C$ is a hyperspecial subgroup of $G$ (i.e.\ it is
of the form ${\tilde \Gbf}(O_K)$ for some reductive model ${\tilde \Gbf}$),
the Hecke algebra $\Hscr(G//C)$ can be identified via the Satake isomorphism
with the algebra of invariants under the
rational Weyl group of $\Gbf$ of the group algebra of the cocharacter group
of a maximal split torus of $\Gbf$ [Ca] 4.1. In particular, it is commutative.

\secstart{}\label{\definedist} Under the equivalence of the categories of
smooth $G$-representations and non-degenerate $\Hscr(G)$-modules the
admissible representations $(\pi,V)$ correspond to those non-degenerate
$\Hscr(G)$-modules such that for any $\phi \in \Hscr(G)$ the operator
$\pi(\phi)$ has finite rank.

In particular, we may speak of the trace of $\pi(\phi)$ if $\pi$ is
admissible. We get a distribution $\phi \asr {\rm Tr}(\pi(\phi))$ which is
denoted by $\chi_{\pi}$ and called the {\it distribution character of
$\pi$}. It is invariant under conjugation.

\secstart{} We keep the notations of \definedist. If
$\{\pi_1,\ldots,\pi_n\}$ is a set of pairwise non-isomorphic 
irreducible admissible representations of $G$, then the set of functionals
$$\{\chi_{\pi_1},\ldots,\chi_{\pi_2}\}$$
is linearly independent (cf.\ [JL] Lemma
7.1). In particular, two irreducible admissible representations with the
same distribution character are isomorphic.

\secstart{}\label{\distchar} Let $(\pi,V)$ be an admissible representation
of $G = {\bf G}(K)$. By a theorem of Harish-Chandra [HC] the distribution
$\chi_{\pi}$ is represented by a locally integrable function on $G$ which is
again denoted by $\chi_{\pi}$, i.e.\ for every $\phi \in \Hscr(G)$ we have
$${\rm Tr}\,\pi(\phi) = \integral_G \chi_{\pi}(g)\,d\phi.$$
The function $\chi_{\pi}$ is locally constant on the set of regular
semisimple (see 5.1.3 for a definition in case $G = \GL_n(K)$) elements in
$G$ (loc.\ cit.), and it is invariant under conjugation. Therefore it
defines a function
$$\chi_{\pi}\colon \{G\}^{\rm reg} \arr \C$$
on the set $\{G\}^{\rm reg}$ of conjugacy classes of regular semisimple
elements in $G$.

\secstart{Proposition} (Lemma of Schur):\label{\schur} {\sl Let $(\pi,V)$
be an irreducible smooth representation of $G$. Every $G$-endomorphism of
$V$ is a scalar.}

\proof: We only consider the case that $\pi$ is admissible (we cannot use
(2.1.17), because its proof uses Schur's lemma hence we should not
invoke (2.1.17) if we do not want to run into a circular argument; a direct
proof of the general case can be found in [Ca] 1.4, it uses the fact that
$\C$ is not countable). For sufficiently small open compact subgroups $C$
of $G$ we have $V^C \not= 0$. Hence we have only to show that every
$\Hscr(G//C)$-endomorphism $f$ of a finite dimensional
irreducible $\Hscr(G//C)$-module $W$ over $\C$ is a scalar. As $\C$ is
algebraically closed, $f$ has an eigenvalue $c$, and $\Ker(f - c\id_W)$ is a
$\Hscr(G//C)$-submodule different from $W$. Therefore $f = c\id_W$.

\secstart{Proposition}\label{\remadmissible}: {\sl Let $(\pi,V)$ be an
  irreducible and smooth complex representation of $G$.
\assertionlist
\assertionitem The representation $\pi$ is admissible.
\assertionitem If $G$ is commutative, it is one-dimensional.}

\proof: The first assertion is difficult and can be found for $G =
\GL_n(K)$ in [BZ1] 3.25. It follows from the fact that every smooth
irreducible representation can be embedded in a representation which is
induced from a smaller group and which is admissible (more precisely it is
supercuspidal, see below). Given (1) the proof of (2) is easy: By (1) we can
assume that $\pi$ is admissible. For any compact open subgroup $C$ of $G$
the space $V^C$ is finite-dimensional and a $G$-submodule. Hence $V = V^C$
for any $C$ with $V^C \not= (0)$ and in particular $V$ is
finite-dimensional. But it is well known that every irreducible
finite-dimensional representation of a commutative group $H$ on a vector
space over an algebraically closed field is
one-dimensional (apply e.g.\ [BouA] chap.\ VIII, \pz 13, Prop.\ 5 to the
group algebra of $H$).

\secstart{Proposition}\label{\linearirred}: {\it Every irreducible smooth
representation of $\GL_n(K)$ is either
one-dimensional or infinite-dimensional. If it is one-dimensional, it is of
the form $\chi \circ \det$ where $\chi$ is a quasi-character of $K\cross$,
i.e.\ a continuous homomorphism $K\cross \arr
\C\cross$.}

\medskip

We leave the proof as an exercise (show e.g.\ that the kernel of a
finite-dimensional representation $\pi\colon \GL_n(K) \arr \GL_n(\C)$ is
open, deduce that $\pi$ is trivial on the subgroup of unipotent upper
triangular matrices $U$, hence $\pi$ is trivial on the subgroup of $\GL_n(K)$
which is generated by all conjugates of $U$ and this is nothing but ${\rm
SL}_n(K)$).

\secstart{}\label{\centchar} Let ${\bf Z}$ be the center of ${\bf G}$. As
$K$ is infinite, ${\bf Z}(K)$ is the center $Z$ of $G$. In the case $G =
GL_n(K)$ we have $Z = K\cross$. For $(\pi,V) \in \Ascr(G)$ we denote by
$\omega_{\pi}\colon Z \arr \C\cross$ its {\it central character}, defined
by
$$\omega_{\pi}(z)\id_V = \pi(z)$$
for $z \in Z$. It exists by the lemma of Schur.

For $G = GL_n(K)$, $\omega_{\pi}$ is a quasi-character of $K\cross$.

\secstart{Proposition}\label{\compactfin} Assume that $G/Z$ is a compact
group. Then every irreducible admissible representation $(\pi,V)$ of $G$ is
finite-dimensional.

\proof: By hypothesis we can find a compact open subgroup $G^0$ of $G$ such
that $G^0Z$ has finite index in $G$ (take for example the group $G^0$
defined in 2.4.1 below). The restriction of an irreducible representation
$\pi$ of $G$ to $G^0Z$ decomposes into finitely many
irreducible admissible representations. By the lemma of Schur, $Z$ acts on each
of these representation as a scalar, hence they are also irreducible
representations of the compact group $G^0$ and therefore they are
finite-dimensional.

\secstart{}\label{\twistbychar} Let $(\pi,V)$ be a smooth representation of
$G$ and let $\chi$ be a quasi-character of $G$. The {\it twisted
representation} $\pi\chi$ is defined as
$$g \asr \pi(g)\chi(g).$$
The $G$-submodules of $(\pi,V)$ are the same as the $G$-submodules of
$(\pi\chi,V)$. In particular $\pi$ is irreducible if and only if $\pi\chi$ is
irreducible. Further, if $C$ is a compact open subgroup
of $G$, $\chi(C) \subset \C\cross$ is finite, and therefore $\chi$ is
trivial on a subgroup $C' \subset C$ of finite index. This shows that $\pi$
is admissible if and only if $\pi\chi$ is admissible.

If $G = \GL_n(K)$ every quasi-character $\chi$ is of the form $\chi' \circ
\det$ where $\chi'$ is a multiplicative quasi-character of $K$
\linearirred, and we write $\pi\chi'$ instead of $\pi\chi$.

\secstart{}\label{\contragred} Let $\pi\colon G \arr \GL(V)$ be a smooth
representation of $G$. Denote by $V\star$ the $\C$-linear dual of $V$. It is a
$G$-module via $(g\lambda)(v) = \lambda(g^{-1}v)$ which is not smooth
if $\dim(V) = \infty$. Define
$$V\vdual = \set{\lambda \in V^*}{${\rm Stab}_{G}(\lambda)$ is open}.$$
This is a $G$-submodule $\pi\vdual$ which is smooth by definition. It
is called the {\it contragredient} of the $G$-module $V$. Further we
have:
\assertionlist
\assertionitem $\pi$ is admissible if and only if $\pi\vdual$ is admissible
and in this case the biduality homomorphism induces an isomorphism $V \arr
(V\vdual)\vdual$ of $G$-modules.
\assertionitem $\pi$ is irreducible if and only if $\pi\vdual$ is
irreducible.
\assertionitem In the case of $G = GL_n(K)$ we can describe the
contragredient also in the following way: If $\pi$ is smooth and
irreducible, $\pi\vdual$ is isomorphic to the representation $g \asr
\pi({}^tg^{-1})$ for $g \in GL_n(K)$.

Assertions (1) and (2) are easy (use that $(V\vdual)^C = (V^C)^*$ for every
compact open subgroup $C$). The last assertion is a theorem of Gelfand and
Kazhdan ([BZ1] 7.3).


\paragraph{Induction and the Bernstein-Zelevinsky classification for GL(n)}

\secstart{}\label{\inductionfromLevi} Fix an ordered partition $\nline =
(n_1,n_2,\ldots,n_r)$ of $n$. Denote by $G_{\nline}$ the algebraic group
$GL_{n_1} \times \cdots \times GL_{n_r}$ considered as a Levi subgroup of
$G_{(n)} = GL_n$. Denote by $P_{\nline} \subset GL_n$ the parabolic subgroup
of matrices of the form
$$\pmatrix{A_1 \cr
& A_2 && * \cr
&&\ldots \cr
&0&&\ldots \cr
&&&& A_r\cr}$$
for $A_i \in GL_{n_i}$ and by $U_{\nline}$ its unipotent radical.
If $(\pi_i,V_i)$ is an admissible representation of $GL_{n_i}(K)$, $\pi_1
\otimes \ldots \otimes \pi_r$ is an admissible reprentation of
$G_{\nline}(K)$ on $W = V_1 \otimes \cdots \otimes V_r$. By extending this
representation to $P_{\nline}$ and by normalized induction we get a
representation $\pi_1 \times \cdots \times \pi_r$ of $GL_n(K)$ whose
underlying complex vector space $V$  is explicitly defined by
$$\bigset{V = }{f\colon GL_n(K) \arr W}{$f$ smooth,
  $f(umg) = \delta_{\nline}^{1/2}(m)(\pi_1 \otimes \cdots \otimes
  \pi_r)(m)f(g)$}{for $u \in U_{\nline}(K)$, $m \in GL_{\nline}(K)$ and $g
  \in GL_n(K)$}.$$
Here we call a map $f\colon GL_n(K) \arr W$ smooth if its stabilizer
$$\set{g \in GL_n(K)}{$f(gh) = f(h)$ for all $h \in GL_n(K)$}$$
is open in $GL_n(K)$ (or equivalently if $f$ is fixed by some open compact
of $\GL_n(K)$ acting by right translation), and $\delta_{\nline}^{1/2}$
denotes the positive square root of the modulus character
$$\delta_{\nline}(m) = \vert \det({\rm Ad}_{U_{\nline}}(m)) \vert.$$
The group $GL_n(K)$ acts on $V$ by right translation.

\secstart{Definition}:\label{\defsupercusp} An irreducible smooth
representation of $GL_n(K)$ is called {\it supercuspidal} if there exists
no proper partition $\nline$ such that $\pi$ is a subquotient of a
representation of the form $\pi_1 \times \cdots \times \pi_r$ where $\pi_i$
is an admissible representation of $GL_{n_i}(K)$. We denote by
$\Ascr^0_n(K) \subset \Ascr_n(K)$ the subset of equivalence classes of
supercuspidal representations of $GL_n(K)$.

\secstart{} Let $\pi_i$ be a smooth representation of $\GL_{n_i}(K)$ for $i
= 1,\ldots,r$. Then $\pi = \pi_1 \times \cdots \times \pi_r$ is a smooth
representation of $\GL_n(K)$ with $n = n_1 + \cdots + n_r$. Further it
follows from the compactness of $\GL_n(K)/P_{\nline}$ that if the $\pi_i$
are admissible, $\pi$ is also admissible ([BZ1] 2.26). Further, by [BZ2] we
have the following

\claim Theorem: If the $\pi_i$ are of finite length (and hence admissible
by \remadmissible) for all $i = 1,\ldots,r$ (e.g.\ if all $\pi_i$ are
irreducible), $\pi_1 \times \cdots \times \pi_r$ is also admissible and of
finite length. Conversely, if $\pi$ is an irreducible admissible
representation of $\GL_n(K)$, there exists a unique partition $n = n_1 +
\cdots + n_r$ of $n$ and unique (up to isomorphism and ordering)
supercuspidal representations $\pi_i$ of $GL_{n_i}(K)$ such that $\pi$ is a
subquotient of $\pi_1 \times \cdots \times \pi_r$.

\secstart{}\label{\definescsupp} If $\pi$ is an irreducible admissible
representation of $\GL_n(K)$ we denote the unique unordered tuple
$(\pi_1,\cdots,\pi_r)$ of supercuspidal representations such that $\pi$ is
a subquotient of $\pi_1 \times \cdots \times \pi_r$ the {\it supercuspidal
support}.

\secstart{Definition}:\label{\defmatrixcoeff} Let $\pi\colon GL_n(K) \arr
GL(V)$ be a smooth representation. For $v \in V$ and $\lambda \in V\vdual$
the map
$$c_{\pi,v,\lambda} = c_{v,\lambda}\colon G \arr \C, \qquad g \asr
\lambda(\pi(g)v)$$
is called the {\it $(v,\lambda)$-matrix coefficient} of $\pi$.

\secstart{}\label{\remmatrixcoeff} Let $(\pi,V)$ be an admissible
representation.
\assertionlist
\assertionitem For $v \in V = (V\vdual)\vdual$ and $\lambda \in
V\vdual$ we have
$$c_{\pi,v,\lambda}(g) = c_{\pi\vdual,\lambda,v}(g^{-1}).$$
\assertionitem If $\chi$ is a quasi-character of $K\cross$ we have
$$c_{\pi\chi,v,\lambda}(g) = \chi(\det(g))c_{\pi,v,\lambda}.$$

\secstart{Theorem}:\label{\describesc} {\sl Let $\pi$ be a smooth
irreducible representation of $\GL_n(K)$. Then the following statements are
equivalent:
\assertionlist
\assertionitem $\pi$ is supercuspidal.
\assertionitem All the matrix coefficents of $\pi$ have compact support modulo
center.
\assertionitem $\pi\vdual$ is supercuspidal.
\assertionitem For any quasi-character $\chi$ of $K\cross$, $\pi\chi$ is
supercuspidal.}

\proof: The equivalence of (1) and (2) is a theorem of Harish-Chandra [BZ1]
3.21. The equivalence of (2), (3) and (4) follows then from \remmatrixcoeff.

\secstart{} For any complex number $s$ and for any admissible
representation we define $\pi(s)$ as the twist of $\pi$ with the character
$\vert\ \vert^s$, i.e.\ the representation $g \asr
\vert\det(g)\vert^s\pi(g)$.

If $\pi$ is supercuspidal, $\pi(s)$ is also supercuspidal. Define a partial
order on $\Ascr^0_n(K)$ by $\pi \leq \pi'$ iff there exists an integer $n
\geq 0$ such that $\pi' = \pi(n)$. Hence every finite interval $\Delta$ is
of the form
$$\Delta(\pi,m) = [\pi, \pi(1), \ldots, \pi(m-1)].$$
The integer $m$ is called the {\it length} of the interval and $nm$ is
called its {\it degree}. We write $\pi(\Delta)$ for
the representation $\pi \times \cdots \times \pi(m-1)$ of $GL_{nm}(K)$.

Two finite intervals $\Delta_1$ and $\Delta_2$ are said to be {\it linked}
if $\Delta_1 \not\subset \Delta_2$, $\Delta_2 \not\subset \Delta_1$, and
$\Delta_1 \cup \Delta_2$ is an interval. We say that $\Delta_1$ {\it
  precedes} $\Delta_2$ if $\Delta_1$ and $\Delta_2$ are linked and if the
minimal element of $\Delta_1$ is smaller than the minimal element of
$\Delta_2$.

\secstart{Theorem} (Bernstein-Zelevinsky classification ([Ze], cf.\ also
[Ro])):\label{\BZclass}
\indention{(2)}
\litem{(1)} {\sl For any finite interval
$\Delta \subset \Ascr^0_n(K)$ of length $m$ the representation
$\pi(\Delta)$ has length $2^{m-1}$. It has a unique irreducible quotient
$Q(\Delta)$ and a unique irreducible subrepresentation $Z(\Delta)$.
\litem{{\rm (2)}} Let $\Delta_1 \subset \Ascr^0_{n_1}(K) ,\ldots, \Delta_r
\subset \Ascr^0_{n_r}(K)$ be finite intervals such that for $i < j$,
$\Delta_i$ does not precede $\Delta_j$ (this is an empty condition if $n_i
\not= n_j$). Then the representation $Q(\Delta_1) \times \cdots \times
Q(\Delta_r)$ admits a unique irreducible quotient
$Q(\Delta_1,\ldots,\Delta_r)$, and the representation $Z(\Delta_1) \times
\cdots \times Z(\Delta_r)$ admits a unique irreducible subrepresentation
$Z(\Delta_1,\ldots,\Delta_r)$.
\litem{{\rm (3)}} Let $\pi$ be a smooth irreducible representation of
$GL_n(K)$. Then it is isomorphic to a representation of the form $Q(\Delta_1,\ldots,\Delta_r)$
(resp.\ $Z(\Delta'_1,\ldots,\Delta'_{r'})$)
for a unique (up to permutation) collection of intervals
$\Delta_1,\ldots,\Delta_r$ (resp.\ $\Delta'_1,\ldots,\Delta'_{r'}$) such
that $\Delta_i$ (resp.\ $\Delta'_i$) does not precede $\Delta_j$ (resp.\
$\Delta'_j$) for $i < j$.
\litem{{\rm (4)}} Under the hypothesis of (2), the representation $Q(\Delta_1)
\times \cdots \times Q(\Delta_r)$ is irreducible if and only if no two of
the intervals $\Delta_i$ and $\Delta_j$ are linked.}

\secstart{}\label{\reformBZclass} For $\pi \in \Ascr^0_n(K)$ the set of
$\pi'$ in $\Ascr^0_n(K)$ which are comparable with $\pi$ with respect to
the order defined above is isomorpic (as an ordered set) to $\Z$, in
particular it is totally ordered. It follows that given a tuple of
intervals $\Delta_i = [\pi_i,\ldots,\pi_i(m_i-1) ]$, $i = 1,\ldots,r$ we can
always permute them such that $\Delta_i$ does not precede $\Delta_j$ for $i
< j$.

Denote by $\Sscr_n(K)$ the set of unordered tuples
$(\Delta_1,\ldots,\Delta_r)$ where $\Delta_i$ is an interval of degree
$n_i$ such that $\sum n_i = n$. Then (2) and (3) of \BZclass\ are
equivalent to the assertion that the maps
$$\eqalign{Q\colon \Sscr_n(K) \arr \Ascr_n(K), &\qquad
(\Delta_1,\ldots,\Delta_r) \asr Q(\Delta_1,\ldots,\Delta_r),\cr
Z\colon \Sscr_n(K) \arr \Ascr_n(K), &\qquad (\Delta_1,\ldots,\Delta_r)
\asr Z(\Delta_1,\ldots,\Delta_r),\cr}$$
are bijections.

The unordered tuple of supercuspidal representations $\pi_i(j)$ for $i =
1,\ldots,r$ and $j = 0,\ldots,m_i-1$ is called the {\it supercuspidal
support}. It is the unique unordered tuple of supercuspidal representations
$\rho_1,\ldots,\rho_s$ such that $\pi = Q(\Delta_1,\ldots,\Delta_s)$ and
$\pi' = Z(\Delta_1,\ldots,\Delta_r)$ is a subquotient of $\rho_1 \times
\cdots \times \rho_s$ ([Ze]).

\secstart{} If $\Rscr_n(K)$ is the Grothendieck group of the category of
admissible representations of $\GL_n(K)$ of finite length and $\Rscr(K) =
\bigoplus_{n\geq 0} \Rscr_n(K)$, then
$$([\pi_1],[\pi_2]) \asr [\pi_1 \times \pi_2]$$
defines a map
$$\Rscr(K) \times \Rscr(K) \arr \Rscr(K)$$
which makes $\Rscr(K)$ into a graded commutative ring ([Ze] 1.9) which is
isomorphic to the ring of polynomials in the indeterminates $\Delta$ for
$\Delta \in \Sscr(K) = \bigcup_{n \geq 1}\Sscr_n(K)$ (loc.\ cit.\ 7.5).

The different descriptions of $\Ascr_n(K)$ via the maps $Q$ and $Z$ define
a map
$$t\colon \Rscr(K) \arr \Rscr(K),\qquad Q(\Delta) \asr Z(\Delta).$$
We have:

\claim Proposition: \indention{(2)}\lnoitem{(1)}\label{\Zelinvolution}The
map $t$ is an involution of the graded ring $\Rscr$.
\litem{(2)} It sends irreducible representations to irreducible
representations.
\litem{(3)} For $\Delta = [\pi,\pi(1),\ldots,\pi(m-1)]$ we have
$$t(Q(\Delta)) = Q(\pi, \pi(1), \ldots, \pi(m-1))$$
where on the right hand side we consider $\pi(i)$ as intervals of length 1.
\litem{(4)} We have
$$\eqalign{t(Q(\Delta_1,\ldots,\Delta_r)) &= Z(\Delta_1,\ldots,\Delta_r),\cr
t(Z(\Delta_1,\ldots,\Delta_r)) &= Q(\Delta_1,\ldots,\Delta_r).\cr}$$

\proof: Assertions (1) and (3) follow from [Ze] 9.15. The second assertion
had been anounced by J.N. Bernstein but no proof has been published. It has
been proved quite recently in [Pr] or [Au1] (see also [Au2]). Assertion (4) is
proved by Rodier in [Ro] th\'eor\`eme 7 under the assumption of (2).

\secstart{} For each interval $\Delta = [\pi, \ldots, \pi(m-1)]$ we set
$$\Delta\vdual = [\pi(m-1)\vdual,\ldots,\pi\vdual] =
[\pi\vdual(1-m),\ldots,\pi\vdual(-1),\pi\vdual].$$
It follows from [Ze] 3.3 and 9.4 (cf.\ also [Tad] 1.15 and 5.6) that we have
$$\eqalign{Q(\Delta_1,\ldots,\Delta_r)\vdual &=
Q(\Delta\vdual_1,\ldots,\Delta\vdual_r),\cr
Z(\Delta_1,\ldots,\Delta_r)\vdual &=
Z(\Delta\vdual_1,\ldots,\Delta\vdual_r).\cr}$$
In particular we see that the involution on $\Rscr$ induced by $[\pi] \asr
[\pi\vdual]$ commutes with the involution $t$ in \Zelinvolution.

\secstart{Example}:\label{\Steinberg} Let $\Delta \subset \Ascr_1(K)$ be
the interval
$$\Delta = (\vert\ \vert^{(1-n)/2}, \vert\ \vert^{(3-n)/2}, \ldots, \vert\
\vert^{(n-1)/2}).$$
The associated representation of the diagonal torus $T \subset
GL_n(K)$ is equal to $\delta_B^{-1/2}$ where $\delta_B(t) = \vert \det {\rm
  Ad}_U(t)\vert_K$ is the modulus character of the adjoint action of $T$ on
the group of unipotent upper triangular matrices $U$ and where $B$ is the subgroup of upper triangular matrices in $GL_n(K)$. Hence we see that
$$\pi(\Delta) = \vert\ \vert^{(1-n)/2} \times \vert\ \vert^{(3-n)/2} \times
\ldots \times \vert\ \vert^{(n-1)/2}$$
consists just of the space of smooth functions on
$B\backslash G$ with the action of $G$ induced by the natural action of $G$
on the flag variety $B\backslash G$. Hence $Z(\Delta)$ is the trivial
representation ${\bf 1}$ of constant functions on $G/B$. The representation
$Q(\Delta) = t(Z({\bf 1}))$ is called the {\it Steinberg representation}
and denoted by ${\rm St}(n)$. It
is selfdual, i.e.\ ${\rm St}(n)\vdual = {\rm St}(n)$ (in fact, it is also
unitary and even square integrable, see the next section). For
$n = 2$ the length of $\pi(\Delta)$ is 2, hence we have ${\rm St}(n) =
\pi(\Delta)/{\bf 1}$.


\paragraph{Square integrable and tempered representations}

\secstart{} We return to the general setting where ${\bf G}$ is an
arbitrary connected reductive group over $K$. Every character
$\alpha\colon {\bf G} \arr \G_m$ defines on $K$-valued points a
homomorphism $\alpha\colon G \arr K\cross$. By composition with the
absolute value $\vert\ \vert_K$ we obtain a homomorphism
$\vert\alpha\vert_K\colon G \arr \R^{>0}$ and we set
$$G^0 = \bigcap_{\alpha}\Ker(\vert\alpha\vert_K).$$
If ${\bf G} = GL_n$ then every $\alpha$ is a power of the determinant,
hence we have
$$\GL_n(K)^0 = \set{g \in GL_n(K)}{$\vert\det(g)\vert_K = 1$}.$$

Let $r$ be a positive real number. We call an admissible representation
$(\pi,V)$ of $G$ {\it essentially $L^r$} if for all $v \in V$ and $\lambda
\in V\vdual$ the matrix coefficient $c_{v,\lambda}$ is $L^r$ on
$G^0$, i.e.\ the integral
$$\integral_{G^0}\vert c_{\lambda,v}\vert^r\,dg$$
exists (where $dg$ denotes some Haar measure of $G^0$).

An admissible representation is called $L^r$ if it is essentially $L^r$ and
if it has a central character \centchar\ which is unitary.

Let $Z$ be the center of $G$. Then the composition $G^0 \arr G \arr G/Z$
has compact kernel and finite cokernel. Hence, if $(\pi,V)$ has a unitary
central character $\omega_{\pi}$, the integral
$$\integral_{Z\backslash G}\vert c_{v,\lambda}\vert^r\,dg$$
makes sense, and $(\pi,V)$ is $L^r$ if and only if this
integral is finite.

\secstart{Proposition}:\label{\Lrimplication} Let $\pi$ be an admissible
representation of $G$ which is $L^r$. Then it is $L^{r'}$ for all $r' \geq
r$.

\proof: This follows from [Si3] 2.5.

\secstart{Definition}: An admissible representation of $G$ is called {\it
essentially square integrable} (resp.\ {\it essentially tempered}) if it is
essentially $L^2$ (resp.\ essentially $L^{2+\eps}$ for all $\eps > 0$). We
have similar definitions by omitting ``essentially''.

By \Lrimplication, any (essentially) square integrable representation is
(essentially) tempered.

\secstart{} The notion of ``tempered'' is explained by the following
proposition (which follows from [Si1] \pz 4.5 and [Si3] 2.6):

\claim Proposition: Let $\pi$ be an irreducible admissible representation
of $G$ such that its central character is unitary. Then the following
assertions are equivalent:
\assertionlist
\assertionitem $\pi$ is tempered.
\assertionitem Each matrix coefficient defines a tempered distribution on
$G$ (with the usual notion of a tempered distribution: It extends
from a linear form on the locally constant functions with compact support
on $G$ to a linear form on the Schwartz space of $G$ (the ``rapidly
decreasing functions on $G$''), see [Si1] for the precise definition in the
$p$-adic setting).
\assertionitem The distribution character of $\pi$ is tempered.

\secstart{Example}: By \describesc\ any supercuspidal representation is
essentially $L^r$ for all $r > 0$. In particular it is essentially square
integrable.

\secstart{}\label{\classL2} If $(\pi,V)$ is any smooth representation of
$G$ which has a central character, then there exists a unique positive real
valued quasi-character $\chi$ of $G$ such that $\pi\chi$ has a unitary central
character (for $G = GL_n(K)$ this is clear as every quasi-character factors
through the determinant \linearirred, for arbitrary reductive groups this
is [Cas] 5.2.5). Hence for $G = GL_n(K)$ the notion of ``essential
square-integrability'' is equivalent to the notion of
``quasi-square-integrability'' in the sense of [Ze]. In particular it
follows from [Ze] 9.3:

\claim Theorem: An irreducible admissible representation $\pi$ of
$\GL_n(K)$ is essentially\break square-integrable if and only if it is of the
form $Q(\Delta)$ with the notations of \BZclass. It is square integrable if
and only if $\Delta$ is of the form $[\rho,\rho(1),\ldots,\rho(m-1)]$ where
the central character of $\rho((m-1)/2)$ is unitary.

\secstart{}\label{\classtempered} We also have the following
characterization of tempered representations in the Bernstein-Zelevinsky
classification (see [Kud] 2.2):

\claim Proposition: An irreducible admissible representation
$Q(\Delta_1,\ldots,\Delta_r)$ of $\GL_n(K)$ is tempered if and only if the
$Q(\Delta_i)$ are square integrable.

\secstart{} If $\pi = Q(\Delta_1,\ldots,\Delta_r)$ is a tempered
representation no two of the intervals $\Delta_i =
[\rho_i,\ldots,\rho_i(m_i-1)]$ are linked as ${\rm cent}(\Delta_i) =
\rho_i((m_i-1)/2)$ has unitary central character and all elements in
$\Delta_i$ different from ${\rm cent}(\Delta_i)$ have a non-unitary central
character. Therefore we have
$$\pi = Q(\Delta_1) \times \cdots \times Q(\Delta_r).$$

\secstart{}\label{\Langclass} Let $\pi = Q(\Delta_1,\ldots,\Delta_r)$ be an
arbitrary irreducible admissible representation. For each $\Delta_i$ there
exists a unique real number $x_i$ such that $Q(\Delta_i)(-x_i)$ is square
integrable. We can order the $\Delta_i$'s such that
$$y_1 := x_1 = \cdots = x_{m_1} > y_2 := x_{m_1+1} = \cdots = x_{m_2} >
\cdots > y_s := x_{m_{s-1}+1} = \cdots = x_r.$$
In this order $\Delta_i$ does not precede $\Delta_j$ for $i < j$ and all
$\Delta_i$'s which correspond to the same $y_j$ are not linked. For $j =
1,\ldots,s$ set
$$\pi_j = Q(\Delta_{m_{j-1}+1})(-y_j) \times \cdots \times
Q(\Delta_{m_j})(-y_j)$$
with $m_0 = 0$ and $m_s = r$. Then all $\pi_j$ are irreducible tempered
representation, and $\pi$ is the unique irreducible quotient of $\pi_1(y_1)
\times \cdots \times \pi_s(y_s)$. This is nothing but the Langlands
classification which can be generalized to arbitrary reductive groups (see
[Si1] or [BW]).


\paragraph{Generic representations}

\secstart{}\label{\fixchar} Fix a non-trivial additive quasi-character
$\psi\colon F \arr \C\cross$ and let $n(\psi)$ be the {\it exponent of
  $\psi$}, i.e.\ the largest integer $n$ such that $\psi(\pi_K^{-n}O_K) =
1$.

\secstart{} Let $U_n(K) \subset \GL_n(K)$ be the subgroup of unipotent upper
triangular matrices and define a one-dimensional representation
$\theta_{\psi}$ of $U_n(K)$ by
$$\theta_{\psi}((u_{ij})) = \psi(u_{12} + \cdots + u_{n-1,n}).$$
If $\pi$ is any representation of $GL_n(K)$ we can consider the space of
homomorphisms of $U_n(K)$-modules
$$\Hom_{U_n(K)}(\pi\restricted{U_n(K)},\theta_{\psi}).$$
If $\pi$ is smooth and irreducible we call $\pi$ {\it generic} if this
space is non-zero.

\secstart{}\label{\propgeneric} In the next few sections we collect some
facts about generic representations of $\GL_n(K)$ which can be found in
[BZ1], [BZ2] and [Ze]. Note that in loc.\ cit.\ the term ``non-degenerate''
is used instead of ``generic''. First of all we have:

\claim Proposition: \indention{(2)}\lnoitem{(1)}The representation $\pi$ is
generic if and only if $\pi\vdual$ is generic.
\litem{(2)} For all multiplicative quasicharacters $\chi\colon K\cross \arr
\C\cross$, $\pi$ is generic if and only if $\chi\pi$ is generic.
\litem{(3)} The property of $\pi$ being generic does not depend of the
choice of the non-trivial additive character $\psi$.

\secstart{} Via the Bernstein-Zelevinsky classification we have the following
characterization of generic representations ([Ze] 9.7):

\claim Theorem: An irreducible admissible representation $\pi =
Q(\Delta_1,\ldots,\Delta_r)$ is generic if and only if no two segments
$\Delta_i$ are linked. In particular we have
$$\pi \cong Q(\Delta_1) \times \cdots \times Q(\Delta_r).$$

\secstart{Corollary}: {\sl Every essentially tempered (and in particular
every supercuspidal) representation is generic.}

\secstart{} If $(\pi,V)$ is generic, it has a Whittaker model: Choose a
$$0 \not= \lambda \in \Hom_{U(K)}(\pi\restricted{U(K)},\theta_{\psi})$$
and define a map
$$\eqalign{V &\arr \set{f\colon GL_n(K) \arr \C}{$f(ug) = \theta(u)f(g)$
    for all $g \in GL_n(K)$, $u \in U(K)$},\cr
v &\asr (g \asr \lambda(\pi(g)v))\cr}.$$
This is an injective homomorphism of $GL_n(K)$-modules if $GL_n(K)$ acts on
the right hand side by right translation, and we call its image the {\it
  Whittaker model of $\pi$ with respect to $\psi$} and denote it by
$\Wscr(\pi,\psi)$.

\secstart{}\label{\genericimportant} The concept of a generic
representation plays a fundamental role in the theory of automorphic forms:
If $\pi$ is an irreducible admissible representation of the adele valued
group $\GL_n(\Abf_L)$ for a number field $L$, it can be decomposed in a
restricted tensor product 
$$\pi = \bigotimes_v \pi_v$$
where $v$ runs through the places of $L$ and where $\pi_v$ is an admissible
irreducible representation of $\GL_n(L_v)$ (see Flath [Fl] for the
details). If $\pi$ is cuspidal, all the $\pi_v$ are generic by Shalika
[Sh].


\paragraph{Definition of L- and epsilon-factors}

\secstart{} Let $\pi$ and $\pi'$ be smooth irreducible representations of
$\GL_n(K)$ and of $GL_{n'}(K)$ respectively. We are
going to define $L$- and $\eps$-factors of the pair $(\pi,\pi')$. We
first do this for supercuspidal (or more generally for generic)
representation and then use the Bernstein-Zelevinsky classification to make the
general definition.

Assume now that our fixed non-trivial additive character $\psi$ \fixchar\
is unitary, i.e.\ $\psi^{-1} = \ygbar$. Let $\pi$ and $\pi'$ be generic
representations of $\GL_n(K)$ and $\GL_{n'}(K)$ respectively. To define
$L$- and $\eps$-factors $L(\pi \times \pi',s)$ and $\eps(\pi \times
\pi',s, \psi)$ we follow [JPPS1].

Consider first the case $n = n'$. Denote by $\Sscr(K^n)$ the set of locally
constant functions $\phi\colon K^n \arr \C$ with compact support. For
elements $W \in \Wscr(\pi,\psi)$, $W' \in \Wscr(\pi',\ygbar)$ in the
Whittaker models and for any
$\phi \in \Sscr(K^n)$ define
$$Z(W,W',\phi,s) = \integral_{U_n(K)\backslash
  \GL_n(K)}W(g)W'(g)\phi((0,\ldots,0,1)g)\vert \det(g)\vert^s\,dg$$
where $dg$ is a $GL_n(K)$-invariant measure on $U_n(K)\backslash
  \GL_n(K)$. This is absolutely convergent if $\Re(s)$ is sufficiently
  large and it is a rational function of $q^{-s}$. The set
$$\set{Z(W,W',\phi,s)}{$W \in \Wscr(\pi,\psi)$, $W' \in \Wscr(\pi',\ygbar)$
  and $\phi \in \Sscr(K^n)$}$$
generates a fractional ideal in $\C[q^s,q^{-s}]$ with a unique generator
$L(\pi \times \pi',s)$ of the form $P(q^{-s})^{-1}$ where $P \in \C[X]$ is
a polynomial such that $P(0) = 1$.

Further $\eps(\pi \times \pi',s,\psi)$ is defined by the equality
$${Z(\Wtilde,\Wtilde',1-s,\hat{\phi}) \over L(\pi\vdual \times \pi'\vdual,
  1-s)} = \omega_{\pi'}(-1)^n\eps(\pi \times \pi',s,\psi)
  {Z(W,W',s,\phi) \over L(\pi \times \pi',s)}.$$
Here we define $\Wtilde$ by $\Wtilde(g) =
W(w_n\,{}^tg^{-1})$ where $w_n \in GL_n(K)$ is the permutation matrix
corresponding to the longest Weyl group element (i.e.\ to the permutation
which sends $i$ to $n+1-i$). Because of \contragred\ this is an element of
$\Wscr(\pi\vdual, \ygbar)$. In the same way we define $\Wtilde' \in
\Wscr(\pi'\vdual,\psi)$. Finally $\hat{\phi}$ denotes the Fourier transform
of $\phi$ with respect to $\psi$ given by
$$\hat{\phi}(x) = \integral_{K^n} \phi(y)\psi({}^ty\,x)\,dy$$
for $x \in K^n$.

Now consider the case $n' < n$. For $W \in \Wscr(\pi,\psi)$, $W' \in
\Wscr(\pi',\ygbar)$ and for $j = 0,1,\ldots,n-n'-1$ define
$$\eqalign{Z(W,W',j,s) = \integral_{U_{n'}(K)\backslash \GL_{n'}(K)}
  \integral_{M_{j \times n'}(K)} & W(\pmatrix{g & 0 & 0\cr x & I_j & 0 \cr
    0 & 0 & I_{n-n'-j} \cr})W'(g)\cr
&\cdot\vert\det(g)\vert^{s-(n-n')/2}\,dx\,dg\cr}$$
where $dg$ is a $GL_{n'}(K)$-invariant measure on $U_{n'}(K)\backslash
  \GL_{n'}(K)$ and $dx$ is a Haar measure on the space of $(j \times
  n')$-matrices over $K$. Again this is absolutely convergent if $\Re(s)$
  is sufficiently large, it is a rational function of $q^{-s}$ and these
  functions generate a fractional ideal with a unique generator
$L(\pi \times \pi',s)$ of the form $P(q^{-s})^{-1}$ where $P \in \C[X]$ is
a polynomial such that $P(0) = 1$. In this case $\eps(\pi \times
\pi',s,\psi)$ is defined by
$${Z(w_{n,n'}\Wtilde,\Wtilde',n-n'-1-j,1-s) \over L(\pi\vdual \times
  \pi'\vdual,1-s)} = \omega_{\pi'}(-1)^{n-1}\eps(\pi \times \pi', \psi
,s){Z(W,W',j,s) \over L(\pi \times \pi',s)}$$
where $w_{n,n'}$ is the matrix $\twomatrix{I_{n'}}{0}{0}{w_{n-n'}} \in
GL_n(K)$.

Finally for $n' > n$ we define
$$L(\pi \times \pi',s) = L(\pi' \times \pi,s), \qquad \eps(\pi \times
\pi',\psi,s) = \eps(\pi' \times \pi,\psi,s).$$

In all cases $L(\pi \times \pi',s)$ does not depend on the choice of
$\psi$, and $\eps(\pi \times
\pi',s,\psi)$ is of the form $cq^{-fs}$ for a non-zero complex number $c$
and an integer $f$ which depend only on $\pi$, $\pi'$ and $\psi$.

This finishes the definition of the $L$- and the $\eps$-factor for the
generic case (and in particular for the supercuspidal case). Note that if
$\pi$ and $\pi'$ are supercuspidal, we have
$$L(\pi \times \pi',s) = \prod_{\chi}L(\chi,s)\lformno$$
where $\chi$ runs over the unramified quasi-characters of $K\cross$ such that
$\chi\pi'\vdual \cong \pi$ and where $L(\chi,s)$ is the $L$-function of a
character as defined in Tate's thesis (see also below). In particular we
have $L(\pi \times \pi',s) = 1$ for $\pi \in \Ascr^0_n(K)$ and $\pi' \in
\Ascr^0_{n'}(K)$ with $n \not= n'$.\sublabel{\describeL}

It seems that there is no such easy way to define $\eps(\pi \times
\pi',\psi,s)$ for supercuspidal $\pi$ and $\pi'$. However, Bushnell and
Henniart [BH] prove that $\eps(\pi \times \pi\vdual,\psi,1/2) =
\omega_{\pi}(-1)^{n-1}$ for every irreducible admissible representation
$\pi$ of $\GL_n(K)$.

\secstart{} From the definition of the $L$- and $\eps$-factor in the
supercuspidal case we deduce the definition of the $L$- and
$\eps$-factor for pairs of arbitrary smooth irreducible representations
$\pi$ and $\pi'$ by the following inductive relations using \BZclass
(cf. [Kud]):
\assertionlist
\assertionitem We have $L(\pi \times \pi',s) = L(\pi' \times \pi,s)$ and
$\eps(\pi \times \pi',\psi,s) = \eps(\pi' \times \pi,\psi,s)$.
\assertionitem If $\pi$ is of the form $Q(\Delta_1,\ldots,\Delta_r)$
\BZclass\ and if $\pi'$ is arbitrary, then
$$\eqalign{L(\pi \times \pi',s) &= \prod_{i=i}^rL(Q(\Delta_i) \times \pi',s)\cr
\eps(\pi \times \pi',\psi,s) &= \prod_{i=i}^r\eps(Q(\Delta_i) \times
\pi',\psi,s).}$$
\assertionitem If $\pi$ is of the form $Q(\Delta)$, $\Delta =
[\sigma,\sigma(r-1)]$ and $\pi' = Q(\Delta')$, $\Delta' = [\sigma',
\sigma'(r'-1)]$ with $r' \geq r$, then
$$\eqalign{L(\pi \times \pi',s) &= \prod_{i=1}^rL(\sigma \times \sigma',
  s+r+r'-1)\cr
\eps(\pi \times \pi',\psi,s) &=
\prod_{i=1}^r\biggl(\Bigl(\prod_{j=0}^{r+r'-2i}\eps(\sigma \times
\sigma', \psi, s + i + j - 1)\Bigr)\cr
&\quad\times\Bigl(\prod_{j=0}^{r+r'-2i-1}{L(\sigma\vdual \times
  \sigma'\vdual, 1-s-i-j) \over L(\sigma \times
  \sigma',s+i+j-1)}\Bigr)\biggr).\cr}$$

\secstart{} Let ${\bf 1}\colon K\cross \arr \C\cross$ be the trivial
multiplicative character. For any smooth irreducible representation $\pi$
of $GL_n(K)$ we define
$$\eqalign{L(\pi,s) &= L(\pi \times {\bf 1},s),\cr
\eps(\pi,\psi,s) &= \eps(\pi \times {\bf 1},\psi,s).\cr}$$
For $n = 1$, $L(\pi,s)$ and $\eps(\pi,\psi,s)$ are the local $L$- and
$\eps$-factors defined in Tate's thesis. For $n > 1$ and $\pi$
supercuspidal, we have $L(\pi,s) = 1$, while $\eps(\pi,\psi,s)$ is
given by a generalized Gauss sum [Bu].

\secstart{}\label{\conductor} Let $(\pi,V)$ be a smooth and irreducible
representation of $\GL_n(K)$. For any non-negative integer $t$ define
$$K_n(t) = \set{\twomatrix{a}{b}{c}{d} \in \GL_n(O_K)}{$c \in M_{1\times
    n-1}(\pi^t_KO_K)$, $d \equiv 1 \pmod{\pi_K^tO_K}$}.$$
In particular, we have $K_n(0) = GL_n(O_K)$. The smallest non-negative
integer $t$ such that $V^{K_n(t)} \not= (0)$ is called the {\it conductor}
of $\pi$ and denoted by $f(\pi)$. By [JPPS1] (cf.\ also [CHK]) it is also
given by the equality
$$\eps(\pi,\psi,s) = \eps(\pi,\psi,0)q^{-s(f(\pi) + nn(\psi))}$$
where $n(\psi)$ denotes the exponent of $\psi$ \fixchar.

\endchapter


\chapter{Explanation of the Galois side}

\paragraph{Weil-Deligne representations}

\secstart{} Let $W_K$ be the Weil group of $K$ \defWeilgroup\ and let
$\varphi_K\colon W_K \arr \Gal(\Kbar/K)$ be the canonical homomorphism.

A {\it representation of $W_K$} (resp.\ {\it of $\Gal(\Kbar/K)$}) is a
continuous homomorphism $W_K \arr \GL(V)$ (resp.\ $\Gal(\Kbar/K) \arr
\GL(V)$) where $V$ is a finite-dimensional complex vector space. Denote by
$\Rep(W_K)$ (resp.\ $\Rep(\Gal(\Kbar/K))$) the category of
representations of the respective group.

Note that a homomorphism of a locally profinite group (e.g.\ $W_K$ or
$\Gal(\Kbar/K)$) into $\GL_n(\C)$ is continuous for the usual topology of
$\GL_n(\C)$ if and only if it is continuous for the discrete topology.

\secstart{}\label{\unramchar} For $w \in W_K$ we set
$$\vert w \vert = \vert w \vert_K = \vert {\rm Art}_K^{-1}(w) \vert_K.$$
Then the map $W_K \arr \C\cross$, $w \asr \vert w \vert^s$ is a
one-dimensional representation (i.e.\ a quasi-character) of $W_K$ for every
complex number $s$. All one-dimensional representations of $W_K$ which are
trivial on $I_K$ (i.e.\ which are unramified) are of this form ([Ta1] 2.3.1).

\secstart{} As $\varphi_K$ is injective with dense image, we can identify
$\Rep(\Gal(\Kbar/K))$ with a full subcategory of $\Rep(W_K)$. A
representation in this subcategory is called of {\it Galois-type}. By [Ta2]
1.4.5 a representation $r$ of $W_K$ is of Galois-type if and only if its
image $r(W_K)$ is finite.

Conversely, by [De2] \pz 4.10 and \unramchar\ every irreducible
representation $r$ of $W_K$ is of the form $r = r' \otimes \vert\ \vert^s$
for some complex number $s$ and for some representation $r'$ of Galois
type.

\secstart{} A representation of Galois-type of $W_K$ is irreducible if and
only if it is irreducible as a representation of $\Gal(\Kbar/K)$. Further,
if $\sigma$ is any irreducible representation of $W_K$, it is of Galois
type if and only if the image of its determinant $\det \circ \sigma$ is a
subgroup of finite order of $\C\cross$.

\secstart{} Let $L$ be a finite extension of $K$ in $\Kbar$. Then we have a
canonical injective homomorphism $W_L \arr W_K$ with finite cokernel. Hence
restriction and induction of representations give functors
$$\eqalign{{\rm res}_{L/K}\colon &\Rep(W_K) \arr \Rep(W_L)\cr
{\rm ind}_{L/K}\colon &\Rep(W_L) \arr \Rep(W_K)\cr}$$
satisfying the usual Frobenius reciprocity.

More precisely, any representation of $W_K$ becomes a representation  of
$W_L$ by restriction $r \asr r\restricted{W_K}$. This defines the map ${\rm
res}_{L/K}$. Conversely, let $r\colon W_L \arr \GL(V)$ be a representation
of $W_L$. Then we define ${\rm ind}_{L/K}(r)$ as the representation of
$W_K$ whose underlying vector space consists of the continuous maps
$f\colon W_K \arr V$ such that $f(xw) = r(x)f(w)$ for all $x \in W_L$ and
$w \in W_K$.

Note that in the context of the cohomology of abstract groups this functor
``induction'' as defined above is often called ``coinduction''.

\secstart{Definition}:\label{\defWDrep} A {\it Weil-Deligne representation
  of $W_K$} is a pair $(r,N)$ where $r$ is a representation of $W_K$ and
  where $N$ is a $\C$-linear endomorphism of $V$ such that
$$r(\gamma)\,N = |{\rm Art}_K^{-1}(\gamma)|_K\,N\,r(\gamma)\lformno$$
for $\gamma \in W_K$.\sublabel{\WDRepequal}

It is called {\it Frobenius semisimple} if $r$ is semisimple.

\secstart{Remark}:\label{\remarkWDrep} Let $(r,N)$ be a Weil-Deligne
representation of $W_K$.
\assertionlist
\assertionitem Let $\gamma \in W_K$ be an element corresponding to a
uniformizer $\pi_K$ via $\Art_K$. Applying \WDRepequal\ we see that $N$ is
conjugate to $qN$, hence every eigenvalue of $N$ must be zero
which shows that $N$ is automatically nilpotent.
\assertionitem The kernel of $N$ is stable under $W_K$, hence if $(r,N)$ is
irreducible, $N$ is equal to zero. Therefore the irreducible Weil-Deligne
representations of $W_K$ are simply the irreducible continuous
representations of $W_K$.

\secstart{}\label{\contragredWD} Let $\rho_1 = (r_1,N_1)$ and $\rho_2 =
(r_2,N_2)$ be two Weil-Deligne representations on complex vector
spaces $V_1$ and $V_2$ respectively.

Their tensor product $\rho_1 \otimes \rho_2 = (r,N)$ is the Weil-Deligne
representation on the space $V_1 \otimes V_2$ given by
$$r(w)(v_1 \otimes v_2) = r_1(w)v_1 \otimes r_2(w)v_2, \qquad N(v_1 \otimes
v_2) = N_1v_1 \otimes v_2 + v_1 \otimes N_2v_2$$
for $w \in W_K$ and $v_i \in V_i$, $i = 1,2$.

Further, $\Hom_{\C}(V_1,V_2)$ becomes the vector space of a Weil-Deligne
representation $\Homline(\rho_1,\rho_2) = (r,N)$ by
$$(r(w)\varphi)(v_1) = r_2(w)(\varphi(r_1(w)^{-1}v_1)), \qquad
(N\varphi)(v_1) = N_2(\varphi(v_1)) - \varphi(N_1(v_1))$$
for $\varphi \in \Hom_{\C}(V_1,V_2)$, $w \in W_K$ and $v_1 \in V_1$.

In particular we get the contragredient $\rho\vdual$ of a Weil-Deligne
representation as the representation $\Homline(\rho,{\bf 1})$ where ${\bf
1}$ is the trivial one-dimensional representation.

\secstart{}\label{\realWDrep} Consider $W_K$ as a group scheme over
$\Q$ (not of finite type) which is the limit of the constant group schemes
associated to the discrete groups $W_K/J$ where $J$ runs through the open
normal subgroups of $I_K$. Denote by $W'_K$ the semi-direct product
$$W'_K = W_K \stimes \G_a$$
where $W_K$ acts on $\G_a$ by the rule $wxw^{-1} = \vert w \vert_Kx$. This
is a group scheme (neither affine nor of finite type) over $\Q$ whose
$R$-valued points for some $\Q$-algebra $R$ without non-trivial
idempotents are given by $W_K \stimes R$, and the law of composition is
given by
$$(w_1,x_1)(w_2,x_2) = (w_1w_2, \vert w_2\vert_K^{-1}x_1 + x_2).$$

A Weil-Deligne representation of $W_K$ is the same as a complex
finite-dimensional representation of the group scheme $W'_K$ whose
underlying $W_K$-representation is semisimple (to see this
use the fact that representations of the additive group
on a finite-dimensional vector space over a field in characteristic zero
correspond to nilpotent endomorphisms).

The group scheme $W'_K$ (or also its $\C$-valued points $W_K \stimes \C$)
is called the {\it Weil-Deligne group}.

\secstart{}\label{\interpretWD} It follows from the Jacobson-Morozov
theorem that we can also interpret a
Weil-Deligne representation as a continuous complex semisimple
representation of the group $W_K \times {\rm SL}_2(\C)$. If $\eta$ is such
a representation, we associate a Weil-Deligne representation $(r,N)$ by the
formulas
$$r(w) = \eta(w,\twomatrix{\vert w \vert^{1 \over 2}}{0}{0}{\vert w
\vert^{-{1 \over 2}}})$$
and
$$\exp(N) = \eta(1, \twomatrix{1}{1}{0}{1}).$$
A theorem of Kostant assures that two representations of $W_F \times {\rm
SL}_2(\C)$ are isomorphic if and only if the corresponding Weil-Deligne
representations are isomorphic (see [Ko] for these facts).


\paragraph{Definition of L- and epsilon-factors}

\secstart{} Let $\rho = ((r,V),N)$ be a Frobenius semisimple Weil-Deligne
representation. Denote by $V_N$ the kernel of $N$ and by $V^{I_K}_N$ the
space of invariants in $V_N$ for the action of the inertia group $I_K$.

The $L$-factor of $\rho$ is given by
$$L(\rho,s) = \det(1 - q^{-s}\Phi\restricted{V^{I_K}_N})^{-1}$$
where $\Phi \in W_K$ is a geometric Frobenius. If $\rho$ and $\rho'$ are
irreducible Weil-Deligne representations of dimension $n$, $n'$ respectively,
we have
$$L(\rho \otimes \rho',s) = \prod_{\chi}L(\chi,s)\lformno$$
where $\chi$ runs through the unramified quasi-characters of $K\cross \cong
W_K^{\rm ab}$ such that $\chi \otimes \rho\vdual = \rho'$ (compare
\describeL). In particular $L(\rho \otimes \rho',s) = 1$ for $n \not=
n'$.\sublabel{\WDLfactor}

Fix a non-trivial additive character $\psi$ of $K$ and let $n(\psi)$ be the
largest integer $n$ such that $\psi(\pi_K^{-n}O_K) = 1$. Further let $dx$
be an additive Haar measure of $K$.

To define $\eps(\rho,\psi,s)$ we first define the $\eps$-factor of the
Weil group representation $(r,V)$. Assume first that $V$ is
one-dimensional, i.e.\ $r = \chi$ is a quasi-character
$$\chi\colon W_K^{\rm ab} \arr \C\cross.$$

Let $\chi$ be unramified (i.e.\ $\chi(I_K) = (1)$ or equivalently $\chi =
\vert\ \vert^s$ for some complex number $s$). Then we set
$$\eps(\chi,\psi,dx) = \chi(w)q^{n(\psi)}\vol_{dx}(O_K) =
q^{n(\psi)(1-s)}\vol_{dx}(O_K)$$
where $w \in W_K$ is an element whose valuation is $n(\psi)$.

If $\chi$ is ramified, let $f(\chi)$ be the conductor of $\chi$, i.e.\ the
smallest integer $f$ such that $\chi({\rm Art}_K(1 + \pi_K^fO_K)) = 1$, and
let $c \in K\cross$ be an element with valuation $n(\psi) + f(\chi)$. Then
we set
$$\eps(\chi,\psi,dx) = \integral_{c^{-1}O_K\cross}\chi^{-1}({\rm
  Art}_K(x))\psi(x)dx.$$

The $\eps$-factors attached to $(r,V)$ with $\dim(V) > 1$ are
characterized by the following theorem of Langlands and Deligne [De2]:

\claim Theorem: There is a unique function $\eps$ which associates with
each choice of a local field $K$, a non-trivial additive character $\psi$ of
$K$, an additive Haar measure $dx$ on $K$ and a representation $r$ of $W_K$
a number $\eps(r,\psi,dx) \in \C\cross$ such that
\assertionlist
\assertionitem If $r = \chi$ is one-dimensional $\eps(\chi,\psi,dx)$ is
defined as above.
\assertionitem $\eps(\ ,\psi,dx)$ is multiplicative in exact sequences of
representations of $W_K$ (hence we get an induced homomorphism $\eps(\
,\psi,dx)\colon {\rm Groth}({\rm Rep}(W_K)) \arr \C\cross$).
\assertionitem For every tower of finite extensions $L'/L/K$ and for
every choice of additive Haar measures $\mu_L$ on $L$ and $\mu_{L'}$ on
$L'$ we have
$$\eps({\rm ind}_{L'/L}[r'], \psi \circ {\rm Tr}_{L/K},
\mu_{L}) = \eps([r'], \psi \circ {\rm Tr}_{L'/K}, \mu_{L'})$$
for $[r'] \in {\rm Groth}({\rm Rep}(W_{L'}))$ with $\dim([r']) = 0$.

\medskip

Note that we have $\eps(\chi,\psi,\alpha dx) =
\alpha\eps(\chi,\psi,dx)$ for $\alpha > 0$ and hence via
inductivity $\eps(r,\psi,\alpha dx) =
\alpha^{\dim(r)}\eps(\chi,\psi,dx)$. In particular if $[r] \in {\rm
  Groth}({\rm Rep}(W_K))$ is of dimension $0$, $\eps([r],\psi,dx)$ is
independent of the choice of $dx$.

Now we can define the $\eps$-factor of the Weil-Deligne
representation $\rho = (r,N)$ as
$$\eps(\rho,\psi,s) = \eps(\vert\ \vert^sr,\psi,dx)\det(-\Phi\restricted{V^{I_K}/V^{I_K}_N})$$
where $dx$ is the Haar measure on $K$ which is self-dual with respect to the
Fourier transform $f \asr \fhat$ defined by $\psi$:
$$\fhat(y) = \integral f(x)\psi(xy)\,dx.$$
In other words ([Ta1] 2.2.2) it is the Haar measure for which $O_K$ gets the
volume $q^{-d/2}$ where $d$ is the valuation of the absolute different of
$K$ (if the ramification index $e$ of $K/\Q_p$ is not divided by $p$, we
have $d = e - 1$, in general $d$ can be calculated via higher ramification
groups [Se2]).

Note that $\eps(\rho,\psi,s)$ is not additive in exact sequences of
Weil-Deligne representations as taking coinvariants is not an exact functor.

\secstart{}\label{\WDconductor} Let $\rho$ be an irreducible Weil-Deligne
representation of dimension $n$, then we can define the {\it conductor
  $f(\rho)$ of $\rho$} by the equality
$$\eps(\rho,\psi,s) = \eps(\rho,\psi,0)q^{-s(f(\rho) + nn(\psi))}$$
where $n(\psi)$ denotes the exponent of $\psi$ \fixchar. This is a
nonnegative integer which can be explicitly expressed in terms of higher
ramification groups (e.g.\ [Se2] VI,\pz 2, Ex.\ 2).

\secstart{} For any $m \geq 1$ we define the
Weil-Deligne representations ${\rm Sp}(m) = ((r,V),N)$ by $V = \C e_0 \oplus
\cdots \oplus \C e_{m-1}$ with
$$r(w)e_i = \vert w \vert^ie_i$$
and
$$Ne_i = e_{i+1}\quad (0 \leq i < m-1), \qquad Ne_{m-1} = 0.$$
In this case we have $V_N = V_N^{I_K} = \C e_{m-1}$ and $\Phi e_i =
q^{-i}e_i$ for a geometric Frobenius $\Phi \in W_K$. Hence the
$L$-factor is given by
$$L(\rho,s) = {1 \over 1-q^{1-s-m}}.$$
Let $\psi$ be an additive character such that $n(\psi) = 0$ and let $dx$ be
the Haar measure on $K$ which is self-dual with respect to Fourier
transform as above. Then we have $\eps(r,\psi,dx) =
q^{-md/2}$ where $d$ is the valuation of the absolute different of
$K$. Hence the $\eps$-factor is given by
$$\eps(\rho,\psi,s) = (-1)^{m-1}q^{-md - (m-2)(m-1) \over 2}.$$

\secstart{}\label{\buildupWeilDeligne} A Frobenius semisimple Weil-Deligne
representation $\rho$ is indecomposable if and only if it has the form
$\rho_0 \otimes {\rm Sp}(m)$ for some $m \geq 1$ and with $\rho_0$
irreducible. Moreover, the isomorphism class of $\rho_0$ and $m$ are
uniquely determined by $\rho$ ([De1] 3.1.3(ii)).

Further (as in every abelian category where all objects have finite length)
every Frobenius semisimple Weil-Deligne representation is the direct sum of
unique (up to order) indecomposable Frobenius semisimple Weil-Deligne
representations.

\endchapter


\chapter{Construction of the correspondence}

\paragraph{The correspondence in the unramified case}

\secstart{Definition}: An irreducible admissible representation $(\pi,V)$
of $\GL_n(K)$ is called {\it unramified}, if the space of fixed vectors under
$C = \GL_n(O_K)$ is non-zero, i.e. if its conductor \conductor\ is zero.

\secstart{Example}: A multiplicative quasi-character $\chi\colon K\cross
\arr \C\cross$ is unramified if and only if $\chi(O\cross_K) = \{1\}$. An
unramified quasi-character $\chi$ is uniquely determined by its value
$\chi(\pi_K)$ which does not depend on the choice of the uniformizing
element $\pi_K$. It is of the form $\vert\ \vert^s$ for a unique $s \in
\C/(2\pi i(\log q)^{-1})\Z$.

\secstart{}\label{\classifyunram} Let $(\chi_1,\ldots,\chi_n)$ be a family
of unramified quasi-characters which we can view as intervals of length zero in
$\Ascr^0_1(K)$. We assume that for $i < j$, $\chi_i$ does not precede
$\chi_j$, i.e.\ $\chi_i^{-1}\chi_j \not= \vert\ \vert_K$. Then
$Q(\chi_1,\ldots,\chi_n)$ is an unramified representation of $\GL_n(K)$.
Conversely we have [Cas2]

\claim Theorem: Every unramified representation $\pi$ of $\GL_n(K)$ is
isomorphic to a representation of the form $Q(\chi_1,\ldots,\chi_n)$ where the
$\chi_i$ are unramified quasi-characters of $K\cross$.

\secstart{} An unramified representation $\pi$ of $\GL_n(K)$ is
supercuspidal if and only if $n = 1$ and $\pi$ is an unramified
quasi-character of $K^{\cross}$.

\secstart{} Let $\pi$ be an unramified representation associated to
unramified quasi-characters $\chi_1,\ldots,\chi_n$. This tuple of
unramified quasi-characters induces a homomorphism
$$T/T_c \arr \C\cross$$
where $T \cong (K\cross)^n$ denotes the diagonal torus of $G$ and
where $T_c \cong (O\cross_K)^n$ denotes the unique maximal compact subgroup
of $T$ of diagonal matrices with coefficients in $O\cross_K$.

Thus the set of unramified
representations may be identified with the set of orbits under the
Weyl group $S_n$ of $\GL_n$ in
$$\That = \Hom(T/T_c, \C\cross) = (\C\cross)^n$$
where the last isomorphism is given by the identification
$$T/T_c = \Z^n, \quad \diag(t_1,\ldots,t_n) \asr
(v_K(t_1),\ldots,v_K(t_n)).$$

\secstart{} To shorten notations set $C = \GL_n(O_K)$. The Hecke algebra
$\Hscr(\GL_n(K)//C)$ is commutative and canonically isomorphic to
the $S_n$-invariants of the group algebra ([Ca] 4.1)
$$\C[X^*(\That)] = \C[X_*(T)] \cong \C[t_1^{\pm 1},\ldots,t_n^{\pm 1}].$$

If $(\pi,V)$ is an unramified representation, $V^C$ is one-dimensional
\unramone, hence we get a canonical homomorphism
$$\lambda_{\pi}\colon \Hscr(\GL_n(K)//C) \arr \End(V^C) = \C.$$
For every $h \in \Hscr(\GL_n(K)//C)$ the map
$$\That/\Omega_{\GL_n} \arr \C, \qquad \pi \asr \lambda_{\pi}(h)$$
can be considered as an element in $\C[X^*(\That)]^{S_n}$ and
this defines the isomorphism
$$\Hscr(\GL_n(K)//C) \arriso \C[X^*(\That)]^{S_n}.$$

\secstart{Definition}: An $n$-dimensional Weil-Deligne representation $\rho
= ((r,V),N)$ is called {\it unramified} if $N = 0$ and if $r(I_K) =
\{1\}$.

\secstart{} Every unramified $n$-dimensional Weil-Deligne representation
$\rho = ((r,\C^n),N)$ is uniquely determined by the $\GL_n(\C)$-conjugacy
class of $r(\Phi) =: g_{\rho}$ for a geometric Frobenius $\Phi$. By
definition this element is semisimple and hence we can consider this as an
$S_n$-orbit of the diagonal torus $(\C\cross)^n$ of $\GL_n(\C)$.
Hence we get a bijection $\rec_n$ between unramified representations of
$GL_n(K)$ and unramified $n$-dimensional Weil-Deligne representations.

This is normalized by the following two conditions:
\indention{(ii)}
\litem{(i)} An unramified quasi-character $\chi$ of $\K\cross$ corresponds
to an unramified quasi-character $\rec_1(\chi)$ of $W_K^{\rm ab}$ via the map
$\Art_K$ from local class field theory.
\litem{(ii)} The representation $Q(\chi_1,\ldots,\chi_n)$ \classifyunram\
corresponds to the unramified Weil-Deligne representation
$$\rec_1(\chi_1) \oplus \cdots \oplus \rec_1(\chi_n).$$

By the inductive definition of $L$- and $\eps$-factors it follows that
the bijection $\rec_n$ satisfies condition (2) of
\localLanglands\ (see also \describeL\ and \WDLfactor). Further condition
(3) for unramified characters and condition (4) are clearly okay, and
condition (5) follows from the obvious fact that if unramified elements in
$\Ascr_n(K)$ or in $\Gscr_n(K)$ correspond to the $S_n$-orbit of
$\diag(t_1,\ldots,t_n)$ their contragredients correspond to the orbit of
$\diag(t_1,\ldots,t_n)^{-1} = \diag(t_1^{-1},\ldots,t_n^{-1})$.

\secstart{} From the global point of view, unramified representations are
the ``normal'' ones: If
$$\pi = \bigotimes_v \pi_v$$
is an irreducible admissible representation of the adele valued
group $\GL_n(\Abf_L)$ for a number field $L$ as in \genericimportant, all
but finitely many $\pi_v$ are unramified.


\paragraph{Some reductions}

\secstart{} In this paragraph we sketch some arguments (mostly due to
Henniart) which show that it suffices to show the existence of a family of
maps $(\rec_n)$ satisfying all the desired properties between the set of
isomorphism classes of supercuspidal representations and the set of
isomorphism classes of irreducible Weil-Deligne representations.
We denote by $\Ascr_n^0(K)$ the subset of $\Ascr_n(K)$ consisting of the
supercuspidal representations of $\GL_n(K)$. Further let $\Gscr_n^0(K)$ be
the set of irreducible Weil-Deligne representations in $\Gscr_n(K)$.

\secstart{Reduction to the supercuspidal case}:\label{\firstreduction} {\sl
In order to prove the local Langlands conjecture \localLanglands, it suffices
to show that there exists a unique collection of bijections
$$\rec_n\colon \Ascr_n^0(K) \arr \Gscr^0_n(K)$$
satisfying \localLanglands\ (1) to (5).}

\reasoning: This follows from \BZclass\ and from
\buildupWeilDeligne. More precisely, for any irreducible admissible
representation $\pi \cong Q(\Delta_1,\ldots,\Delta_r)$, with $\Delta_i =
[\pi_i,\pi_i(m_i-1)]$ and $\pi_i \in \Ascr^0_{n_i}(K)$ define
$$\rec_{n_1m_1 + \cdots + n_rm_r}(\pi) = \bigoplus_{i=1}^r\rec_{n_i}(\pi_i)
\otimes {\rm Sp}(m_i).$$
Properties \localLanglands\ (1) to (5) follow then (nontrivially) from the
inductive description of the $L$- and the $\eps$-factors.

\secstart{Reduction to an existence statement}:\label{\existreduction} If
there exists a collection of bijections $(\rec_n)_n$ as in \firstreduction,
it is unique. This follows from the fact that representations $\pi \in
\Ascr^0_n(K)$ are already determined inductively by their $\eps$-factors in
pairs and that by \localLanglands(1) $\rec_1$ is given by class field theory. More precisely, we have the following theorem of Henniart [He3]:

\claim Theorem: Let $n \geq 2$ be an integer and let $\pi$ and $\pi'$ be
representations in $\Ascr^0_n(K)$. Assume that we have an equality
$$\eps(\pi \times \tau,\psi,s) = \eps(\pi' \times \tau,\psi,s)$$
for all integers $r = 1,\ldots,n-1$ and for every $\tau \in
\Ascr^0_r(K)$. Then $\pi \cong \pi'$.

\secstart{Injectivity}:\label{\injectivity} Every collection of maps
$\rec_n$ as in \firstreduction\ is automatically injective: If $\chi$ is a
quasi-character of $K\cross$, its $L$-function $L(\chi,s)$ is given by
$$L(\chi,s) = \cases{(1-\chi(\pi)q^s)^{-1},&if $\chi$ is unramified,\cr
1,&if $\chi$ is ramified.}$$
In particular, it has a pole in $s = 0$ if and only if $\chi = 1$. Hence by
\describeL\ and \WDLfactor\ we have for $\pi, \pi' \in \Ascr^0_n(K)$:
$$\eqalign{\rec_n(\pi) = \rec_n(\pi') &\iff \hbox{$L(\rec_n(\pi)\vdual \otimes
    \rec_n(\pi'),s)$ has a pole in $s = 0$}\cr
&\iff \hbox{$L(\pi\vdual \times \pi',s)$ has a pole in $s = 0$}\cr
&\iff \pi = \pi'.\cr}$$

\secstart{Surjectivity}: {\sl In order to prove the local Langlands
conjecture it suffices to show that there exists a collection of maps
$$\rec_n\colon \Ascr_n^0(K) \arr \Gscr^0_n(K)$$
satisfying \localLanglands\ (1) to (5).}

\reasoning: Because of the preservation of $\eps$-factors in pairs it
follows from \WDconductor\ and \conductor\ that $\rec_n$ preserves
conductors. But by the numerical local Langlands theorem of Henniart [He2]
the sets of elements in $\Ascr^0_n(K)$ and $\Gscr^0_n(K)$ which have the same
given conductor and the same central character are finite and have the same
number of elements. Hence the bijectivity of $\rec_n$ follows from its
injectivity \injectivity.


\paragraph{A rudimentary dictionary of the correspondence}

\secstart{} In this section we give some examples how certain properties of
admissible representations can be detected on the corresponding
Weil-Deligne representation and vice versa. Throughout $(\pi,V_{\pi})$
denotes an admissible irreducible representation of $\GL_n(K)$, and $\rho =
((r,V_r),N)$ the $n$-dimensional Frobenius-semisimple Weil-Deligne
representation associated to it via the local Langlands correspondence
\localLanglands.

\secstart{} First of all, we have of course:

\claim Proposition: The admissible representation $\pi$ is supercuspidal if
and only if $\rho$ is irreducible (or equivalently iff $r$ is irreducible
and $N = 0$).

\secstart{}\label{\basicdic} Write $\pi = Q(\Delta_1,\ldots,\Delta_s)$ in the
Bernstein-Zelevinsky classification \BZclass, where $\Delta_i =
[\pi_i,\ldots,\pi_i(m_i-1)]$ is an interval of supercuspidal
representations of $\GL_{n_i}(K)$.

By \firstreduction\ we have
$$\rho = \bigoplus_{i=1}^s(\rec_{n_i}(\pi_i) \otimes {\rm Sp}(m_i)).$$
Set $\rho_i = ((r_i,V_{r_i}), 0) = \rec_{n_i}(\pi_i)$. The underlying
representation of the Weil group of $\pi_i \otimes {\rm Sp}(m_i)$ is then
given by
$$r_i \oplus r_i(1) \oplus \cdots \oplus r_i(m_i-1)$$
where $r(x)$ denotes the representation $w \asr r(w)\vert w \vert^x$
for any representation $r$ of $W_K$ and any real number $x$. We have
$(r_i(j),0) = \rec_{n_i}(\pi_i(j))$.

Further, if $N_i$ is the nilpotent endomorphism of $\rho_i \otimes {\rm
Sp}(m_i)$, its conjugacy class (which we can consider as a non-ordered
partition of $n_im_i$ by the Jordan normal form) is given by the partition
$$n_im_i = \underbrace{m_i + \cdots + m_i}_{\hbox{$n_i$-times}}.$$
Hence we get:

\claim Proposition: The underlying $W_K$-representation $r$ of $\rho$
depends only on the supercuspidal support $\tau_1,\ldots,\tau_t$ of
$\pi$ \reformBZclass. More precisely, we have an isomorphism of Weil-Deligne
representations
$$(r,0) \cong \rec(\tau_1) \oplus \ldots \oplus \rec(\tau_t).$$
The conjugacy class of $N$ is given by the degree $n_i$ of $\pi_i$ and the
length $m_i$ of the intervals $\Delta_i$ as above. In particular, we have
$N = 0$ if and only if all intervals $\Delta_i$ are of length 1.

\secstart{Example}: The Steinberg representation ${\rm St}(n)$ \Steinberg\
corresponds to the Weil-Deligne representation $\vert\ \vert^{(1-n)/2}{\rm
Sp}(n)$.

\secstart{} Recall from \classifyunram\ that $\pi$ is unramified if and
only if all intervals $\Delta_i$ are of length 1 and consist of an
unramified quasi-character of $K\cross$. Hence \basicdic\ shows that $\pi$
is unramified if and only if $\rho$ is unramified. We used this already in
(4.1).

\secstart{} The ``arithmetic information'' of $W_K$ is encoded in the
inertia subgroup $I_K$. The quotient $W_K/I_K$ is the free group generated by
$\Phi_K$ and hence ``knows'' only the number $q$ of elements
in the residue field of $K$. Therefore Weil-Deligne representations $\rho =
(r,N)$ with $r(I_K) = 1$ should be particularly simple. We call such
representations {\it $I_K$-spherical}. Then $r$ is a semisimple
representation of $<\Phi_K> \cong \Z$. Obviously, every finite-dimensional
semisimple representation of $\Z$ is the direct sum of one-dimensional
representations. Hence $r$ is the direct sum of quasi-characters of $W_K$
which are necessarily unramified.

On the $\GL_n(K)$-side let $\Iscr \subset \GL_n(K)$ be an Iwahori subgroup,
i.e.\ $\Iscr$ is an open compact subgroup of $\GL_n(K)$ which is conjugated
to the group of matrices $(a_{ij}) \in \GL_n(O_K)$ with $a_{ij} \in
\pi_KO_K$ for $i > j$. We say that $\pi$ is {\it $\Iscr$-spherical} if the
space of $\Iscr$-fixed vectors is non-zero. By a theorem of Casselman ([Ca]
3.8, valid for arbitrary reductive groups - with the appropriate
reformulation) an irreducible admissible representation is
$\Iscr$-spherical if and only if its supercuspidal support consists of
unramified quasi-characters. Altogether we get:

\claim Proposition: We have equivalent assertions:
\assertionlist
\assertionitem The irreducible admissible representation $\pi$ is
$\Iscr$-spherical.
\assertionitem The supercuspidal support of $\pi$ consists of unramified
quasi-characters.
\assertionitem The corresponding Weil-Deligne representation $\rho$ is
$I_K$-spherical.

\medskip

By \Heckemodule\ the irreducible admissible $\Iscr$-spherical
representations are nothing but the finite-dimensional irreducible
$\Hscr(\GL_n(K)//\Iscr)$-modules. The structure of the $\C$-algebra
$\Hscr(\GL_n(K)//\Iscr)$ is known in terms of generators and relations
([IM]) and depends only on the isomorphism class of $\GL_n$ over some
algebraically closed field (i.e.\ the based root datum of $\GL_n$) and on
the number $q$.

\secstart{} Finally, we translate several notions which have been defined
for admissible representations of $\GL_n(K)$ into properties of
Weil-Deligne representations:

\claim Proposition: Let $\pi$ be an irreducible admissible representation
of $\GL_n(K)$ and let $\rho = (r,N)$ be the corresponding Weil-Deligne
representation.
\assertionlist
\assertionitem We have equivalent statements
\subindention{{\rm (iii)}}
\llitem{{\rm (i)}} $\pi$ is essentially square-integrable.
\llitem{{\rm (ii)}} $\rho$ is indecomposable.
\llitem{{\rm (iii)}} The image of the Weil-Deligne group $W'_F(\C)$ under
$\rho$ is not contained in any proper Levi subgroup of $\GL_n(\C)$.
\assertionitem We have equivalent statements
\subindention{{\rm (iii)}}
\llitem{{\rm (i)}} $\pi$ is tempered.
\llitem{{\rm (ii)}} Let $\eta$ be a representation of $W_K
\times \SL_2(\C)$ associated to $\rho$ (unique up to isomorphism)
\interpretWD. Then $\eta(W_F)$ is bounded.
\llitem{{\rm (iii)}} Let $\eta$ be as in (ii) and let $\Phi \in W_K$ a
geometric Frobenius. Then $\eta(\Phi)$ has only eigenvalues of absolute
value 1.
\assertionitem The representation $\pi$ is generic if and only if $L(s,{\rm
Ad} \circ \rho)$ has no pole at $s = 1$ (here ${\rm Ad}\colon \GL_n(\C)
\arr \GL(M_n(\C))$ denotes the adjoint representation).

\proof: (1): The equivalence of (i) and (ii) follows from \classL2\ and
\buildupWeilDeligne, the equivalence of (ii) and (iii) is clear as any
factorization through a Levi subgroup $\GL_{n_1}(\C) \times \GL_{n_2}(\C)
\subset \GL_n(\C)$ would induce a decomposition of $\rho$.

(3): This is [Kud] 5.2.2.

(2): The equivalence of (ii) and (iii) follows from the facts that the
   image of the inertia group $I_K$ under $\eta$ is finite, as $I_K$ is
   compact and totally disconnected, and that a subgroup $H$ of semisimple
   elements in $\GL_n(\C)$ is bounded if and only if every element of $H$
   has only eigenvalues of absolute value 1 (use the spectral norm).

Now $\rho$ is indecomposable if and only if $\eta$ is indecomposable. To
   prove the equivalence of (i) and (iii) we can therefore assume by
   \classtempered\ and \basicdic\ that $\rho$ is indecomposable, i.e.\ that
   $\pi = Q(\Delta)$ is essentially square integrable. Let $x \in \R$ be the
   unique real number such that $Q(\Delta)(x)$ is square integrable
   \classL2. Then the description of $\rec(Q(\Delta))$ in \basicdic\ shows
   that
$$\vert\det(\eta(w))\vert = \vert\det(r(w))\vert = \vert w \vert^{nx}.$$
Now $\pi$ is square integrable if and only if its central character is
unitary. But by property (4) of the local Langlands classification the
central character of $\pi$ is given by $\det \circ r$. Hence (iii) is
equivalent to (i) by the following lemma whose proof we leave as an
exercise:

\claim Lemma: Let $\pi$ be a supercuspidal representation of $\GL_n(K)$ and
denote by $\omega_{\pi}$ its central character. For an integer $m \geq
1$ let $\delta$ be the interval $[\pi((1-m)/2),\pi((m-1)/2)]$. Let $\eta$ be
a representation of $W_K \times {\rm SL}_2(\C)$ associated to
$(r,N) = \rec(Q(\Delta))$ \interpretWD. Then for $w \in W_K$ all absolute
values of eigenvalues of $\eta(w)$ are equal to $\vert\omega_{\pi}({\rm
Art}_K^{-1}(w))\vert^{1/n}$. In particular all eigenvalues of $\eta(w)$
have the same absolute value.

\hint: First show the result for $m = 1$ where there is no difference
between $\eta$ and $r$. Then the general result can be checked by making
explicit the Jacobson-Morozov theorem in the case of $\GL_{nm}(K)$.


\paragraph{The construction of the correspondence after Harris and Taylor}

\secstart{}\label{\strategy} Fix a prime $\ell \not= p$ and an isomorphism
of an algebraic closure $\Qdbar_{\ell}$ of $\Q_{\ell}$ with $\C$. Denote by
$\kappa$ the residue field of $O_K$ and by $\kgbar$ an algebraic closure of
$\kappa$. For $m \geq 0$ and $n \geq 1$ let $\Sigma_{K,n,m}$ be the unique (up
to isomorphism) one-dimensional special formal $O_K$-module of $O_K$-height $n$
with Drinfeld level $\pfr^m_K$-structure over $\kbar$. Its deformation
functor on local Artinian $O_K$-algebras with residue field $\kgbar$ is
prorepresented by a complete noetherian local $O_K$-algebra $R_{K,n,m}$
with residue field $\kgbar$. Drinfeld showed that $R_{K,n,m}$ is regular and
that the canonical maps $R_{K,n,m} \arr R_{K,n,m+1}$ are finite and flat. The
inductive limit (over $m$) of the formal vanishing cycle sheaves of
$\Spf(R_{K,n,m})$ with coefficients in $\Qdbar_{\ell}$ gives a collection
$(\Psi^i_{K,\ell,h})$ of infinite-dimensional $\Qdbar_{\ell}$-vector spaces
with an admissible action of the subgroup of $GL_h(K) \times D\cross_{K,1/n}
\times W_K$ consisting of elements $(\gamma,\delta,\sigma)$ such that
$$|{\rm Nrd}\delta||\det\gamma|^{-1}|{\rm Art}_K^{-1}\sigma| = 1.$$
For any irreducible representation $\rho$ of $D\cross_{K,1/n}$ set
$$\Psi^i_{K,\ell,n}(\rho) =
\Hom_{D\cross_{K,1/n}}(\rho,\psi^i_{K,\ell,n}).$$
This is an admissible $(GL_n(K) \times W_K)$-module. Denote by
$[\Psi_{K,\ell,n}(\rho)]$ the virtual representation
$$(-1)^{n-1}\sum_{i=0}^{n-1}(-1)^i[\Psi^i_{K,\ell,n}(\rho)].$$
Then the first step is to prove:

\claim Construction theorem: Let $\pi$ be an irreducible supercuspidal
representation of $\GL_n(K)$. Then there is a (true) representation
$$r_{\ell}(\pi)\colon W_K \arr GL_n(\Qdbar_{\ell}) = \GL_n(\C)$$
such that in the Grothendieck group
$$[\Psi_{K,\ell,n}(\JL(\pi)\vdual)] = [\pi \otimes r_{\ell}(\pi)]$$
where $\JL$ denotes the Jacquet-Langlands bijection between irreducible
representations of $D\cross_{K,1/n}$ and essentially square integrable
irreducible admissible representations of $GL_n(K)$.

\medskip

Using this theorem we can define $\rec_n = \rec_{K,n}\colon \Ascr^0_n(K) \arr
\Gscr_n(K)$ by the formula
$$\rec_n(\pi) = r_{\ell}(\pi\vdual \otimes (|\ |_K \circ
\det)^{(1-n)/2}).$$

That this map satisfies \localLanglands\ (1) - (5) follows from
compatibility of $r_{\ell}$ with many instances of the global Langlands
correspondence. The proof of these compatibilities and also the proof of
the construction theorem follow from an analysis of the bad reduction of
certain Shimura varieties. I am not going into any details here and refer
to [HT].

\secstart{} In the rest of this treatise we explain the ingredients of the
construction of the collection of maps $(\rec_n)$.

\endchapter


\chapter{Explanation of the correspondence}

\paragraph{Jacquet-Langlands theory}

\secstart{}\label{\describeskew} We collect some facts about skew fields
with center $K$ (see e.g.\ [PR] as a reference).

Let ${\rm Br}(K)$ be the Brauer group of
$K$. As a set it can be identified with the set of isomorphism classes of
finite-dimensional division algebras over $K$ with center $K$. For $D, D'
\in {\rm Br}(K)$, $D \otimes D'$ is again a central simple
algebra over $K$, hence it is isomorphic to a matrix algebra $M_r(D'')$ for
some $D'' \in {\rm Br}(K)$. If we set $D\cdot D' := D''$, this defines the
structure of an abelian group on ${\rm Br}(K)$.

This group is isomorphic to $\Q/\Z$ where the homomorphism $\Q/\Z \arr {\rm
  Br}(K)$ is given as follows: For a rational number $\lambda$ with $0 \leq
  \lambda < 1$ we write $\lambda = s/r$ for integers $r$, $s$ which are
  prime to each other and with $r > 0$ (and we make the convention $0 =
  0/1$). Then the associated skew field $D_{\lambda}$ is given by
$D_{\lambda} = K_r[\Pi]$ where $K_r$ is the (unique up to isomorphism)
  unramified extension of $K$ of degree $r$ and where $\Pi$ is an
  indeterminate satisfying the relations $\Pi^r = \pi_k^s$ and $\Pi a =
  \sigma_K(a)\Pi$ for $a \in K_r$.

We call $r$ the index of $D_{\lambda}$. It is the order of $D_{\lambda}$ as
an element in the Brauer group and we have
$$\dim_K(D_{\lambda}) = r^2.$$

If $B$ is any simple finite-dimensional $K$-algebra with center
$K$, it is isomorphic to $M_r(D)$ for some skew field $D$ with center
$K$. Further, $B \tensor{K} L$ is a simple $L$-algebra with center $L$ for
any extension $L$ of $K$. In particular $B \tensor{K} \Kbar$ is isomorphic
to an algebra of matrices over $\Kbar$ as there do not exist any
finite-dimensional division algebras over algebraically closed fields
$\Kbar$ except $\Kbar$ itself.

Conversely, if $D$ is a skew field with center $K$ which is
finite-dimensional over $K$ we can associate the invariant ${\rm inv}(D)
\in \Q/\Z$: As $D \tensor{K} \Kbar$ is isomorphic to
some matrix algebra $M_r(\Kbar)$, we have $\dim_K(D) =
r^2$. The valuation $v_K$ on $K$ extends uniquely to $D$ by the formula
$$v_D(\delta) = {1 \over r}v_K({\rm Nrd}_{D/K}(\delta))$$
for $\delta \in D$. Moreover $D$ is complete in the topology given by this
valuation. It follows from the
definition of $v_D$ that the ramification index $e$ of $D$ over $K$ is
smaller than $r$. Set
$$O_D = \set{\delta \in D}{$v_D(\delta) \geq 0$}, \qquad \Pfr_D =
\set{\delta \in D}{$v_D(\delta) > 0$}.$$
Clearly $\Pfr_D$ is a maximal right and left ideal of $O_D$ and the
quotient $\kappa_D = O_D/\Pfr_D$ is a skew field which is a finite extension of
$\kappa$, hence it is finite and has to be commutative. Let $L \subset D$
be the unramified extension corresponding to the extension $\kappa_D$ of
$\kappa$. As no skew field with center $K$ of dimension $r^2$ can contain a
field of $K$-degree bigger than $r$ we have for the inertia index $f$ of
$D$ over $K$
$$f = [\kappa_D : \kappa] = [L : K] \leq r.$$
Hence the formula $r^2 = ef$ shows that $e = f = r$. Further we have seen
that $D$ contains a maximal unramified subfield. The extension $L/K$ is
Galois with cyclic Galois-group generated by the Frobenius automorphism
$\sigma_K$. By the Skolem-Noether theorem (e.g.\ [BouA] VIII, \pz 10.1),
there exists an element $\delta \in D\cross$ such that
$\sigma_K(x) = \delta x\delta^{-1}$ for all $x \in L$.
Then
$$\inv(D) = v_D(\delta) \in {1 \over r}\Z/\Z \subset \Q/\Z$$
is the invariant of $D$.

\secstart{} For every $D \in {\rm Br}(K)$ we can consider its units as an
algebraic group over $K$. More precisely, we define for every $K$-algebra
$R$
$$D\cross(R) = (D \tensor{K} R)\cross.$$
This is an inner form of $GL_{n,K}$ if $n$ is the index of $D$.

\secstart{} Let $D \in {\rm Br}(K)$ be a division algebra with center $K$ of
index $n$. Let $\{d\}$ be a $D\cross$-conjugacy class of elements in
$D\cross$. The image of $\{d\}$ in $D \tensor{K} \Kbar \cong M_n(\Kbar)$ is
a $GL_n(\Kbar)$-conjugacy class $\{d\}'$ of elements in $GL_n(\Kbar)$ which
does not depend on the choice of the isomorphism $D \tensor{K} \Kbar \cong
M_n(\Kbar)$ as any automorphism of $M_n(\Kbar)$ is an inner
automorphism. Further, $\{d\}'$ is fixed by the natural action of
$\Gal(\Kbar/K)$ on conjugacy classes of $GL_n(\Kbar)$. Hence its
similarity invariants in the sense of [BouA] chap.\ 7, \pz 5 are polynomials
in $K[X]$ and it follows that there is a unique $GL_n(K)$-conjugacy class
of elements $\alpha(\{d\})$ in $GL_n(K)$ whose image in $GL_n(\Kbar)$ is
$\{d\}'$. Altogether we get a canonical injective map $\alpha$ from the set
of $D\cross$-conjugacy classes $\{D\cross\}$ in $D\cross$ into the set of
$GL_n(K)$-conjugacy classes $\{GL_n(K)\}$ in $GL_n(K)$.

The image of $\alpha$ consists of the set of conjugacy classes of elliptic
elements in $GL_n(K)$. Recall that an element $g \in GL_n(K)$ is called
{\it elliptic} if it is contained in a maximal torus $T(K)$ of $GL_n(K)$
such that $T(K)/K\cross$ is compact. Equivalently, $g$ is elliptic if
and only if $K[g]$ is a field.

We call a conjugacy class $\{g\}$ in $GL_n(K)$ {\it semisimple} if it
consists of elements which are diagonalizable over $\Kbar$ or,
equivalently, if $K[g]$ is a product of field extensions for $g \in
\{g\}$. A conjugacy class $\{g\}$ is called {\it regular semisimple} if it
is semisimple and if all eigenvalues of elements in $\{g\}$ in $\Kbar$ are
pairwise different. Note that every elliptic element is semisimple. We make
the same definitions for conjugacy classes in $D\cross$, or equivalently we
call a conjugacy class of $D\cross$ {\it semisimple} (resp.\ {\it regular
semisimple}) if its image under $\alpha\colon \{D\cross\} \arr \{G\}$ is
semisimple (resp.\ regular semisimple).

\secstart{} Denote by $\Ascr^2(G)$ the set of isomorphism classes of
irreducible admissible essentially square integrable representations of
$G$. We now have the following theorem which is due to Jacquet and
Langlands in the case $n = 2$ and due to Rogawski and Deligne, Kazhdan and
Vigneras in general ([Rog] and [DKV]):

\claim Theorem: Let $D$ be a skew field with center $K$ and with index
$n$. It exists a bijection, called {\it Jacquet-Langlands correspondence},
$$\JL\colon \Ascr^2(D\cross) \leftrightarrow \Ascr^2(\GL_n(K))$$
which is characterized on characters by
$$\chi_{\pi} = (-1)^{n-1}\chi_{\JL(\pi)}.\lformno$$
Further $\JL$ satisfies the following conditions:\sublabel{\compchar}
\assertionlist
\assertionitem We have equality of central characters
$$\omega_{\pi} = \omega_{\JL(\pi)}.$$
\assertionitem We have an equality of $L$-functions and of
$\eps$-functions up to a sign
$$L(\pi,s) = L(\JL(\pi),s), \qquad \eps(\pi,\psi,s) =
\eps(\JL(\pi),\psi),s)$$
(for the definition of $L$- and $\eps$-function of irreducible
admissible representations of $D\cross$ see e.g.\ [GJ]).
\assertionitem The Jacquet-Langlands correspondence is compatible with
twist by characters: If $\chi$ is a multiplicative quasi-character of
$K$, we have
$$\JL(\pi(\chi \circ {\rm Nrd})) = \JL(\pi)(\chi \circ
\det).$$
\assertionitem It is compatible with contragredient:
$$\JL(\pi\vdual) = \JL(\pi)\vdual.$$

\secstart{Remark}: Note that for $G = D\cross$ every admissible
representation is essentially square integrable as $D\cross/K\cross$ is
compact.


\paragraph{Special p-divisible O-modules}

\secstart{} Let $R$ be a ring. A {\it $p$-divisible group} over $R$ is an
inductive system $G = (G_n,i_n)_{n\geq 1}$ of finite locally free
commutative group schemes $G_n$ over $\Spec(R)$ and group scheme
homomorphisms $i_n\colon G_{n} \ar G_{n+1}$ such that for all integers $n$
there is an exact sequence
$$0 \arr G_{1} \arvarover(20){i_{n-1} \circ \cdots \circ i_1} G_n
\arrover{p} G_{n-1} \arr 0$$
of group schemes over $\Spec(R)$. We have the obvious notion of a
homomorphism of $p$-divisible groups. This way we get a $\Z_p$-linear category.

As the $G_n$ are finite locally free, their underlying schemes are by
definition of the form $\Spec(A_n)$ where $A_n$ is an $R$-algebra which is
a finitely generated locally free $R$-module. In particular, it makes sense
to speak of the rank of $A_n$ which we also call the rank of $G_n$. From
the exact sequence above it follows that $G_n$ is of rank $p^{nh}$ for some
non-negative locally constant function $h\colon \Spec(R) \arr \Z$ which is
called the {\it height of G}.

\secstart{} Let $G = (G_n)$ be a $p$-divisible group over some ring $R$ and
let $R'$ be an $R$-algebra. Then the inductive system of $G_n \otimes_R R'$
defines a $p$-divisible group over $R'$ which we denote by $G_{R'}$.

\secstart{}\label{\dimpdiv} Let $G = (G_n)$ be a $p$-divisible group over a
ring $R$. If there exists some integer $N \geq 1$ such
that $p^NR = 0$, the Lie algebra $\Lie(G_n)$ is a locally free $R$-module
for $n \geq N$ whose rank is independent of $n \geq N$. We call this rank
the {\it dimension of $G$}. More generally, if $R$ is $p$-adically complete
we define the dimension of $G$ as the dimension of the $p$-divisible group
$G_{R/pR}$ over $R/pR$.

\secstart{} Let $R$ be an $O_K$-algebra. A {\it special $p$-divisible
  $O_K$-module over $R$} is a pair $(G,\iota)$ where $G$ is a $p$-divisible
  group over $R$ and where $\iota\colon O_K \arr \End(G)$ is a homomorphism
  of $\Z_p$-algebras such that for all $n \geq 1$ the $O_K$-action induced
  by $\iota$ on $\Lie(G_n)$ is the same as the $O_K$-action which is
  induced from the $R$-module structure of $\Lie(G)$ via the $O_K$-module
  structure of $R$. In other words the induced homomorphism $O_K
  \tensor{\Zp} O_K \arr \End(\Lie(G_n))$ factorizes through the
  multiplication $O_K \tensor{\Zp} O_K \arr O_K$.

The height ${\rm ht}(G)$ of a special $p$-divisible $O_K$-module
$(G,\iota)$ is always divisible by $[K : \Q_p]$ and we call
${\rm ht}_{O_K}(G) := [K : \Qp]^{-1}{\rm ht}(G)$ the {\it $O_K$-height of
  $(G,\iota)$}.

\secstart{} If $(G = (G_n),\iota)$ is a special $p$-divisible $O_K$-module
over an $O_K$-algebra $R$ and if $R \arr R'$ is an $R$-algebra, we get an
induced $O_K$-action $\iota'$ on $G_{R'}$ and the pair $(G_{R'},\iota')$ is
a special $p$-divisible $O_K$-module over $R'$ which we denote by
$(G,\iota)_{R'}$.

\secstart{}\label{\formnot} Let $k$ be a perfect extension of the residue
class field $\kappa$ of $O_K$. Denote by $W(k)$ the ring of Witt vectors of
$k$. Recall that this is the unique (up to unique isomorphism inducing the
identity on $k$) complete
discrete valuation ring with residue class field $k$ whose maximal ideal is
generated by $p$. Further $W(k)$ has the property that for any complete
local noetherian ring $R$ with residue field $k$ there is a unique local
homomorphism $W(k) \arr R$ inducing the identity on $k$.

In particular, we can consider $W(\kappa)$. It can be identified with
the ring of integers of the maximal unramified extension of $\Q_p$ in
$K$ (use the universal property of the ring of Witt vectors). Set
$$W_K(k) = W(k) \tensor{W(\kappa)} O_K.$$
This is a complete discrete valuation ring of mixed characteristic with residue
field $k$ which is a formally unramified $O_K$-algebra (i.e.\ the image of
$\pfr_K$ generates the maximal ideal of $W_K(k)$). There exists a unique
continuous automorphism $\sigma_K$ of $W(k)$ which induces the
automorphism $x \asr x^q$ on $k$. We denote the induced automorphism
$\sigma_K \otimes \id_{O_K}$ again by $\sigma_K$.

\secstart{Proposition}:\label{\formperfekt} {\sl The category of special
$p$-divisible $O_K$-modules
$(G,\iota)$ over $k$ and the category of triples $(M,F,V)$ where $M$ is a
free $W_K(k)$-module of rank equal to the $O_K$-height and $F$ (resp.\ $V$)
is a $\sigma$- (resp.\ $\sigma^{-1}$-) linear map such that $FV = VF =
\pi_K\id_M$ are equivalent. Via this equivalence there is a canonical
functorial isomorphism
$$M/VM \cong \Lie(G).$$
We call $(M,F,V) = M(G,\iota)$ the Dieudonn\'e module of $(G,\iota)$.}

\proof: To prove this we use covariant Dieudonn\'e theory for $p$-divisible
groups as in [Zi1] for example. Denote by $\sigma$ the usual Frobenius of
the ring of Witt vectors. Covariant Dieudonn\'e theory tells us that there
is an equivalence of the category of $p$-divisible groups over $k$ with the
category of triples $(M',F',V')$ where $M'$ is a free $W(k)$-module of rank
equal to the height of $G$ and with a $\sigma$-linear (resp.\ a
$\sigma^{-1}$-linear) endomorphism $F'$ (resp.\ $V'$) such that $F'V' = V'F' =
p\id_{M'}$ and such that $M'/V'M' = \Lie(G)$. Let us call this functor
$M'$. Let $(G,\iota)$ be a special $p$-divisible $O_K$-module. Then the
Dieudonn\'e module $M'(G)$ is a $W(k) \tensor{\Zp} O_K$-module, the operators
$F'$ and $V'$ commute with the $O_K$-action and the induced homomorphism
$O_K \tensor{\Zp} k \arr \End(M'/V'M')$ factors through the multiplication
$O_K \tensor{\Zp} k \arr k$. We have to construct from these data a triple
$(M,F,V)$ as in the claim of the proposition. To do this write
$$W(k) \tensor{\Zp} O_K = W(k) \tensor{\Zp} W(\kappa) \tensor{W(\kappa)}
O_K = \prod_{\Gal(\kappa/\F_p)} W_K(k).$$
By choosing the Frobenius $\sigma = \sigma_{\Q_p}$ as a generator of
$\Gal(\kappa/\Fp)$ we can identify this group with $\Z/r\Z$ where $p^r = q$.
We get an induced decomposition $M' = \oplus_{i \in \Z/r\Z} M'_i$ where the
$M_i$ are $W_K(k)$-modules defined by
$$M_i = \set{m \in M'}{$(a \otimes 1)m = (1
    \otimes \sigma^{-i}(a))m$ for all $a \in W(\kappa) \subset O_K$}.$$
The operator $F'$ (resp.\ $V'$) is homogeneous of degree $-1$ (resp.\ $+1$)
with respect to this decomposition. By the condition on the $O_K$-action on
the Lie algebra we know that $M'_0/VM'_{r-1} = M'/V'M'$ and hence that
$VM'_{i-1} = M'_i$ for all $i \not= 0$. We set $M = M'_0$ and $V =
(V')^r\restricted{M'_0}$. It follows that we have $M/VM = M'/V'M'$. Further
the action of $\pi_K$ on $M/VM = M'/VM'$ equals the scalar multiplication
with the image of $\pi_K$ under the map $O_K \arr \kappa \arr k$ but this
image is zero. It follows that $VM$ contains $\pi_KM$ and hence we can
define $F = V^{-1}\pi_K$. Thus we constructed the triple $(M,F,V)$ and it
is easy to see that this defines an equivalence of the category of triples
$(M',F',V')$ as above and the one of triples $(M,F,V)$ as in the claim.

\secstart{} Let $(G,\iota)$ be a special $p$-divisible $O_K$-module over a
ring $R$. We call it {\it \'etale} if it is an inductive system of finite
\'etale group schemes. This is equivalent to the fact that its Lie algebra
is zero. If $p$ is invertible in $R$, $(G,\iota)$ will be always \'etale.
 
Now assume that $R = k$ is a perfect field of characteristic $p$ and let
$(M,F,V)$ be its Dieudonn\'e module. Then $(G,\iota)$ is \'etale if and
only if $M = VM$.

In general there is a unique decomposition $(M,F,V) = (M_{\et},F,V) \oplus
(M_{\inf},F,V)$ such that $V$ is bijective on $M_{\et}$ and such that
$V^NM_{\inf} \subset \pi_KM_{\inf}$ for large $N$ (define $M_{\et}$ (resp.\
$M_{\inf}$) as the projective limit over $n$ of $\bigcap_m V^m(M/\pi_K^nM)$
(resp.\ of $\bigcup_m \Ker(V^m\restricted{M/\pi_K^nM})$)). We call the
$W_K(k)$-rank of $M_{\et}$ the {\it \'etale $O_K$-height of $(M,F,V)$ {\rm
    or} of $(G,\iota)$}.

We call $(G,\iota)$ {\it formal} or also {\it infinitesimal} if its \'etale
$O_K$-height is zero.

\secstart{Proposition}:\label{\uniqueformmodule} {\sl Let $k$ be an
algebraically closed field of characteristic $p$. For all non-negative
integers $h \leq n$ there exists up to isomorphism exactly one special
$p$-divisible $O_K$-module of $O_K$-height $n$, \'etale $O_K$-height $h$
and of dimension one. Its Dieudonn\'e module $(M,F,V)$ is the free
$W_K(k)$-module with basis $(d_1,\ldots,d_h,e_1,\ldots,e_{n-h})$ such that
$V$ is given by
$$\eqalign{Vd_i &= d_i, \qquad\qquad i = 1,\ldots,h\cr
Ve_i &= e_{i+1}, \qquad\qquad i = 1,\ldots,n-h-1\cr
Ve_{n-h} &= \pi_Ke_1.\cr}$$
This determines also $F$ by the equality $F = V^{-1}\pi_K$.}

\medskip

The key point to this proposition ist the following lemma due to
Dieudonn\'e:

\claim Lemma: Let $M$ be a free finitely generated $W_K(k)$-module and let
$V$ be a $\sigma_K^a$-linear bijection where $a$ is some integer different
from zero. Then there exists a $W_K(k)$-basis $(e_1,\ldots,e_n)$ of $M$
such that $Ve_i = e_i$.

\medskip

A proof of this lemma in the case of $K = \Qp$ can be found in [Zi1]
6.26. The general case is proved word by word in the same way if one
replaces everywhere $p$ by $\pi_K$.

\proof of the Proposition: Let $(G,\iota)$ be a special $p$-divisible
$O_K$-module as in the proposition and let $(M,F,V)$ be its Dieudonn\'e
module. We use the decomposition $(M,F,V) = (M_{\et},F,V) \oplus
(M_{\inf},F,V)$ and can apply the lemma to the \'etale part. Hence we can
assume that $h = 0$ (note that $M_{\inf}/VM_{\inf} = M/VM$). By definition
of $M_{\inf}$, $V$ acts nilpotent on $M/\pi_KM$. We get a decreasing
filtration $M/\pi_KM \supset V(M/\pi_KM) \supset \cdots \supset
V^N(M/\pi_KM) = (0)$. The successive quotients have dimension 1 because
$\dim_k(M/VM) = 1$. Hence we see that $V^nM \subset \pi_KM$. On the other
hand we have
$$\length_{W_K(k)}(M/V^nM) = n\length_{W_K(k)}(M/VM) = n =
\length_{W_K(k)}(M/\pi_KM)$$
which implies $V^nM = \pi_KM$. Hence we can apply the
lemma to the operator $\pi_K^{-1}V^n$ and we get a basis of elements $f$
satisfying $V^nf = \pi_Kf$. Choose an element $f$ of this basis which does
not lie in $VM$. Then the images of $e_i := V^{i-1}f$ for $i = 1,\ldots,n$
in $M/\pi_KM$ form a basis of the $k$-vector space $M/\pi_KM$. Hence the
$e_i$ form a $W_K(k)$-basis of $M$, and $V$ acts in the desired form.

\secstart{Definition}: We denote the unique formal $p$-divisible
$O_K$-module of height $h$ and dimension $1$ over an algebraically closed
field $k$ of characteristic $p$ by $\Sigma_{h,k}$.

\secstart{}\label{\endoalgebra} Denote by $D'_{K,1/h}$ the ring of
endomorphisms of $\Sigma_{h,k}$ and set $D_{K,1/h} = D'_{K,1/h} \tensor{\Z}
\Q$. Then this is ``the'' skew field with center $K$ and invariant $1/h \in
\Q/\Z$ \describeskew. This follows from the following more general proposition:

\claim Proposition: Denote by $L$ the field of fractions of $W_K(k)$ and
fix a rational number $\lambda$. We write $\lambda = r/s$ with integers $r$
and $s$ which are prime to each other and with $s > 0$ (and with the
convention $0 = 0/1$). Denote by $N_{\lambda} = (N,V)$ the pair consisting
of the vector space $N = L^s$ and of the $\sigma_K^{-1}$-linear bijective
map $V$ which acts on the standard basis via the matrix
$$\pmatrix{0 & 0 & \ldots & 0 & \pi_K^r\cr
1 & 0 && \ldots & 0\cr
0 & 1 & 0 & \ldots & 0 \cr
&&\ldots\cr
0 & \ldots & 0 & 1 & 0 \cr}.$$
Then $\End(N_{\lambda}) = \set{f \in \End_L(N)}{$f \circ V = V \circ f$}$
is the skew field $D_{\lambda}$ with center $K$ and invariant equal to the
image of $\lambda$ in $\Q/\Z$ (cf.\ \describeskew).

\proof: We identify $\F_{q^s}$ with the subfield of $k$ of elements $x$
with $x^{q^s} = x$. This contains the residue field $\kappa$ of $O_K$ and
we get inclusions
$$O_K \subset O_{K_s} := W(\F_{q^s}) \tensor{W(\kappa)} O_K \subset
W_K(k)$$
and hence
$$K \subset K_s \subset L.$$
These extensions are unramified, $[K_s : K] = s$, and $K_s$ can be
described as the fixed field of $\sigma_K^s$ in $L$.

To shorten notations we set $A_{\lambda} = \End(N_{\lambda})$. As
$N_{\lambda}$ does not have any non-trivial $V$-stable subspaces (cf.\ [Zi1]
6.27), $A_{\lambda}$ is a skew field and its center contains $K$. For a
matrix $(u_{ij}) \in \End(L^s)$ an easy explicit calculation shows that
$(u_{ij}) \in A_{\lambda}$ if and only if we have the relations
$$\matrix{u_{11} &= &\sigma_K^{-1}(u_{ss}),&\qquad \cr
u_{i+1,j+1} &= &\sigma_K^{-1}(u_{ij}),&& 1 \leq i,j \leq s-1, \cr
u_{1,j+1} &= &\pi_K^r\sigma_K^{-1}(u_{sj}), && 1 \leq j \leq s-1, \cr
u_{i+1,j} &= &\pi_K^{-r}\sigma_K^{-1}(u_{is}), && 1 \leq i \leq s-1. \cr}$$
It follows that $\sigma_K^s(u_{ij}) = u_{ij}$ for all $i,j$, and hence
$u_{ij} \in K_s$. Further, sending a matrix $(u_{ij}) \in A_{\lambda}$ to
its first column defines a $K_s$-linear isomorphism $A_{\lambda} \cong
K_s^s$, hence $\dim_{K_s}(A_{\lambda}) = s$.

The $K_s$-algebra homomorphism
$$\varphi\colon K_s \tensor{K} A_{\lambda} \arr M_s(K_s), \qquad \alpha
\otimes x \asr \alpha x$$
is a homomorphism of $M_s(K_s)$-left modules and hence it is surjective as the
identity matrix is in its image. Therefore $\varphi$ is bijective. In
particular, $K_s$ is the center of $K_s \tensor{K} A_{\lambda}$ and hence
the center of $A_{\lambda}$ is equal to $K$.

Now define
$$\Pi = \pmatrix{0 & 0 & \ldots & 0 & \pi_K^r\cr
1 & 0 && \ldots & 0\cr
0 & 1 & 0 & \ldots & 0 \cr
&&\ldots\cr
0 & \ldots & 0 & 1 & 0 \cr} \in A_{\lambda}.$$
Then we have the relations $\Pi^s = \pi_K^r$ and $\Pi d = \sigma_K(d)\Pi$ for
$d \in A_{\lambda}$. We get an embedding $D_{\lambda} = K_s[\Pi] \air
A_{\lambda}$ by
$$\eqalign{\Pi &\asr \Pi \cr
K_s \ni \alpha &\asr \pmatrix{\sigma_K^{-1}(\alpha) \cr
& \sigma_K^{-2}(\alpha) \cr
&& \cdots \cr
&&& \alpha \cr} \in A_{\lambda} \subset M_s(K_s)\cr}$$
which has to be an isomorphism because both sides have the same $K$-dimension.

\secstart{} Over a complete local noetherian ring $R$ with perfect residue
field $k$ we have the following alternative description of a special
formal $p$-divisible $O_K$-module due to Zink [Zi2]. For this we need a
more general definition of the Witt ring.

Let $R$ be an arbitrary commutative ring with 1. The Witt ring $W(R)$ is
characterized by the following properties:
\indention{(b)}
\litem{(a)} As a set it is given by $R^{\N_0}$, i.e.\ elements of $W(R)$
can be written as infinite tuples $(x_0,x_1,\ldots,x_i,\ldots)$.
\litem{(b)} If we associate to each ring $R$ the ring $W(R)$ and to each
homomorphism of rings $\alpha\colon R \arr R'$ the map
$$W(\alpha)\colon
(x_0,x_1,\ldots) \asr (\alpha(x_0),\alpha(x_1),\ldots),$$
then we obtain a functor from the category of rings into the category of rings.
\litem{(c)} For all integers $n \geq 0$ the so called Witt polynomials
$$\eqalign{w_n\colon W(R) &\arr R \cr
(x_0,x_1,\ldots) &\asr x_0^{p^n} + px_1^{p^{n-1}} + \ldots + p^nx_n\cr}$$
are ring homomorphisms.

For the existence of such a ring see e.g.\ [BouAC] chap.\ IX, \pz 1. If we
endow the product $R^{\N_0}$ with the usual ring structure the map
$$\xline \asr (w_0(\xline),w_1(\xline),\ldots)$$
defines a homomorphism of rings
$$W_*\colon W(R) \arr R^{\N_0}.$$

The ring $W(R)$ is endowed with two operators $\tau$ and $\sigma$ which are
characterized by the property that they are functorial in $R$ and that they
make the following diagrams commutative
$$\matrix{W(R) & \arvarover(35){\tau} & W(R) \cr
\addleft{W_*} && \addright{W_*} \cr
R^{\N_0} & \arvarover(35){\xline \asr (0,px_0,px_1,\ldots)} & R^{\N_0},
\cr}$$
$$\matrix{W(R) & \arvarover(35){\sigma} & W(R) \cr
\addleft{W_*} && \addright{W_*} \cr
R^{\N_0} & \arvarover(35){\xline \asr (x_1,x_2,\ldots)} &
R^{\N_0}. \cr}$$
The operator $\tau$ can be written explicitly by $\tau(x_0,x_1,\ldots) =
(0,x_0,x_1,\ldots)$ and it is called {\it Verschiebung of $W(R)$}. It is an
endomorphism of the additive group of $W(R)$. If $R$ is of characteristic
$p$ (i.e.\ $pR = 0$), $\sigma$ can be described as $(x_0,x_1,\ldots) \asr
(x_0^p,x_1^p,\ldots)$. For an arbitrary ring, $\sigma$ is a ring
endomorphism and it is called {\it Frobenius of $W(R)$}.

There are the following relations for $\sigma$ and $\tau$:
\indention{(iii)}
\litem{(i)} $\sigma \circ \tau = p\cdot \id_{W(R)}$,
\litem{(ii)} $\tau(x\sigma(y)) = \tau(x)y$ for $x,y \in W(R)$,
\litem{(iii)} $\tau(x)\tau(y) = p\tau(xy)$ for $x,y \in W(R)$,
\litem{(iv)} $\tau(\sigma(x)) = \tau(1)x$ for $x \in W(R)$, and we have
$\tau(1) =p$ if $R$ is of characteristic $p$.

We have a surjective homomorphism of rings
$$w_0\colon W(R) \arr R, \qquad (x_0,x_1,\ldots) \asr x_0,$$
and we denote its kernel $\tau(W(R))$ by $I_R$. We have $I_R^n =
\tau^n(W(R))$ and $W(R)$ is complete with respect to the $I_R$-adic
topology.

If $R$ is a local ring with maximal ideal $\mfr$, $W(R)$ is local as well
with maximal ideal $\set{(x_0,x_1,\ldots) \in W(R)}{$x_0 \in \mfr$}$.

\secstart{}\label{\equivdisplay} Now we can use Zink's theory of displays
to give a description of special formal $p$-divisible groups in terms of
semi-linear algebra. Let $R$ be a complete local noetherian $O_K$-algebra
with perfect residue field $k$. We extend $\sigma$ and $\tau$ to $W(R)
\tensor{\Zp} O_K$ in an $O_K$-linear way. Then we get using [Zi2]:

\claim Proposition: The category of special
formal $p$-divisible $O_K$-modules of height $h$ over $R$ is equivalent to
the category of tuples $(P,Q,F,V^{-1})$ where
\bulletlist
\bulletitem $P$ is a finitely generated $W(R) \tensor{\Zp} O_K$-module
which is free of rank $h$ over $W(R)$,
\bulletitem $Q \subset P$ is a $W(R) \tensor{\Zp} O_K$-submodule which contains
$I_RP$, and the quotient $P/Q$ is a direct summand of the $R$-module
$P/I_RP$ such that the induced action of $R \tensor{\Zp} O_K$ on $P/Q$
factorizes through the multiplication $R \otimes O_K \arr R$,
\bulletitem $F\colon P \arr P$ is a $\sigma$-linear map,
\bulletitem $V^{-1}\colon Q \arr P$ is a $\sigma$-linear map whose image
generates $P$ as a $W(R)$-module,

{\sl satisfying the following two conditions:
\indention{{\rm (b)}}
\litem{{\rm (a)}} For all $m \in P$ and $x \in W(R)$ we have the relation
$$V^{-1}(\tau(x)m) = xF(m).$$
\litem{{\rm (b)}} The unique $W(R) \tensor{\Zp} O_K$-linear map
$$V^{\#}\colon P \arr W(R) \tensor{\sigma,W(R)} P$$
satisfying the equations
$$\eqalign{V^{\#}(xFm) &= p\cdot x \otimes m \cr
V^{\#}(xV^{-1}n) &= x \otimes n\cr}$$
for $x \in W(R)$, $m \in P$ and $n \in Q$ is topologically nilpotent, i.e.\
the homomorphism $V^{N\#}\colon P \arr W(R) \tensor{\sigma^N,W(R)} P$ is
zero modulo $I_R + pW(R)$ for $N$ sufficiently large.}

\secstart{}\label{\definefullset} To define the notion of a Drinfeld level
structure we need the following definition: Let $R$ be a ring and let $X =
\Spec(A)$ where $A$ is finite locally free over $R$ of rank $N \geq 1$. For
any $R$-algebra $R'$ we denote by $X(R')$ the
set of $R$-algebra homomorphisms $A \arr R'$ (or equivalently of all
$R'$-algebra homomorphisms $A \tensor{R} R' \arr R'$). The multiplication
with an element $f \in A \tensor{R} R'$ defines an $R'$-linear endomorphism
of $A \tensor{R} R'$. As $A$ is finite locally free we can speak of the
determinant of this endomorphism which is an element ${\rm Norm}(f)$ in $R'$.

We call a finite family of elements $\varphi_1,\ldots,\varphi_N \in X(R')$
a {\it full set of sections of $X$ over $R'$} if we have for every
$R'$-algebra $T$ and for all $f \in A \tensor{R} T$ an equality in $T$
$${\rm Norm}(f) = \prod_{i=1}^N\varphi_i(f).$$

\secstart{} Let $R$ be an $O_K$-algebra and let $G$ be a special $p$-divisible
$O_K$-module over $R$. We assume that its $O_K$-height $h$ is constant on
$\Spec(R)$, e.g.\ if $R$ is a local ring (the only case which will be used
in the sequel). The $O_K$-action on $G$ defines for every integer $m \geq
1$ the multiplication with $\pi_K^m$
$$[\pi_K^m]\colon G \arr G.$$
This is an endomorphism of $p$-divisible groups whose kernel
is a finite locally free group scheme $G[\pi_K^m]$ over $\Spec(R)$ of rank
$q^{mh}$.

Let $R'$ be an $R$-algebra. A {\it Drinfeld $\pfr_K^m$-structure on $G$
over $R'$} is a homomorphism of $O_K$-modules
$$\alpha\colon (\pfr^{-m}/O_K)^h \arr G[\pi_K^m](R')$$
such that the finite set of $\alpha(x)$ for $x \in (\pfr^{-m}_K/O_K)^h$
forms a full set of sections.

\secstart{} It follows from the definition \definefullset\ that if
$\alpha\colon (\pfr^{-m}/O_K)^h \arr G[\pi_K^m](R')$ is a Drinfeld
$\pfr^m_K$-structure over $R'$ then for any $R'$-algebra $T$ the composition
$$\alpha_T\colon (\pfr^{-m}/O_K)^h \arrover{\alpha} G[\pi_K^m](R') \arr
G[\pi_K^m](T),$$
where the second arrow is the canonical one induced by functoriality from
$R' \arr T$, is again a Drinfeld $\pfr_K^m$-structure.

\secstart{}\label{\universalDrinfeld} Being a Drinfeld $\pfr_K^m$-structure
is obviously a closed property. More precisely: Let $\alpha\colon
(\pfr^{-m}/O_K)^h \arr G[\pi_K^m](R')$ be a homomorphism of abelian
groups. Then there exists a (necessarily unique) finitely generated ideal
$\afr \subset R'$ such that a homomorphism of $O_K$-algebras $R' \arr T$
factorizes over $R'/\afr$ if and only if the composition $\alpha_T$ of
$\alpha$ with the canonical homomorphism $G[\pi_K^m](R') \arr
G[\pi_K^m](T)$ is a Drinfeld $\pfr_K^m$-structure over $T$.

It follows that for every special formal $p$-divisible $O_K$-module
$(G,\iota)$ over
some $O_K$-algebra $R$ the functor on $R$-algebras which associates to
each $R$-algebra $R'$ the set of Drinfeld $\pfr_K^m$-structures on
$(G,\iota)_{R'}$ is representable by an $R$-algebra ${\rm
  DL}_m(G,\iota)$ which is of finite presentation as $R$-module. Obviously
${\rm DL}_0(G,\iota) = R$.

\secstart{}\label{\lowerlevel} Let $\alpha\colon (\pfr^{-m}/O_K)^h \arr
G[\pi_K^m](R')$ be a Drinfeld $\pfr_K^m$-structure over $R'$. As $\alpha$
is $O_K$-linear, it induces for all $m' \leq m$ a homomorphism
$$\alpha[\pi_K^{m'}]\colon (\pfr^{-m'}/O_K)^h \arr G[\pi_K^{m'}](R').$$

\claim Proposition: This is a Drinfeld $\pfr_K^{m'}$-structure.

\medskip

For the proof of this non-trivial fact we refer to [HT] 3.2 (the hypothesis
in loc.\ cit.\ that $\Spec(R')$ is noetherian with a dense set of points
with residue field algebraic over $\kappa$ is superfluous as we can always
reduce to this case by [EGA] IV, \pz 8 and \universalDrinfeld). If $R'$ is
a complete local noetherian ring with perfect residue class field (this is
the only case which we will use in the sequel) the proposition follows from
the fact that we can represent $(G,\iota)$ by a formal group law and that in
this case a Drinfeld level structure as defined above is the same as a
Drinfeld level structure in the sense of [Dr].

\secstart{}\label{\levelflat} Let $(G,\iota)$ be a special formal
$p$-divisible $O_K$-module over an $O_K$-algebra $R$. By \lowerlevel\ we
get for non-negative integers
$m \geq m'$ canonical homomorphisms of $R$-algebras
$${\rm DL}_{m'}(G,\iota) \arr {\rm DL}_m(G,\iota).$$
It follows from [Dr] 4.3 that these homomorphisms make ${\rm
  DL}_m(G,\iota)$ into a finite locally free module over ${\rm
  DL}_{m'}(G,\iota)$.

\secstart{Example}:\label{\trivDrinfeld} If $R$ is an $O_K$-algebra of
characteristic $p$ and if $G$ is a special formal $p$-divisible $O_K$-module
of $O_K$-height $h$ and of dimension 1, then the trivial homomorphism
$$\alpha^{\rm triv}\colon (\pfr_K^{-m}/O_K)^h \arr G[\pfr_K^m], \qquad x
\asr 0$$
is a Drinfeld $\pfr_K^m$-structure. If $R$ is reduced, this is the only one.


\paragraph{Deformation of p-divisible O-modules}

\secstart{} In this paragraph we fix an algebraically closed field $k$ of
characteristic $p$ together with a homomorphism $O_K \arr
k$. Further we fix integers $h \geq 1$ and $m \geq 0$. By
\uniqueformmodule\ and by \trivDrinfeld\ there exists up to isomorphism
only one special formal $p$-divisible $O_K$-module $\Sigma_h$ of height $h$
and dimension 1 with Drinfeld $\pfr_K^m$-structure $\alpha^{\rm triv}$ over
$k$. We denote the pair $(\Sigma_h,\alpha^{\rm triv})$ by $\Sigma_{h,m}$.

Let $\Cscr$ be the category of pairs $(R,s)$ where $R$ is a complete
local noetherian $O_K$-algebra and where $s$ is an isomorphism of the
residue class field of $R$ with $k$. The morphisms in $\Cscr$ are
local homomorphisms of $O_K$-algebras inducing the identity on $k$.

\secstart{Definition}:\label{\definedeformation} Let $(R,s) \in \Cscr$ be a
complete local noetherian $O_K$-algebra. A triple $(G,\alpha,\varphi)$
consisting of a formal special $p$-divisible $O_K$-module $G$ over $R$, of
a Drinfeld $\pfr_K^m$-structure $\alpha$ of $G$ over $R$ and of an isomorphism
$$\varphi\colon \Sigma_{h,m} \arriso (G \tensor{R} k, \alpha_k)$$
is called a {\it deformation of $\Sigma_{h,m}$ over $R$}.

A triple $(R_{h,m}, \Sgtilde_{h,m}, \varphi)$ consisting of a complete local
noetherian ring $R_{h,m}$ with residue field $k$ and of a deformation
$(\Sgtilde_{h,m},\varphi)$ of $\Sigma_{h,m}$ is called {\it universal
  deformation of $\Sigma_{h,m}$} if for every deformation
$(G,\alpha,\varphi)$ over some $R \in \Cscr$ there exists a unique morphism
$R_{h,m} \arr R$ in $\Cscr$ such that $(\Sgtilde_{h,m},\varphi)_R$ is
isomorphic to $(G,\alpha,\varphi)$.

A universal deformation is unique up to unique isomorphism if it exists.

\secstart{Proposition}:\label{\univdefexist} {\sl We keep the notations of
\definedeformation.
\assertionlist
\assertionitem A universal deformation $(R_{h,m}, \Sgtilde_{h,m}, \varphi)$
of $\Sigma_{h,m}$ exists.
\assertionitem For $m = 0$ the complete local noetherian ring $R_{h,0}$ is
isomorphic to the power series ring $W_K(k)\dlbrack
t_1,\ldots,t_{h-1}\drbrack$.
\assertionitem For $m \geq m'$ the canonical homomorphisms $R_{h,m'} \arr
R_{h,m}$ are finite flat. The rank of the free $R_{h,0}$-module $R_{h,m}$
is $\#GL_h(O_K/\pfr_k^m)$.
\assertionitem The ring $R_{h,m}$ is regular for all $m \geq 0$.}

\proof: Assertion (1) follows from a criterion of Schlessinger [Sch] using
rigidity for $p$-divisible groups (e.g. [Zi1]) and the fact that the canonical
functor from the category of special $p$-divisible $O_K$-modules over ${\rm
  Spf}(R_{h,m})$ to the category of special $p$-divisible $O_K$-modules
over $\Spec(R_{h,m})$ is an equivalence of categories (cf. [Me] II, 4). The
second assertion follows easily from general deformation theory of
$p$-divisible groups (for an explicit description of the universal
deformation and a proof purely in terms of linear algebra one can use [Zi2]
and \equivdisplay). Finally, (3) and (4) are more involved (see [Dr] \pz 4,
note that (3) is essentially equivalent to \levelflat).

\secstart{}\label{\firstaction} Let $D_{1/h}$ be ``the'' skew field with
center $K$ and invariant $1/h$. The ring $R_{h,m}$ has a continuous action
of the ring of units $O\cross_{D_{1/h}}$ of the integral closure
$O_{D_{1/h}}$ of $O_K$ in $D_{1/h}$: Let $\Sgtilde_{g,m} =
(G,\alpha,\varphi)$ be the universal special formal $p$-divisible
$O_K$-module with Drinfeld $\pfr_K^m$-structure over $R_{h,m}$. For $\delta
\in O\cross_{D_{1/h}}$ the composition
$$\Sigma_{h,m} \arrover{\delta} \Sigma_{h,m} \arrover{\varphi} (G,\alpha)
\tensor{R_{h,m}} k$$
is again an isomorphism if we consider $\delta$ as an automorphism of
$\Sigma_{h,m}$ (which is the same as an automorphism of $\Sigma_{h,0}$) by
\endoalgebra. Therefore $(G,\alpha,\varphi \circ \delta)$ is a deformation
of $\Sigma_{h,m}$ over $R_{h,m}$ and by the definition of a universal
deformation this defines a continuous automorphism $\delta\colon R_{h,m}
\arr R_{h,m}$.

\secstart{}\label{\secondaction} Similarly as in \firstaction\ we also get a
continuous action of $GL_h(O_K/\pfr_K^m)$ on $R_{h,m}$: Again let
$\Sgtilde_{g,m} = (G,\alpha,\varphi)$ be the universal special formal
$p$-divisible $O_K$-module with Drinfeld $\pfr_K^m$-structure over
$R_{h,m}$. For $\gamma \in \GL_h(O_K/\pfr_K^m)$, $\alpha \circ \gamma$ is
again a Drinfeld $\pfr_K^m$-structure, hence $(G,\alpha \circ \gamma,
\varphi)$ is a deformation of $\Sigma_{h,m}$ and defines a continuous
homomorphism
$$\gamma\colon R_{h,m} \arr R_{h,m}.$$

\secstart{}\label{\action} By combining \firstaction\ and \secondaction\ we
get a continuous left action of
$$\GL_h(O_K) \times O\cross_{D_{1/h}} \aerr GL_h(O_K/\pfr_K^m) \times
O\cross_{D_{1/h}}$$
on $R_{h,m}$. Now we have the following lemma ([HT] p.\ 52)

\claim Lemma: This action can be extended to a continuous left action of
$GL_h(K) \times D\cross_{1/h}$ on the direct system of the $R_{h,m}$ such
that for $m_2 >> m_1$ and for $(\gamma,\delta) \in \GL_h(K) \times
D\cross_{1/h}$ the diagram
$$\matrix{R_{h,m_1} & \arvarover(35){(\gamma,\delta)} & R_{h,m_2} \cr
\auu && \auu \cr
W(k) & \arvarover(35){\sigma_K^{v_K(\det(\gamma)) - v_K({\rm
      Nrd}(\delta))}} & W(k) \cr}$$
commutes.


\paragraph{Vanishing cycles}

\secstart{}\label{\notatvan} Let $W$ be a complete discrete valuation ring
with maximal ideal $(\pi)$, residue field $k$ and field of fractions
$L$. Assume that $k$ is algebraically closed (or more generally separably
closed). The example we will use later on is the ring $W = W_K(k)$ for an
algebraically closed field $k$ of characteristic $p$. Set $\eta =
\Spec(L)$, $\hgbar = \Spec(\Lbar)$ and $s = \Spec(k)$.

We will first define vanishing and nearby cycles for an algebraic
situation (cf.\ [SGA 7] exp.\ I, XIII). Then we will generalize to the
situation of formal schemes.

\secstart{} Let $f\colon X \arr \Spec(W)$ be a scheme of finite type over
$W$ and define $X_{\hgbar}$ and $X_s$ by cartesian diagrams
$$\matrix{X_s & \arrover{i} & X & \allover{{\bar j}} & X_{\hgbar} \cr
\addleft{f_s} && \addleft{f} && \addright{f_{\hgbar}} \cr
s & \arr & \Spec(W) & \al \eta \al & \hgbar.\cr}$$
The formalism of vanishing cycles is used to relate the cohomology of $X_s$
and of $X_{\hgbar}$ together the action of the inertia group on the
cohomology of $X_{\hgbar}$.

Fix a prime $\ell$ different from the characteristic $p$ of $k$ and let
$\Lambda$ be a finite abelian group which is annihilated by a power of
$\ell$. For all
integers $n \geq 0$ the sheaf
$$\Psi^n(\Lambda) = i^*R^n{\bar j}_*\Lambda$$
is called the {\it sheaf of vanishing cycles} of $X$ over $W$. Via
functoriality it carries an action of $\Gal(\Lbar/L)$. Note that by
hypothesis $\Gal(\Lbar/L)$ equals the inertia group of $L$.

\secstart{} If $f\colon X \arr \Spec(W)$ is proper, the functor $i^*$
induces an isomorphism (proper base change)
$$i^*\colon H^p(X,R^q{\bar j}\Lambda) \arriso
H^p(X_s,i^*R^q{\bar j}_*\Lambda)$$
and the Leray spectral sequence for ${\bar j}$ can be written as
$$H^p(X_s,\Psi^n(\Lambda)) \implies H^{p+q}(X_{\hgbar},\Lambda).$$
This spectral sequence is $\Gal(\Lbar/L)$-equivariant.
This explains, why the vanishing cycles ``measure'' the difference of the
cohomology of the special and the generic fibre.

\secstart{}\label{\vanishsmooth} Let $\Lambda$ be as above. If $X$ is
proper and smooth over $W$, there is no difference between the cohomology
of the generic and the special fibre. More precisely, we have ([SGA 4] XV,
2.1):

\claim Proposition: Let $f\colon X \arr W$ be a smooth and proper
morphism. Then we have for all $n \geq 0$ a canonical isomorphism
$$H^n(X_{\hgbar},\Lambda) \arriso H^n(X_s,\Lambda).$$
This isomorphism is $\Gal(\Lbar/L)$-equivariant where the action on the
right hand side is trivial. Further $\psi^n(\Lambda) = 0$ for $n \geq 1$
and $\psi^0 = \Lambda$.

\secstart{} Now we define vanishing cycles for formal schemes: We keep the
notations of \notatvan. We call a topological $W$-algebra $A$ {\it
  special}, if there exists an ideal $\afr \subset A$ (called an {\it
  ideal of definition of $A$}\/) such that $A$ is complete with respect to
the $\afr$-adic topology and such that $A/\afr^n$ is a finitely
generated $W$-algebra for all $n \geq 1$ (in fact the same condition for $n
= 2$ is sufficient ([Ber2] 1.2)). Equivalently, $A$ is topologically
$W$-isomorphic to a quotient of the topological $W$-algebra
$$W\{T_1,\ldots,T_m\}\dlbrack S_1,\ldots,S_n\drbrack = W\dlbrack
S_1,\ldots,S_n\drbrack\{T_1,\ldots,T_m\},$$
in particular $A$ is again noetherian. Here if $R$ is a topological
ring, $R\{T_1,\ldots,T_m\}$ denotes the subring of power series
$$\sum_{\nline \in \N_0^n}c_{\nline}T^{\nline}$$
in $R\dlbrack T_1,\ldots,T_M \drbrack$ such that for every neighborhood $V$
of $0$ in $R$ there is only a finite number of coefficients $c_{\nline}$
not belonging to $V$. If $R$ is complete Hausdorff with respect to the
$\afr$-adic topology for an ideal $\afr$, we have a canonical isomorphism
$$R\{T_1,\ldots,T_m\} \arriso \limproj_{n} (R/\afr^n)[T_1,\ldots,T_m].$$

Note that all the rings $R_{h,m}$ defined in \univdefexist\ are special
$W_K(k)$-algebras.

\secstart{} Let $A$ be a special $W$-algebra and let $\Xscr =
\Spf(A)$. Denote by $\Xscr^{\rm rig}$ its associated rigid space over
$L$. It can be constructed as follows: For
$$A = W\{T_1,\ldots,T_m\}\dlbrack S_1,\ldots,S_n\drbrack$$
we have $\Xscr^{\rm rig} = E^m \times D^n$ where
$E^m$ and $D^n$ are the closed resp.\ open polydiscs of radius $1$ with
center at zero in $L^m$ resp.\ in $L^n$. If now $A$ is some quotient of the
special $W$-algebra $A' = W\{T_1,\ldots,T_m\}\dlbrack
S_1,\ldots,S_n\drbrack$ with kernel $\afr' \subset A'$, $\Xscr^{\rm rig}$
is the closed rigid analytic subspace of $\Spf(A')^{\rm rig}$ defined by
the sheaf of ideals $\afr'\Oscr_{\Spf(A')^{\rm rig}}$. This defines a
functor from the category of special $W$-algebras to the category of
rigid-analytic spaces over $L$ (see [Ber2] \pz 1 for the precise definition).

\secstart{} Let $A$ be a special $W$-algebra and let $\afr \subset A$ be
its largest ideal of definition. A topological $A$-algebra $B$ which is
special as $W$-algebra (or shorter a special $A$-algebra - an abuse of
language which
is justified by [Ber2] 1.1) is called {\it \'etale over $A$} if $B$ is
topologically finitely generated as $A$-algebra (which is equivalent to the
fact that for every ideal of definition $\afr'$ of $A$, $\afr'B$ is an
ideal of definition of $B$) and if the morphism of
commutative rings $A/\afr \arr B/\afr B$ is \'etale in the usual sense.

The assignment $B \asr B/\afr B$ defines a functor
from the category of \'etale special $A$-algebras to the category of \'etale
$A/\afr$-algebras which is an equivalence of categories and hence we get an
equivalence of \'etale sites
$$(A)_{\et} \sim (A/\afr)_{\et}.$$
Note that for every ideal of definition $\afr'$ of $A$ the \'etale sites
$(A/\afr)_{\et}$ and $(A/\afr')_{\et}$ coincide. We just chose
the largest ideal of definition to fix notations.

On the other hand, if we write $\Xscr = \Spf(A)$ and $\Yscr = \Spf(B)$ for
an \'etale special $A$-algebra $B$ we get a (quasi-)\'etale morphism of
rigid analytic spaces
$$\Yscr^{\rm rig} \arr \Xscr^{\rm rig}.$$

By combining this functor with a quasi-inverse of $B \asr B/\afr B$ we get
a morphism of \'etale sites
$$s\colon (\Xscr^{\rm rig} \tensor{L} \Lbar)_{\et} \arr (\Xscr^{\rm
  rig})_{\et} \arr (\Xscr)_{\et} \arriso (\Xscr_{\rm red})_{\et}$$
with $\Xscr_{\rm red} = \Spec(A/\afr)$.

Let $\Lambda$ be a finite abelian group which is annihilated by a power of
$\ell$. For $n \geq 0$ the sheaves
$$\psi^n(\Lambda) := R^ns_*\Lambda$$
are called {\it vanishing cycle sheaves}.

\secstart{}\label{\smoothaction} Let $A$ be a special $W$-algebra with
largest ideal of definition $\afr$, $\Xscr =
\Spf(A)$. Then the group $\Aut_W(\Xscr)$ of automorphisms of $\Xscr$ over
$W$ (i.e.\ of continuous $W$-algebra automorphisms $A \ar A$) acts on
$\psi^n(\Lambda)$. Further we have the following result of Berkovich
[Ber2]:

\claim Proposition: Assume that $\psi^i(\Z/\ell\Z)$ is constructible for
all $i$. Then there exists an integer $n \geq 1$ with the following
property: Every element $g \in \Aut_W(\Xscr)$ whose image in
$\Aut_W(A/\afr^n)$ is the identity acts trivially on $\psi^i(\Z/\ell^m\Z)$
for all integers $i,m \geq 0$.


\paragraph{Vanishing cycles on the universal deformation of special p-divisible
  O-modules}

\secstart{} Let $k$ be an algebraic closure of the residue field $\kappa$
of $O_K$. Then $W = W_K(k)$ is the ring of integers of $\Khat^{\rm nr}$,
the completion of the maximal unramified extension of $K$. Further denote
by $I_K$ the inertia group and by $W_K$ the Weil group of $K$.

\secstart{} Consider the system $P$ of special formal schemes
$$\ldots \arr \Spf(R_{m,h}) \arr \Spf(R_{m-1,h}) \arr \ldots \arr
\Spf(R_{0,h}).$$
By applying the functor $(\ )^{\rm rig}$ we get a system $P^{\rm rig}$ of
rigid spaces
$$\ldots \arr \Spf(R_{m,h})^{\rm rig} \arr \Spf(R_{m-1,h})^{\rm rig} \arr
\ldots \arr \Spf(R_{0,h})^{\rm rig}.$$
these systems $P$ and $P^{\rm rig}$ have an action by the group $\GL_h(O_K)
\times O\cross_{D_{1/h}}$ \firstaction\ and \secondaction.
Denote by $\Psi^i_m(\Lambda)$ the vanishing cycle sheaf for $\Spf(R_{h,m})$
with coefficients in some finite abelian $\ell$-primary group $\Lambda$ and
set
$$\Psi^i_m = (\limproj_n \Psi^i_m(\Z/\ell^n\Z)) \tensor{\Z_{\ell}}
\Qdbar_{\ell}.$$
Note that we have
$$\limproj_n \Psi^i_m(\Z/\ell^n\Z) = H^i((\Spf R_m)^{\rm rig}
\tensor{\hat K^{\rm nr}} {\hat \Kbar}, \Z_{\ell}),$$
in particular these $\Z_{\ell}$-modules carry an action by $I_K =
\Gal({\hat \Kbar}/{\hat K^{\rm nr}})$.
Further write
$$\Psi^i = \limind_m \Psi^i_m.$$
Via our chosen identification $\Qdbar_{\ell} \cong \C$ we can consider
$\Psi^i_m$ and $\Psi^i$ as $\C$-vector spaces which carry an action of
$$\GL_h(O_K) \times O\cross_{D_{1/h}} \times I_K.$$

\secstart{Lemma}: {\sl We have the following properties of the $(\GL_h(O_K)
\times O\cross_{D_{1/h}} \times I_K)$-modules $\Psi^i_m$ and $\Psi^i$.
\assertionlist
\assertionitem The $\Psi^i_m$ are finite-dimensional $\C$-vector spaces.
\assertionitem We have $\Psi^i_m = \Psi^i = 0$ for all $m \geq 0$ and for
all $i > h-1$.
\assertionitem The action of $GL_h(O_K)$ on $\Psi^i$ is admissible.
\assertionitem The action of $O\cross_{D_{1/h}}$ on $\Psi^i$ is smooth.
\assertionitem The action of $I_K$ on $\Psi^i$ is continuous.}

\proof: For the proof we refer to [HT] 3.6. We only remark that (3) -- (5)
follow from general results of Berkovich [Ber2] and [Ber3] if we know
(1). To show (1) one uses the fact that the system of formal schemes $P$
comes from an inverse system of proper schemes of finite type over $W$
(cf.\ the introduction) and a comparison theorem of Berkovich which relates
the vanishing cycles of a scheme of finite type over $W$ with the
vanishing cycle sheaves for the associated formal scheme.

\secstart{} Let $A_h$ be the group of elements $(\gamma,\delta,\sigma) \in
  \GL_h(K) \times D\cross_{1/h} \times W_K$ such that
$$v_K(\det(\gamma)) = v_K({\rm Nrd}(\delta)) + v_K(\Art_K^{-1}(\sigma)).$$
The action of $\GL_h(K) \times D\cross_{1/h}$ on the system $(R_{h,m})_m$
\action\ gives rise to an action of $A_K$ on $\Psi^i$.

Moreover, if $(\rho,V_{\rho})$ is an irreducible admissible representation of
$D\cross_{1/h}$ over $\C$ (and hence necessarily finite-dimensional
\compactfin) then we set
$$\Psi^i(\rho) = \Hom_{O\cross_{1/h}}(\rho, \Psi^i).$$
This becomes naturally an admissible $\GL_h(K) \times W_K$-module if we
define for $\phi \in \Psi^i(\rho)$ and for $x \in V_{\rho}$
$$((\gamma,\sigma)\phi)(x) =
(\gamma,\delta,\sigma)\phi(\rho(\delta)^{-1}x)$$
where $\delta \in D\cross_{1/h}$ is some element with $v_K({\rm
  Nrd}(\delta)) = v_K(\det(\gamma)) - v_K(\Art_K^{-1}(\sigma))$.


\bibliography

\indention{[BouAC]\ }
\litem{[AT]} E. Artin, J. Tate: {\it Class Field Theory}, Benjamin (1967).
\litem{[Au1]} A.-M. Aubert: {\it Dualit\'e dans le groupe de
Grothendieck de la cat\'egorie des repr\'esentations lisses de longueur
finie d'un groupe r\'eductif $p$-adique}, Trans. Amer. Math. Soc. {\bf
347}, no. 6, 2179--2189 (1995).
\litem{[Au2]} A.-M. Aubert: {\it Erratum: Duality in the Grothendieck
group of the category of finite-length smooth representations of a $p$-adic
reductive group}, Trans. Amer. Math. Soc. {\bf 348}, no. 11, 4687--4690
(1996).
\litem{[AW]} M.F. Atiyah, C.T.C. Wall: {\it Cohomology of Groups}, in J.W.S
Cassels, A. Fr\"ohlich (Ed.): Algebraic Number Theory, Acadamic Press, (1967).
\litem{[Ber1]} V.G. Berkovich: {\it Vanishing cycles for formal schemes},
Inv. Math. {\bf 115}, 539--571 (1994).
\litem{[Ber2]} V.G. Berkovich: {\it Vanishing cycles for formal schemes II},
Inv. Math. {\bf 125}, 367--390 (1996).
\litem{[Ber3]} V.G. Berkovich: {\it Vanishing cycles}, Appendix II to:
M. Harris, R. Taylor: ``On the geometry and cohomology of some simple
Shimura varieties'', preprint Harvard, 1999.
\litem{[BHK]} C.J. Bushnell, G. Henniart, P.C. Kutzko: {\it Correspondance
  de Langlands locale pour $\GL_n$ et conducteurs de paires}, Ann. Sci. ENS
$4^e$ s\'erie {\bf 31}, pp. 537--560 (1998).
\litem{[BouA]} N. Bourbaki: {\it Algebra}, chap. IV-VII, Springer (1988);
chap. VIII, Hermann (1958).
\litem{[BouAC]} N. Bourbaki: {\it Alg\`ebre commutative, chapitres 8 et 9},
Masson, Paris, 1983.
\litem{[BouL]} N. Bourbaki: {\it Groupes et Alg\`ebres de Lie}, chap. VII
et VIII, Masson, 1990.
\litem{[Bu]} C. Bushnell: {\it Gauss sums and local constants for GL(N)},
$L$-functions and Arithmetic (J. Coates, M.J. Taylor, eds.), London
Math. Soc. Lecture Notes {\bf 153}, Cambridge University press,
61--73 (1991).
\litem{[BW]} A. Borel, N. Wallach: {\it Continuous cohomology, discrete
subgroups, and representations of reductive groups}, Ann. of
Math. Stud. {\bf 94}, Princeton University Press (1980).
\litem{[BZ1]} I.N. Bernstein, A.V. Zelevinsky: {\it Representations of the
  group $\GL(n,F)$ where $F$ is a non-archimedean local field},
Russ. Math. Surveys {\bf 31}:3, 1--68 (1976).
\litem{[BZ2]} I.N. Bernstein, A.V. Zelevinsky: {\it Induced representations
  of reductive $\pfr$-adic groups Id}, Ann. Sci. ENS $4^e$ s\'erie {\bf 10},
441--472 (1977).
\litem{[Ca]} P. Cartier: {\it Representations of $\pfr$-adic groups: a
  survey}, Proc. of Symp. in Pure Math. {\bf 33}, part 1, 111--155
(1979).
\litem{[Cas1]} W. Casselman: {\it Introduction to the theory of admissible
  representations of $p$-adic reductive groups}, draft 1995
(www.math.ubc.ca/people/faculty/cass/\break research.html).
\litem{[Cas2]} W. Casselman: {\it The unramified principal series of
$\pfr$-adic groups I, The spherical function}, Compositio Math. {\bf
40}, no. 3, 387--406 (1980). 
\litem{[CF]} J.W.S Cassels, A. Fr\"ohlich: {\it Algebraic Number Theory},
Acad. Press (1967).
\litem{[De1]} P. Deligne: {\it Formes modulaires et repr\'esentations de
${\rm GL}(2)$}, in ``Modular functions of one variable II'' (Proc.
Internat. Summer School, Univ. Antwerp, Antwerp, 1972), Lecture Notes in
Mathematics {\bf 349}, Springer, 55--105 (1973).
\litem{[De2]} P. Deligne: {\it Les constantes des \'equations fonctionelles
  des fonctions $L$}, in ``Modular functions of one variable II'' (Proc.
Internat. Summer School, Univ. Antwerp, Antwerp, 1972), Lecture Notes in
Mathematics {\bf 349}, Springer, 501--596 (1973).
\litem{[Dr1]} V.G. Drinfeld: {\it Elliptic modules}, Mat. USSR Sbornik {\bf
  23}, pp. 561--592 (1974).
\litem{[Dr2]} V.G. Drinfeld: {\it Elliptic modules II}, Mat. USSR Sbornik {\bf
  31}, pp. 159--170 (1977).
\litem{[EGA]} A. Grothendieck, J. Dieudonn\'e: El\'ements de g\'eom\'etrie
alg\'ebrique, I Grund\-lehren der Mathematik {\bf 166} (1971) Springer, II-IV
Publ. Math. IHES {\bf 8} (1961) {\bf 11} (1961) {\bf 17} (1963) {\bf 20}
(1964) {\bf 24} (1965) {\bf 28} (1966) {\bf 32} (1967).
\litem{[Fl]} D. Flath: {\it Decomposition of representations into tensor
  products}, Proc. Symp. Pure Math. {\bf 33}, part 1, pp. 179--183 (1979).
\litem{[GJ]} R. Godement, H. Jacquet: {\it Zeta functions of simple
  algebras}, Lecture Notes in Mathematics {\bf 260}, Springer, 1972.
\litem{[Has]} H. Hasse: {\it Die Normenresttheorie relativ-abelscher
  Zahlk\"orper als Klassen\-k\"orpertheorie im Kleinen}, J. reine
angew. Math. {\bf 162}, pp. 145--154 (1930).
\litem{[HC]} Harish-Chandra: {\it Admissible invariant distributions on
reductive $p$-adic groups}, Notes by S. DeBacker and P.J. Sally, Jr.,
University Lecture Series {\bf 16}, AMS (1999).
\litem{[He1]} G. Henniart: {\it La conjecture de Langlands locale pour
  $\GL(3)$}, Mem. Soc. Math. France {\bf 11-12}, pp. 1--186 (1984).
\litem{[He2]} G. Henniart: {\it La conjecture de Langlands num\'erique for
  $\GL(n)$}, Ann. Sci. ENS $4^e$ s\'erie {\bf 21}, pp. 497--544 (1988).
\litem{[He3]} G. Henniart: {\it Caract\'erisation de la correspondance de
  Langlands locale par les facteurs $\eps$ de pairs}, Inv.
Math. {\bf 113}, pp. 339--350 (1993).
\litem{[He4]} G. Henniart: {\it Une preuve simple des conjectures de
  Langlands pour $GL(n)$ sur un corps $p$-adique}, Inv. Math. {\bf 139},
pp. 439--455 (2000).
\litem{[HT]} M. Harris, R. Taylor: {\it On the geometry and cohomology of
  some simple Shimura varieties}, preprint Harvard (1999)
(www.math.harvard.edu/\~\ rtaylor).
\litem{[IM]} N. Iwahori, H. Matsumoto: {\it On some Bruhat decompositions
and the structure of the Hecke rings of $p$-adic Chevalley groups},
Publ. Math. IHES {\bf 25}, 5--48 (1965).
\litem{[Ja]} H. Jacquet: {\it Generic Representations}, in
``Non-Commutative Harmonic Analysis'', Lecture Notes in Mathematics {\bf 587},
91--101 (1977).
\litem{[JL]} H. Jacquet, R.P. Langlands: {\it Automorphic forms on
  $GL(2)$}, Lecture Notes in Mathematics {\bf 114}, Springer (1970).
\litem{[JPPS1]} H. Jacquet, I.I. Piatetski-Shapiro, J. Shalika: {\it
  Conducteur des repr\'esen\-tations du groupe lin\'eaire}, Math. Annalen
{\bf 256}, 199--214 (1981).
\litem{[JPPS1]} H. Jacquet, I.I. Piatetski-Shapiro, J. Shalika: {\it
Rankin-Selberg convolutions}, Amer. J. Math. {\bf 105}, pp. 367--483
(1983).
\litem{[Kn]} A.W. Knapp: {\it Local Langlands correspondence: The
  archimedean case}, Proc. Symp. Pure Math. {\bf 55}, part 2, pp. 393--410
(1994).
\litem{[Ko]} B. Kostant: {\it The principal three-dimensional subgroups and
the Betti numbers of a complex simple Lie group}, Amer. J. Math. {\bf 81},
973--1032 (1959).
\litem{[Kud]} S. Kudla: {\it The Local Langlands Correspondence: The
  Non-Archimedean Case}, Proc. of Symp. Pure Math. {\bf 55}, part 2,
pp. 365--391 (1994).
\litem{[Kut]} P. Kutzko: {\it The local Langlands conjecture for $\GL(2)$ of
  a local field}, Ann. of Math. {\bf 112}, pp. 381--412 (1980).
\litem{[LRS]} G. Laumon, M. Rapoport, U. Stuhler: {\it $\Dscr$-elliptic
  sheaves and the Langlands correspondence}, Invent. Math. {\bf 113},
pp. 217--338 (1993).
\litem{[Me]} W. Messing: {\it The crystals associated to Barsotti-Tate
  groups: with an application to abelian schemes}, Lecture Notes in
Mathematics {\bf 264}, Springer (1972).
\litem{[M\oe]} C. M\oe glin: {\it Representations of $\GL_n(F)$ in the
non-archimedean case}, Proc. of Symp. Pure Math. {\bf 61},
pp. 303--319 (1997).
\litem{[Neu]} J. Neukirch: {\it Algebraic number theory}, Grundlehren der
Mathematischen Wissenschaften {\bf 322}. Springer (1999).
\litem{[PR]} V. Platonov, A. Rapinchuk: {\it Algebraic Groups and Number
Theory}, Pure and Applied Mathematics {\bf 139}, Academic Press (1994).
\litem{[Pr]} K. Procter: {\it Parabolic induction via Hecke algebras and
the Zelevinsky duality conjecture}, Proc. London Math. Soc. (3) {\bf 77},
no. 1, 79--116 (1998).
\litem{[Ro]} F. Rodier: {\it Repr\'esentations de $\GL(n,k)$ o\`u $k$ est un
  corps $p$-adique},\break S\'eminaire Bourbaki exp.\ 587 (1981-1982),
Ast\'erisque {\bf 92-93}, 201--218 (1982).
\litem{[Rog]} J. Rogawski: {\it Representations of ${\rm GL}(n)$ and
division algebras over a $p$-adic field}, Duke Math. J. {\bf 50}, no. 1,
161--196 (1983). 
\litem{[Sch]} Schlessinger: {\it Functors of Artin rings}, 
Trans. Amer. Math. Soc. {\bf 130}, 208--222 (1968).
\litem{[Se1]} J.P. Serre: {\it Local Class Field Theory}, in J.W.S Cassels,
A. Fr\"ohlich (Ed.): Algebraic Number Theory, Acadamic Press, (1967).
\litem{[Se2]} J.P. Serre: {\it Local Fields}, Graduate Texts in Mathematics
{\bf 67}, second corrected printing, Springer (1995).
\litem{[SGA 4]} A. Grothendieck, M. Artin, et. al.:
{\it S\'eminaire de G\'eometrie Alg\'ebrique du Bois-Marie, Th\'eorie des
  Topos et cohomologie \'etale des sch\'emas}\/ (1963--64), Lecture Notes
in Mathematics {\bf 269}, {\bf 270}, {\bf 305}, Springer (1972--73).
\litem{[SGA 7]} A. Grothendieck, M. Raynaud, et. al.:
{\it S\'eminaire de G\'eometrie Alg\'ebrique du Bois-Marie, Groupes de
  Monodromie en G\'eometrie Alg\'ebrique}\/ (1967--68), Lecture Notes in
Mathematics {\bf 288}, {\bf 340}, Springer (1972--73).
\litem{[Sh]} J.A. Shalika: {\it The multiplicity one theorem for $\GL(n)$},
Ann. of Math. {\bf 100}, pp. 171--193 (1974).
\litem{[Si1]} A.J. Silberger: {\it The Langlands quotient theorem for
$p$-adic groups}, Math. Ann. {\bf 236}, 95--104 (1978).
\litem{[Si2]} A.J. Silberger: {\it Introduction to harmonic analysis of
reductive $p$-adic groups}, Mathematical Notes of Princeton University
Press {\bf 23}, Princeton (1979).
\litem{[Si3]} A.J. Silberger: {\it Asymptotics and Integrability Properties
for Matrix Coefficients of Admissible Representations of reductive $p$-adic
Groups}, J. of Functional Ana\-lysis {\bf 45}, 391--402 (1982).
\litem{[Ta1]} J. Tate: {\it Fourier Analysis in Number Fields and Hecke's
Zeta-Functions}, in J.W.S Cassels, A. Fr\"ohlich (Ed.): Algebraic Number
Theory, Acadamic Press, (1967).
\litem{[Ta2]} J. Tate: {\it Number Theoretic Background}, Automorphic Forms,
Representations, and L-Functions, part 2, Proc. Symp. Pure Math. {\bf 33},
AMS, pp. 3--22 (1993).
\litem{[Tad]} M. Tadi\'c: {\it Classification of unitary representations in
irreducible representations of general linear group (non-Archimedean
case)}, Ann. Sci. \'Ecole Norm. Sup. (4) {\bf 19}, no. 3, 335--382 (1986).
\litem{[Ti]} J. Tits: {\it Reductive groups over local fields},
Proc. Symp. Pure Math. {\bf 33}, part 1, pp. 29--69 (1979).
\litem{[Tu]} J. Tunnell: {\it On the local Langlands conjecture for
  $\GL(2)$}, Invent. Math. {\bf 46}, pp. 179--200 (1978).
\litem{[Vi]} M.-F. Vign\'eras: {\it Repr\'esentations $\ell$-modulaires
d'un groupe r\'eductif $p$-adique avec $\ell \not= p$}, Progress in
Mathematics {\bf 137}, Birkh\"auser, 1996.
\litem{[Ze]} A.V. Zelevinsky: {\it Induced representations of reductive
  p-adic groups II, on irreducible representations of GL(n)}, Ann. Sci. ENS
$4^e$ s\'erie {\bf 13}, pp. 165--210, 1980.
\litem{[Zi1]} T. Zink: {\it Cartiertheorie formaler Gruppen}, Teubner-Texte
zur Mathematik {\bf 68}, 1984.
\litem{[Zi2]} T. Zink: {\it The display of a formal p-divisible group},
preprint 1998, Universit\"at Bielefeld,
http://www.mathematik.uni-bielefeld.de/\~\ zink

\bye